\documentclass[11pt,a4paper,reqno]{amsart}
\usepackage[T1]{fontenc}
\usepackage{amsmath}
\usepackage{amsthm}
\usepackage{amsfonts}
\usepackage{amssymb}
\usepackage[pdftex]{graphicx,color,epsfig}
\usepackage{amsbsy}
\usepackage{mathrsfs}
\usepackage{bbm}
\usepackage{titletoc}
\newtheorem{theorem}{Theorem}[section]

\newtheorem{lemma}[theorem]{Lemma}

\newtheorem{proposition}[theorem]{Proposition}
\newtheorem{corollary}[theorem]{Corollary}

\theoremstyle{remark}
\newtheorem{remark}[theorem]{\it \bf{Remark}\/}
\newenvironment{acknowledgement}{\noindent{\bf Acknowledgement.~}}{}
\numberwithin{equation}{section}
\catcode`@=11
\def\section{\@startsection{section}{1}%
  \z@{1.5\linespacing\@plus\linespacing}{.5\linespacing}%
  {\normalfont\bfseries\large\centering}}
\catcode`@=12
\newcommand{\be}{\begin{equation}}
\newcommand{\ee}{\end{equation}}
\newcommand{\bea}{\begin{eqnarray}}
\newcommand{\eea}{\end{eqnarray}}
\newcommand{\bee}{\begin{eqnarray*}}
\newcommand{\eee}{\end{eqnarray*}}

\def\pa{\partial}
\def\na{\nabla}
\def\Eq{\mbox{\rm Eq}}

\def\RR{\mathbb{R}}

\def\ds{\displaystyle}
\def\ni{\noindent}
\def\bs{\bigskip}
\def\ms{\medskip}
\def\ss{\smallskip}
\def\eps{\varepsilon}
\def\fref#1{{\rm (\ref{#1})}}
\def\pref#1{{\rm \ref{#1}}}

\def\calH{{\mathcal H}}

\def\calC{{\mathcal C}}

\def\calA{{\mathcal A}}
\def\calE{{\mathcal E}}
\def\calO{{\mathcal O}}

\def\calD{{\mathcal D}}

\def\calT{{\mathcal T}}

\def\X{\dot{H}^1_{r}}%
\catcode`@=11
\def\supess{\mathop{\operator@font Sup\,ess}}
\catcode`@=12
\def\un{{\mathbbmss{1}}}
\title[]{A new variational approach to the stability of gravitational systems}
\author[M. Lemou]{Mohammed Lemou}
\address{CNRS and IRMAR, Universit\'e de Rennes 1, France}
\email{mohammed.lemou@univ-rennes1.fr}
\author[F. M\'ehats]{Florian M\'ehats}
\address{IRMAR, Universit\'e de Rennes 1, France}
\email{florian.mehats@univ-rennes1.fr}
\author[P. Rapha\"el]{Pierre Rapha\"el}
\address{IMT, Universit\'e Paul Sabatier, Toulouse, France}
\email{pierre.raphael@math.univ-toulouse.fr}
\begin{document}

\begin{abstract}

We consider the three dimensional gravitational Vlasov Poisson system which describes the mechanical state of a stellar system subject to its own gravity. A well-known conjecture in astrophysics is that the steady state solutions which are nonincreasing functions of their microscopic energy are nonlinearly stable by the flow. This was proved at the linear level by several authors based on the pioneering work by Antonov in 1961. Since then, standard variational techniques based on concentration compactness methods as introduced by P.-L. Lions in 1983 have led to the nonlinear stability of subclasses of stationary solutions of ground state type.

In this paper, inspired by pioneering works from the  physics litterature \cite{lynden}, \cite{WZS}, \cite{aly}, we use the monotonicity of the Hamiltonian under generalized symmetric rearrangement transformations to prove that non increasing steady solutions are local minimizer of the Hamiltonian under equimeasurable constraints, and extract compactness from suitable minimizing sequences.  This implies the nonlinear stability of nonincreasing anisotropic steady states under radially symmetric perturbations.

\end{abstract}


\maketitle
\titlecontents{section} 
[1.5em] 
{\vspace*{0.1em}\bf} 
{\contentslabel{2.3em}} 
{\hspace*{-2.3em}} 
{\titlerule*[0.5pc]{\rm.}\contentspage\vspace*{0.1em}}

\vspace*{-2mm}



\section{Introduction and main results}



\subsection{Setting of the problem}


We consider  the three dimensional gravitational Vlasov-Poisson system
\be
\label{vp}
\left   \{ \begin{array}{llll}
         \ds \pa_t f+v\cdot\nabla_x f-\nabla
         \phi_f \cdot\nabla_v f=0, \qquad (t,x,v)\in \RR_+\times \RR^3\times \RR^3\\[3mm]
          \ds f(t=0,x,v)=f_0(x,v)\geq 0,\\[3mm]
         \end{array}
\right.
\ee
where, throughout this paper,
\be
\label{phif}
 \rho_f(x)=\int_{\RR^3} f(x,v)\,dv\quad \mbox{and}\quad  \phi_f(x)=-\frac{1}{4\pi |x|}\ast \rho_f
\ee
are the density and the gravitational Poisson field associated to $f$. This nonlinear transport equation is a well known model in astrophysics for the description of the mechanical state of a stellar system subject to its own gravity and the dynamics of galaxies, see for instance \cite{binney,fridman}. \\

Unique global classical solutions for initial data $f_0\in \calC^1_c$, $f_0\geq 0$, where $\calC^1_c$ denotes the space of compactly supported and continuously differentiable functions, have been shown to exist in \cite{LP,Pf,S1}  and to propagate the corresponding regularity. Two fundamental properties of the nonlinear transport flow (\ref{vp}) are then first the preservation of the total Hamiltonian 
\be
\label{defhamiltonian}
\calH(f(t))=  
          \ds \frac{1}{2}\int_{\RR^6}
         |v|^2  f(t,x,v)dxdv-\frac{1}{2}\int_{\RR^3}|\nabla
          \phi_f(t,x)|^2dx=\calH (f(0)),
\ee
and second the preservation of all the so-called Casimir functions: $\forall G \in \calC^1([0,+\infty), \RR^+)$ such that $G(0)=0$, 
\be
\label{loi1}
\qquad \int_{\RR^6}G(f(t,x,v))\,dxdv=\int_{\RR^6}G(f_0(x,v))\,dxdv\,.
\ee
This last property induces a continuum of conservation laws and is the major difference between this kind of problem and other nonlinear dispersive problems like nonlinear wave or Schr\"odinger equations.


\subsection{Nonlinear stability of steady state solutions}


A classical problem which has attracted a considerable amount of work both in the astrophysical \cite{A1,A2,kandrup1,kandrup2,lynden,lynden2, WZS, aly} and mathematical communities, is the question of the nonlinear stability of stationary states. If we restrict our study to radially symmetric stationary states --that is a priori depending on $(|x|,|v|,x\cdot v)$ only--, Jean's theorem \cite{Batt} ensures that they can be described as functions of their own microscopic energy and their angular momentum:
\be
\label{Pdefejdoj}
e(x,v)=\frac{|v|^2}{2}+\phi_Q(x),\qquad \ell(x,v)=|x\times v|^2,
\ee
        
\be     
\label{nvknvdn}
        Q(x,v)=F\left(e(x,v),\ell(x,v)\right).
        \ee
The existence of such steady states has been discussed in \cite{Batt} for a large class of smooth functions $F$. A well-known conjecture in astrophysics, \cite{binney}, is now that among these stationary solutions, those who are  {\it nonincreasing} functions of their microscopic energy $e$ are nonlinearly stable by the Vlasov Poisson flow, explicitly:\\

\ni
{\it Conjecture}: Non increasing anisotropic galaxies $F=F(e,\ell)$ with $\frac{\partial F}{\partial e}<0$ on the support of $Q$ are stable by spherically symmetric perturbations for the flow (\ref{vp}). Non increasing isotropic spherical galaxies $F=F(e)$  with $\frac{\partial F}{\partial e}<0$ on the support of $Q$ are orbitally stable against general perturbations for the flow (\ref{vp}).\\
   
Remarkably enough, this conjecture has been proved at the linear level  by Doremus, Baumann and Feix \cite{doremus} (see also \cite{gillon,kandrup1,SDLP} for related works), following the pioneering work by Antonov in the 60's \cite{A1,A2}. These results are based  on some coercivity properties of the linearized Hamiltonian under constraints formally arising from the linearization of the Casimir conservation laws (\ref{loi1}), see Lynden-Bell \cite{lynden}.

 At the nonlinear level, the general problem is open. However, the nonlinear stability of  a large class of stationary solutions of so-called {\it ground state} type including the polytropic states has been obtained using variational methods in \cite{wol}, \cite{Guo,GuoRein2,GuoRein,Guo2,Dol}, completed by \cite{S2}. In \cite{LMR-note,LMR1,LMR5},  see also \cite{SS}, we observed that a direct application of Lion's concentration compactness techniques \cite{PLL1,PLL2}, implies that or a large class of convex functions $j$, the two parameters --according to the scaling symmetry of (\ref{vp})-- minimization problem 
 \be
\label{2c}
I(M_1,M_j)=\inf_{|f|_{L^1}=M_1, \ |j(f)|_{L^1}=M_j} \calH (f), \qquad M_1,M_j>0
\ee 
 is attained up to symmetries on a steady state solution to (\ref{vp}) of the form (\ref{nvknvdn}), and all minimizing sequences to (\ref{2c}) are relatively compact up to a translation shift in the natural energy space $$\calE=\{f\geq 0 \ \ \mbox{with} \ \ |f|_{\calE}=|f|_{L^1}+|f|_{L^{\infty}}+||v|^2f|_{L^1}<+\infty\}.$$ The so-called Cazenave, Lions \cite{CL} theory of orbital stability then immediately implies the orbital stability of the corresponding ground state steady solution, \cite{LMR1}. In fact, this last step requires the knowledge of the uniqueness of the minimizer to (\ref{2c}) which is a delicate open problem in general, see \cite{S2}, but this difficulty was overcome in \cite{LMR5}.\\

Other non variational approaches based on linearization techniques have also been explored in \cite{W,GR4}. Recently, Guo and Lin  \cite{GuoLin} proved the radial stability of the so called King model $F(e)=(\exp(e_0-e)-1)_+$ which is not in the class of ground states as obtained in the framework of (\ref{2c}). Adapting a robust approach developped by Lin and Strauss in their study of the Vlasov Maxwell system, \cite{lin,lin-strauss1,lin-strauss2}, the authors use the infinity of conservation laws provided by the nonlinear transport  to construct a sufficient large approximation of the kernel of the linearized operator close to the steady state. This allows them to recover a coercivity statement of the linearized energy using Antonov's coercitivity property which after linearization and control of higher order terms for the King model yields the claimed stability in the radial class.


\subsection{Additional conserved quantities in the radial setting}

 
 Our main purpose in this paper is to describe a generalized variational approach for the nonlinear stability of steady states which fully takes into account the nonlinear transport structure of the problem, and in particular the continuum of constraints at hand from (\ref{loi1}).
 
 First recall that in general, the full set of invariant quantities conserved by the nonlinear transport flow \fref{vp} depends on the initial data and its possible symmetries. From now and for the rest of this paper, we shall restrict our attention to {\it spherically symmetric} solutions $f(x,v)=f(|x|,|v|,x\cdot v)$ where we will systematically abuse notations and identify $f$ with its image through various diffeomorphisms. We then let $\calE_{rad}$ be the  space of spherically symmetric distribution functions of finite energy
\be
\label{Erad}
\calE_{rad}=\{f\in \calE, \ \ \ f\mbox{ spherically symmetric}\},
\ee
and recall that if $f$ is spherically symmetric, then $\rho_f(x)=\rho_f(|x|)$ and $\phi_f(x)=\phi_f(|x|)$. This implies in particular from a direct computation that the momentum $$\ell=|x\times v|^2$$ is conserved by the characteristic flow associated to \fref{vp}, and hence a larger class of Casimir conservation laws (\ref{loi1}) holds: 
\be
\label{loi2}
\int_{\RR^6}G(f(t,x,v),|x\times v|^2)dxdv=\int_{\RR^6}G(f_0(x,v),|x\times v|^2)dxdv
\ee
for all $G\in\mathcal C^1([0,+\infty)\times[0,+\infty),\RR^+)$ with $G(0,\ell)=0$, $\forall \ell\geq0$.\\

Let us reformulate \fref{loi2} in terms of {\it equimeasurability} properties of $f$ and $f_0$. Performing the change of variables $$r=|x|,\ \ w=|v|, \ \ x\cdot v=|x||v|\cos\theta, \ \ r,w>0, \ \ \theta\in ]0,\pi[,$$ the Lebesgue measure is mapped onto: $$\int_{\RR^6}f(x,v)dxdv=8\pi^2\int_{r=0}^{\infty}\int_{w=0}^{+\infty}\int_{\theta=0}^{\pi}f(r,w,\cos\theta)r^2w^2\sin\theta drdwd\theta.$$ We then perform the second change of variables $$r=r, \ \ u=w\,\mbox{sign}(\cos\theta), \ \ \ell=r^2w^2\sin^2\theta$$ and get from Fubini:
\be
\label{changevariables}
\int_{\RR^6}f(x,v)dxdv=\int_{\ell=0}^{+\infty}\left[\int_{(r,u)\in \Omega_\ell}f(r,u,\ell)d\nu_\ell\right]d\ell
\ee with 
\be
\label{defomegal}
\Omega_\ell=\{(r,u)\in \RR^+\times \RR \ \ \mbox{with} \ \ r^2u^2>\ell\}
\ee
and
 \be
 \label{nul}
 d\nu_\ell= 4 \pi^2  \un_{r^2u^2> \ell}\,  (r^2u^2 - \ell)^{-1/2} r |u| dr du.
 \ee 
 We then define the distribution function of $f$ at given kinetic momentum $\ell$: $\forall \ell>0$, $\forall s\geq 0$, 
 \be
\label{defmuf}
\mu_f(s, \ell)= \nu_\ell\left\{(r,u)\in \Omega_\ell,\ \ f(r,u,\ell)>s \right\},
\ee
or equivalently 
 \be
\label{defmuf-rul}
\mu_f(s, \ell)=4 \pi^2   \int_{r=0}^{+\infty} \int_{u=-\infty}^{+\infty}\un_{f(r,u,\ell)>s}  (r^2u^2 - \ell )^{-1/2} r |u|  \un_{r^2u^2>\ell}\, dr du \,.
\ee
We now define the set of distribution functions which are equimeasurable to $f$ at given $\ell$ by:
\be
\label{defeqf}
\Eq(f)=\{g\geq 0\ \ \mbox{spherically symmetric},  \ \ \forall s>0, \ \ \mu_f(s,\ell)=\mu_g(s,\ell) \ \ a.e.  \  \ \ell\}.
\ee
We then have from standard arguments:

\begin{lemma}[Characterization of $\Eq(f)$]
\label{propeqf}
Let $f\in L^1\cap L^{\infty}$, nonnegative and spherically symmetric, then the following are equivalent:\\
(i) $g \in \Eq(f)$;\\
(ii) $\forall G(h,\ell)\geq 0$, $\mathcal C^1$  with $G(0,\ell)=0$, there holds:
$$\int_{\RR^6}G(f(x,v),|x\times v|^2)dxdv=\int_{\RR^6}G(g(x,v),|x\times v|^2)dxdv.$$
\end{lemma}

Lemma \ref{propeqf} allows us to reformulate the conservation laws of the full Casimir class \fref{loi2} in the radial setting as follows:
\be
\label{conslooi3}
\forall t\geq0,  \ \ f(t)\in \Eq(f_0).
\ee


\subsection{Assumption (A) on the steady state}


 Before stating the results,  let us fix our assumptions on the steady state $Q$.
  \begin{itemize}
\item[(i)] $Q$ is a continuous, nonnegative,  non zero, compactly supported steady state solution of the Vlasov-Poisson system \fref{vp}.\\
\item[(ii)] There exists a continuous function $F:\RR\times \RR_+\to \RR_+$ such that 
\be
\label{Q}
\forall (x,v)\in \RR^6, \ \ Q(x,v)=F\left(\frac{|v|^2}{2}+\phi_Q(x), |x\times v|^2\right).
\ee
\item[(iii)] There exists $e_0<0$ such that:\\
$$\begin{array}{l}
\calO=\{(e,\ell)\in \RR\times\RR_+\,:\,\, F(e,\ell)>0\}\subset ]-\infty,e_0[\times \RR_+,\\[3mm]
F \mbox{ is $\calC^1$  on }\calO,\mbox{ with }\ds\frac{\partial F}{\partial e}<0.
\end{array}
$$
\end{itemize}

\begin{remark}
Note that $\frac{\partial F}{\partial e}$ may be infinite at the boundary of $\calO$, as is the case for polytropic ground states $F(e,\ell)=(e_0-e)_+^{q} \ell^{\kappa},$ for some $ 0<q<1$ and $ \kappa \geq 0$. \end{remark}

Below we list a number of physically relevant models  for which our non linear stability result applies. All these  examples are extracted from \cite{binney} to which we refer for a detailed
physical description of  various gravitational models.

\bs

\ni
{\bf Examples.}
\begin{itemize}
\item[--]  {\em Polytropes and double-power models: } The polytropes correspond to the following form of $F$:
$$F(e,\ell)= (e_{0}-e)_{+}^{q} \ell ^{\kappa}, \qquad 0<q< 7/2,\qquad \kappa\geq 0,$$
where $e_{0}<0$ is a  constant  threshold energy.  A generalization of these polytropes is provided by the so-called double-power model \cite{binney}:
$$F(e, \ell)=  \sum_{0\leq i,j\leq N} \alpha_{ij} (e_{0}-e)_{+}^{q_i}\ell^{\kappa_j} ,$$
where  $\alpha_{ij}$ are  nonnegative constants.
\item[--] {\em  Michie-King models:}
$$F(e, \ell)= \exp(-\ell/2r_{a}^2) \left(\exp(e_{0}-e) -1\right)_{+},$$
where $e_{0}<0$ and the constant $r_{a}>0$ is the anisotropy radius \cite{binney}.
When  $r_{a}$ goes to infinity, this model reduces to the {\em King} model.
\item[--] {\em Osipkov-Merritt models:}
$$F(e, \ell)= G\left(e_{0}-e+\frac{\ell}{2r_{a}^2}\right),$$
where $e_{0}<0$,  $r_{a}>0$ are constants, and $G$ is a nonincreasing $C^1$ function  such that $G(t)=0$ for all $t\leq 0$. 
\end{itemize}


\subsection{Statement of the results}


From \fref{conslooi3}, a natural generalization of (\ref{2c}) in the radial setting is to minimize the Hamiltonian under constraints of given equimeasurability. This is a very natural strategy to prove stability in a nonlinear transport setting which goes back in fluid mechanics to the celebrated works of Arnold, see e.g \cite{arnold1}, \cite{arnold2}, \cite{arnoldbook}, Marchioro and Pulvirenti \cite{MPbook}, \cite{MP} , Wolansky  and Ghil \cite{WG}, and references therein, and is also very much present in the  physics litterature, see in particular Lynden-bell \cite{lynden}, Gardner \cite{gardner}, Wiechen, Ziegler, Schindler \cite{WZS}, Aly \cite{aly} and references therein. The mathematical implementation of the corresponding variational problem is however confronted to the description of bounded sequences in $\Eq(f_0)$ and a possible lack of compactness in general, see for example Alvino, Trombetti and Lions \cite{ATL} for an introduction to this kind of problem.\\

Our first result is the characterization of non increasing states as local minimizers of the Hamiltonian in $\calE_{rad}$ under a constraint of equimeasurability:

 \begin{theorem}[Local variational characterization of $Q$]
\label{theo1}
There exists a constant $C_0>0$ such that the following holds. For all $R>0$, there exists $\delta_0(R)>0$ such that, for all $f\in \calE_{rad}\cap \Eq(Q)$ satisfying
\be
\label{a-theo1}
|f-Q|_{\calE}\leq R,\qquad |\na\phi_f-\na\phi_Q|_{L^2}\leq \delta_0(R),
\ee
we have
\be
\label{minoration}
\calH(f)-\calH(Q)\geq  C_0 |\na \phi_f-\na \phi_Q|_{L^2}^2\,.
\ee
If in addition $\calH(f)=\calH(Q)$, then $f=Q$.
\end{theorem}
 
Theorem \ref{theo1} was first obtained by Guo, Rein \cite{GR4} for a perturbation $f$ near $Q$ \footnote{and not only $\phi_f$ near $\phi_Q$ which is an issue for the proof of Theorem \ref{theo2}} in the specific case of the isotropic King model, and for isotropic relativistic models $F(e,l)=F(e)$ with locally bounded derivative $F'(e)$ in \cite{HR}, and this excludes any singularity at the boundary --as many polytropic models would have--.\\

Let us stress onto the fact that Theorem \ref{theo1} by itself alone is too weak to yield a stability statement including the full set of radial pertubations. Hence the importance of Theorem  \ref{theo1} relies in fact mostly on its proof. Indeed, a new important feature of our analysis is to use a {\it monotonicity property} of the Hamiltonian under a generalized Schwarz symmetrization which is not the standard radial rearrangement but {\it a rearrangement with respect to a given microscopic energy} $\frac{|v|^2}{2}+\phi(x)$, at fixed angular momentum $|x\times v|^2$, see Proposition \ref{Lemmarearrangement} for a precise definition and Proposition \ref{propclemonotonie} for the monotonicity statement. This monotonicity is very much a consequence of the "bathtub" principle for symmetric rearrangements, see Lieb and Loss \cite{lieb-loss}, and was already observed in the physics litterature, see Gardner \cite{gardner}, Aly \cite{aly}. It produces a reduced functional $\mathcal
  J(\phi_f)$ which {\it depends on the Poisson field $\phi_f$ only and not the full distribution function}. The outcome is a lower bound 
\be
\label{cioyeo}
\mathcal H(f)-\mathcal H(Q)\geq \mathcal {J}(\phi_f)-\mathcal J(\phi_Q).
\ee 
Interestingly enough, the reduced functional $\mathcal J$ was first introduced on physical ground as a {\it generalized potential energy} in the pioneering works by Lynden-Bell \cite{lynden}, see also  Wiechen, Ziegler, Schindler \cite{WZS}. It now turns out from explicit computation that the critical points of $\mathcal J$ are the Poisson field of steady states, and that the Hessian of $\mathcal J$ near the Poisson field of a nondecreasing steady state can be directly connected to {\em the Hartree-Fock exchange operator} \cite{lynden},  which is coercive from Antonov's stability criterion, see section \ref{sectionj}, and hence $\phi_Q$ itself is a local minimizer of $\mathcal J$.\\
 
 The important outcome of the structure \fref{cioyeo} is that by reducing the problem to a problem on the Poisson field only, we are able to extract compactness in the radial setting from {\it any} minimizing sequence whose Hamiltonian converges to $Q$ {\it without the assumption of equimeasurability}, thanks to the smoothing and compactness provided by the radial Poisson equation. This allows us to prove the following compactness result which is the heart of our analysis.\\
 
 Given $f\in \calE_{rad}$, we consider the family of its Schwarz symmetrizations $f^*(\cdot,\ell)$, $\ell >0$, as defined in Proposition \ref{defprop}. We then claim:
  
 \begin{theorem}[Compactness of local minimizing sequences]
\label{theo2}
There exists $\delta>0$ such that the following holds. Let $f_n$ be a sequence of functions of $\calE_{rad}$, bounded in $L^\infty$, such that 
\be
\label{condfn1}
|\na\phi_{f_n}-\na\phi_Q|_{L^2}<\delta,
\ee
and
\be
\label{condfn2}
\underset{\footnotesize n\to +\infty}{\mbox{\rm lim sup}}\,\calH(f_n)\leq \calH(Q),\qquad f_n^* \to Q^*  \mbox { in }  L^1(\RR_+\times \RR_+)\quad \mbox{as }n\to +\infty
\ee
then 
\be
\label{convfn}
f_n\to Q \mbox{ in }L^1(\RR^6),\qquad |v|^2f_n \to |v|^2Q\mbox{ in }L^1(\RR^6).
\ee
\end{theorem}

Theorem \ref{theo2} is the key to the radial Cazenave-Lions'  theory of orbital stability \cite{CL} and implies that any compactly supported non increasing steady state $Q$ as defined by \fref{Q}, is nonlinearly stable under the action of the Vlasov-Poisson flow with respect to spherical perturbations. We thus obtain the main result of this paper:

\begin{theorem}[Nonlinear stability of $Q$ under the nonlinear flow (\ref{vp})]
\label{theo3}
For all $M$ large enough  and for all $\eps>0$, there exists $\eta>0$ such that the following holds true. Let $f_0\in \calE_{rad}\cap\calC^1_c$ with
\be
\label{eta}
|f_0-Q|_{L^1}<\eta,\quad |f_0|_{L^\infty}< M,\quad \calH(f_0)<\calH(Q)+\eta,
\ee
then the corresponding global strong solution $f(t)$ to \fref{vp} satisfies: $\forall t\geq 0$, 
\be
\label{eta2}
|f(t)-Q|_{L^1}+||v|^2(f(t)-Q)|_{L^1}<\eps, \ \ |f(t)|_{L^{\infty}}<M.
\ee
\end{theorem}

\bs
\ni
{\bf Comments on Theorem \ref{theo3}}

\ms
\ni
{\it 1. Linear versus nonlinear stability}. A natural strategy to pass from linear to nonlinear stability is to try to linearize the problem and estimate higher order terms as perturbations. This turns out to be quite delicate in general and the control of higher order terms may be challenging, see \cite{W}, \cite{GuoLin} for a treatment of the King model, \cite{LMR7} for the polytropic case. Our analysis avoids this classical difficulty using two facts. We first derive a {\it global} monotonicity property which is fundamentally a nonlinear property and does not rely on any linearization procedure, Proposition \ref{propclemonotonie}, and which reduces the problem to understanding a simpler functional on the Poisson field $\phi$ only. For this functional, we do apply a linearization procedure that is a Taylor expansion near $\phi_Q$, but we avoid the computation of higher order terms thanks to {\it compactness} properties of the Hessian, see \fref{conti2}, \fref{comp}.

\ms
\ni
{\it 2. Comparison with previous nonlinear stability results}. In view of the nonlinear stability result obtained for ground state type minimizers of (\ref{2c}) which are not restricted to the radial class, one may ask whether a generic steady solution of the form \fref{Q} can in fact be obtained as a ground state  for (\ref{2c}). This is a nontrivial issue which is connected to the notion of equivalence of ensemble in statistical physics. In a forthcoming work \cite{LMR6} and following pioneering ideas from Lieb and Yau \cite{LiebYau}, we will exhibit a large class of monotonic functions $F$ for which the equivalence of ensemble actually holds. There are however of course many well known examples where this equivalence of ensembles fails. Note also that physical investigations around these minimization problems can be found in \cite{chavanis} and the references therein.

\ms
\ni
{\it 3. Comparison with 2D incompressible Euler}. The conservation of equimeasurability properties by the nonlinear transport flow has also been used in the literature to prove the stability of steady states for the 2D incompressible Euler flow, see for example Marchioro, Pulvirenti \cite{MP} and references therein. For a discussion on variational problems with equimeasurability constraints in fluid dynamics, one can also refer to Serre \cite{serre}. Our result generalizes this approach to the Vlasov-Poisson system which is however more delicate due to the non trivial structure of both the Hamiltonian and the steady states solutions.\\

The conjecture of stability of nonincreasing radially symmetric steady states is hence proved for radial perturbations. Note that the result is expected to be optimal for anisotropic galaxies with a non trivial dependence on $\ell$ as some numerical simulations suggest the possible instability of anisotropic models against general perturbations, see \cite{binney}. One important open problem after this work is certainly the general setting of nonradial perturbations for spherical models.


\subsection{Strategy of the proof}


Let us give a brief insight into the proof of the variational characterization of $Q$ given by Theorem \ref{theo1} and the lower bound \fref{cioyeo} which are key  features of our analysis. It follows in three main steps.\\

\ni
{\bf Step 1.} Rearrangement with respect to a given Poisson field.\\

\ni
Let a Poisson field $\phi$ and a radially symmetric distribution function $f\in \mathcal E_{rad}$, we aim at defining the Schwarz symmetrization of $f$ with respect to the microscopic energy $e=\frac{|v|^2}{2}+\phi(x)$ at each given kinetic momentum $\ell$. In other words, given $\ell=|x\times v|^2>0$, we are looking for a function $f^{*\phi}(x,v)=G\left(\frac{|v|^2}{2}+\phi(x),\ell\right)$ which is a nonincreasing function of $e$ and which is equimeasurable to $f$ in the sense of \fref{defmuf}, \fref{defeqf} i.e.: $$\forall t>0, \ \ \mu_f(t,\ell)=\mu_{f^{*\phi}}(t,\ell) \ \ a.e \ \ \ell>0.$$ As a simple change of variables formula similar to \fref{changevariables} reveals, the choice of $f^{*\phi}$ is essentially unique and given by:
\be
\label{deffstarphi}
f^{*\phi}(x,v)=f^*\left(a_{\phi}\left(\frac{|v|^2}{2}+\phi(x), |x \times v|^2\right), |x \times v|^2\right)\un_{\frac{|v|^2}{2}+\phi(x)<0}
\ee
where $f^*$ is the standard Schwarz symmetrization of $f$ at given $\ell$ -see Proposition \ref{defpropschwartz} for a precise statement- and $a_{\phi}$ is the Jacobian of the change of variables, explicitly:
\be
\label{defaphi}
a_{\phi}(e,\ell)=  8 \pi^2\sqrt{2}  \int_0^{+\infty}  \left( e- \phi(r) - \frac{\ell}{2r^2}\right)_+^{1/2} dr.
\ee
Note that the steady state $Q$ being by assumption a nonincreasing function of its microscopic energy, it is automatically a fixed point for this transformation --see Corollary \ref{coroll}--:
\be
\label{afaxux}
Q=Q^{*\phi_Q}.
\ee

\ni
{\bf Step 2.} Monotonicity of the Hamiltonian under the $f^{*\phi_f}$ rearrangement.\\

\ni
The key property which can be found in the physics litterature, see in particular Aly \cite{aly}, is now the monotonicity of the Hamiltonian \fref{defhamiltonian} under the generalized rearrangement \fref{deffstarphi}, see Proposition \ref{propclemonotonie}: 
\be
\label{moninoe}
\forall f\in \mathcal E_{rad}, \ \ \mathcal H(f)\geq \mathcal H(f^{* \phi_f}).
\ee
Pick then $f$ as in the hypothesis of Theorem \ref{theo1} so that $f^*=Q^*$, then a slightly more careful analysis of the monotonicity formula \fref{moninoe} implies a lower bound of the Hamiltonian by a functional {\it which depends on the Poisson field $\phi_f$ only}: 
\be
\label{defclse}
\mathcal H(f) -\mathcal H(Q)\geq \mathcal J(\phi_f)-\mathcal J(\phi_Q)
\ee
with 
$$
\mathcal J(\phi_f)=\mathcal H(Q^{*\phi_f})+\frac{1}{2}\int_{\RR^3}|\nabla \phi_{Q^{*\phi_f}}-\nabla \phi_f|^2\,.
$$
This dependence in $\phi$ only which displays nice {\it compactness properties} in the radial setting is the key to the proof of the convergence of minimizing sequences, Theorem \ref{theo2}. Another important feature in the proof of Theorem \ref{theo2} will be to not only use the monotonicity \fref{moninoe}, but to observe that some {\it norm} is controlled by $\mathcal H(f)-\mathcal H(f^{*\phi_f})$, see section \ref{sectionthmmain}.\\

\ni
{\bf Step 3.} Coercivity of the Hessian of $\mathcal J$ at $\phi_Q$.\\

\ni
From \fref{defclse}, the lower bound \fref{minoration} now follows from the lower bound: 
\be
\label{covhohvo}
\mathcal J(\phi_f)-\mathcal J(\phi_Q)\geq C|\nabla \phi_f-\nabla \phi_Q|_{L^2}^2
\ee
in the vicinity of $\phi_Q$. This local coercivity lower bound relies on an explicit computation of the Taylor expansion of $\mathcal J$ at $\phi_Q$, Proposition \ref{lemA1}. The steady state equation \fref{afaxux} implies that $\phi_Q$ is a critical point of $\mathcal J$, while the Hessian at $\phi_Q$ is intimately related to the  Lynden-Bell Hartree-Fock exchange operator \cite{lynden} , which coercivity was essentially proved 40 years ago by Antonov, \cite{A1}, \cite{A2}, see Proposition \ref{propotransport}, using in particular the fact that {\it in radial symmetry}, the kernel of the linearized transport operator close to $Q$ is explicit. Note that the rigorous derivation of the first two derivatives of $\mathcal J$ at $\phi_Q$ requires a detailed study of the regularity properties of the Jacobian $a_{\phi}$ given by \fref{defaphi} which a priori displays a $\sqrt{\cdot}$ regularity only.\\

This paper is organized as follows. In section \ref{sectnotation}, we introduce the Schwarz symmetrization with respect to the microscopic energy $\frac{|v|^2}{2}+\phi(x)$, at fixed kinetic momentum $|x\times v|^2$, and prove some natural continuity property of the corresponding object $f^{*\phi}$, Proposition \ref{Lemmarearrangement}, and of the Jacobian function $a_{\phi}$, Lemmas \ref{lemmadefaphi}, \ref{contaphi}, \ref{cont-aphi-l-phi}. In section \ref{sectionstar}, we prove the key monotonicity Proposition \ref{propclemonotonie} which reduces the analysis to coercivity properties of the functional $\mathcal J(\phi)$ near $\phi_Q$, Proposition \ref{propA}. We then first conclude the proofs of Theorems \ref{theo1}, \ref{theo2}, \ref{theo3}. assuming Proposition \ref{propA} which is eventually proved in section \ref{sectionj}.
\bs

\begin{acknowledgement}
The authors would like to thank P.-E. Jabin for stimulating discussions about this work, and are endebted to J.-J. Aly for having kindly guided them through the physics reference on the subject and in particular the pioneering important works \cite{lynden}, \cite{WZS}, \cite{aly}. M. Lemou was supported by the Agence Nationale de la Recherche, ANR Jeunes Chercheurs MNEC. F. M\'ehats was supported by the Agence Nationale de la Recherche, ANR project QUATRAIN. P. Rapha\"el was supported by the Agence Nationale de la Recherche, ANR Projet Blanc OndeNonLin and ANR Jeune Chercheur SWAP.
\end{acknowledgement}


\section{Symmetric rearrangement with respect to a given microscopic energy}
\label{sectnotation}


Our aim in this section is to introduce the symmetric rearrangement $f^{*\phi}$ of a distribution function $f\in \calE_{rad}$ with respect to a given microscopic energy $\frac{|v|^2}{2}+\phi(x)$. This notion generalizes the standard Schwarz symmetrization and is the well fitted object for the study of the minimization problem (\ref{minoration}). This symmetrization involves the use  the Jacobian function $a_{\phi}$ given by \fref{defaphi}. We start with proving some continuity and differentiability properties of this functional which will be used all along the paper, Lemma \ref{lemmadefaphi}, \ref{contaphi}, \ref{cont-aphi-l-phi}, and then define $f^{*\phi}$ and give its first properties, Proposition \ref{Lemmarearrangement}.


\subsection{Definition and differentiability in $e$ of the Jacobian function $a_{\phi}$}
\label{sectnotationaphiconitnuity}


Our aim in this subsection is to study the Jacobian $a_{\phi}(e,\ell)$ given by \fref{defaphi}, which appears in the definition of the generalized Schwarz symmetrization \fref{symgenphi}. The class of Poisson potentials $\phi$ which is well fitted for the analysis is
\be
\label{Phirad}
\Phi_{rad} =\left\{ \phi: \RR^3\rightarrow \RR\ \mbox{such that there exists} \  f\in \calE_{rad} \ \ \mbox{with}\  \phi=\phi_f  \right\}
\ee
Let us start with some properties of the so called effective potential appearing in the definition \fref{defaphi} which are elementary but crucial to obtain uniform bounds on $a_{\phi}$ and its various derivatives.

\begin{figure}
\begin{picture}(0,0)%
\epsfig{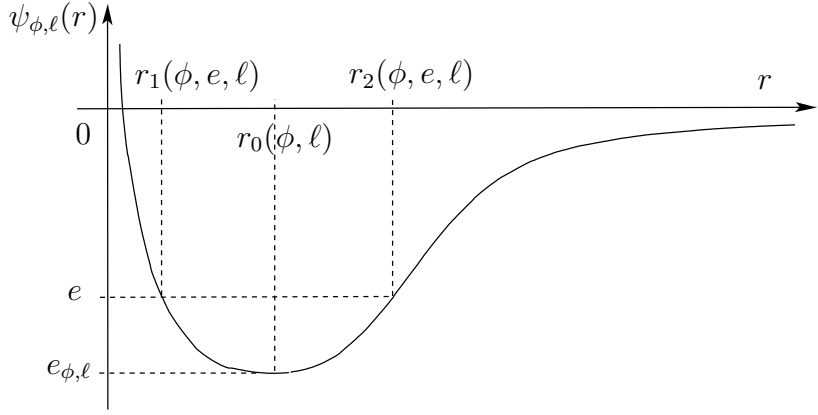}%
\end{picture}%
\setlength{\unitlength}{4144sp}%
\begingroup\makeatletter\ifx\SetFigFont\undefined%
\gdef\SetFigFont#1#2#3#4#5{%
  \reset@font\fontsize{#1}{#2pt}%
  \fontfamily{#3}\fontseries{#4}\fontshape{#5}%
  \selectfont}%
\fi\endgroup%
\begin{picture}(4799,2449)(1,-2093)
\put(  1,209){\makebox(0,0)[lb]{\smash{{\SetFigFont{12}{14.4}{\rmdefault}{\mddefault}{\updefault}{\color[rgb]{0,0,0}$\psi_{\phi,\ell}(r)$}%
}}}}
\put(748,-136){\makebox(0,0)[lb]{\smash{{\SetFigFont{12}{14.4}{\rmdefault}{\mddefault}{\updefault}{\color[rgb]{0,0,0}$r_1(\phi,e,\ell)$}%
}}}}
\put(2014,-133){\makebox(0,0)[lb]{\smash{{\SetFigFont{12}{14.4}{\rmdefault}{\mddefault}{\updefault}{\color[rgb]{0,0,0}$r_2(\phi,e,\ell)$}%
}}}}
\put(4434,-170){\makebox(0,0)[lb]{\smash{{\SetFigFont{12}{14.4}{\rmdefault}{\mddefault}{\updefault}{\color[rgb]{0,0,0}$r$}%
}}}}
\put(224,-1866){\makebox(0,0)[lb]{\smash{{\SetFigFont{12}{14.4}{\rmdefault}{\mddefault}{\updefault}{\color[rgb]{0,0,0}$e_{\phi,\ell}$}%
}}}}
\put(396,-496){\makebox(0,0)[lb]{\smash{{\SetFigFont{12}{14.4}{\rmdefault}{\mddefault}{\updefault}{\color[rgb]{0,0,0}$0$}%
}}}}
\put(349,-1428){\makebox(0,0)[lb]{\smash{{\SetFigFont{12}{14.4}{\rmdefault}{\mddefault}{\updefault}{\color[rgb]{0,0,0}$e$}%
}}}}
\put(1351,-511){\makebox(0,0)[lb]{\smash{{\SetFigFont{12}{14.4}{\rmdefault}{\mddefault}{\updefault}{\color[rgb]{0,0,0}$r_0(\phi,\ell)$}%
}}}}

\end{picture}%

\caption{Profile of the effective potential $\psi_{\phi,\ell}(r)$}
\label{figure}

\end{figure}

\begin{lemma}[Structure of the effective potential for $\phi \in \Phi_{rad}$]
\label{lemmapotentials}
Let $\phi=\phi_f \in \Phi_{rad}$, be non zero. For $\ell>0$,  consider the effective potential
\be
\label{psiel}
\psi_{\phi,\ell}(r) = \phi(r) + \frac{\ell}{2r^2}, \ \ \ \ \ \  r>0.
\ee
(i) Structure of $\psi_{\phi,\ell}$: $\psi_{\phi,\ell}\in \mathcal C^1(\RR^3\backslash \{0\})$ and
\be
\label{ephil}
e_{\phi, \ell} = \inf_{r\geq 0} \left[ \psi_{\phi,\ell}(r)\right] ,
\ee
is attained at a unique $r_0(\phi,\ell)$. $\psi_{\phi,\ell}$ is strictly decreasing on $(0,r_0(\phi,\ell))$ and strictly increasing on $(r_0(\phi,\ell),+\infty)$ with 
\be
\label{ochohieo}
\lim_{r\to 0}\psi_{\phi,\ell}(r)=+\infty, \ \ \lim_{r\to +\infty}\psi_{\phi,\ell}(r)=0.
\ee
Moreover, the function $\ell\mapsto e_{\phi,\ell}$ is continuous on $\RR_+^*$, with the uniform bound: 
\be
\label{mine2}
\forall \ell>0, \qquad  \max\left( \phi(0),  -\frac{|f|_{L^1}^2}{2\ell}\right)\leq e_{\phi, \ell}<0.
\ee
(ii) Level sets of $\psi_{\phi,\ell}$: for $e_{\phi,\ell}<e<0$, let
\be
\label{r1}
r_1(\phi,e,\ell)=\inf\left\{r\geq 0 \ st. \ e-\psi_{\phi,\ell}(r)>0\right\} , 
\ee
\be
\label{r2}
 r_2(\phi,e,\ell)=\sup\left\{r\geq 0 \ st. \ e-\psi_{\phi,\ell}(r)>0\right\}.
\ee
Then $r_1(\phi,e,\ell), r_2(\phi,e,\ell)$ are $\mathcal C^1$ functions of $e$ with uniform bounds: $\forall e_{\phi,\ell}< e<0$:
\be
\label{majr2-simple}
0<\frac {\ell}{2|f|_{L^1}}\leq r_1(\phi,e,\ell) < r_2(\phi,e,\ell) \leq \frac {-|f|_{L^1}}{e}.
\ee
(iii) Concavity lower bound: there holds the uniform concavity lower bound
$\forall e\in (e_{\phi,\ell},0)$, $\forall r\in  [r_1(e,\phi,\ell),r_2(e,\phi,\ell)]$,
\be
\label{clvm}
e-\psi_{\phi,\ell}(r)\geq \frac{\ell}{2r^2\,r_1r_2}(r-r_1(\phi,e,\ell))(r_2(\phi,e,\ell)-r).
\ee
\end{lemma}
On Figure \ref{figure}, we summarize the properties of $\psi_\ell$ described above.

\begin{remark} In the sequel and when there is no ambiguity, we will avoid the $(\phi,e,\ell)$ dependence and note $r_0,r_1,r_2$.
\end{remark}

\begin{proof} The proof is elementary but relies on a crucial way on the {\it positivity} of $\Delta\phi_f$.\\
Let us recall  the standard interpolation estimate for $f\in\calE$:
\be
\label{interpolation}
|\na\phi_f|_{L^2}^2\leq C||v|^2f|_{L^1}^{1/2}|f|_{L^1}^{7/6}|f|_{L^\infty}^{1/3}\,.
\ee
Let $\phi=\phi_f\in \Phi_{rad}$, then by interpolation and Sobolev embedding,  $\rho_f\in L^{5/3}(\RR^3)$ and thus $\phi_f\in W^{2,5/3}_{loc}(\RR^3)\subset \calC^0(\RR^3)$ and $\phi_f\in \mathcal C^1(\RR^3\backslash \{0\})$ by elliptic regularity and the radial assumption, from which $\psi_{\phi,\ell}\in \mathcal C^1(\RR^3\backslash \{0\})$.\\
We now integrate the radial Poisson equation and get:
\be
\label{poisson}
r^2\phi_{f}'(r)=4\pi \int_0^rs^2\rho_f(s)ds\geq 0, \qquad  \lim_{r\to +\infty} r\phi_f(r)=-|f|_{L^1}.
\ee
Note that the second identity is obtained by integrating the first one as follows:
\be
\label{be}
r\phi_f(r)=-4\pi\int_0^rs^2\rho_f(s)ds-4\pi r\int_r^{+\infty}s\rho_f(s)ds.
\ee
We deduce that $\phi=\phi_f$ is continuous, nondecreasing and nonpositive on $[0,+\infty[$ with
\be
\label{minphi}
\phi(r) \geq \max\left(\phi(0), - \frac{ |f|_{L^1}}{r}\right), \ \ \ \forall r\geq 0, 
\ee
and  there exists $\tilde{r}>0$,  such that 
\be
\label{majphi}
\phi(r)\leq -\frac{|f|_{L^1}}{2r},  \ \ \   \forall \ r\geq \tilde{r}.
\ee
Thus \fref{minphi}, \fref{majphi} imply \fref{ochohieo}. From \fref{majphi}, $e_{\phi,\ell}$ given by \fref{ephil} satisfies
$$e_{\phi, \ell} \leq \inf_{r\geq \tilde{r}} \left[ -\frac{|f|_{L^1}}{2r}+ \frac{\ell}{2r^2}\right] <0,$$
since by assumption $f\neq 0$,
 and hence $e_{\phi,\ell}$ is attained at some $r_0=r_0(e,\phi,\ell)$. Thus from  \fref{minphi}:

$$
e_{\phi, \ell} = \phi(r_0)+\frac{\ell}{2 r_0^2}\geq\max\left( \phi(0),  -\frac{|f|_{L^1}}{r_0}+ \frac{\ell}{2 r_0^2}\right)\geq\max\left( \phi(0),  -\frac{|f|_{L^1}^2}{2\ell}\right)\ \ \ \forall \ell >0,$$ and \fref{mine2} is proved.

Observe now from  \fref{poisson} again that:
\be
\label{psi'}
\psi_{\phi,\ell}'(r) =\phi'(r) -\frac{\ell}{r^3}, \ \ \mbox{and}\ \ (r^2\psi_{\phi,\ell}'(r))'= r^2\rho_f +\frac{\ell}{r^2}>0,
\ee
and hence from $\psi'_{\ell}(r_0(e,\phi,\ell))=0$:
\be
\label{nume}
\forall r>0, \ \ r^2\psi_{\phi,\ell}'(r) = \int_{r_0}^{r}\left( r^2\rho_f(r) + \frac{\ell}{r^2}\right)dr,
\ee
which yields the uniqueness of the minimum $r_0>0$ and the claimed monotonicity properties of $\psi_{\phi,\ell}$.

Together with \fref{ochohieo}, we conclude from \fref{nume} that $r_1,r_2$ given by \fref{r1}, \fref{r2} are well defined for $e_{\phi,\ell}<e<0$, and are $\mathcal C^1$ functions of $e$ from the implicit function theorem. To prove the uniform bound \fref{majr2-simple}, we observe from \fref{minphi}:
$$\left\{r\geq 0; \ st. \ e-\phi(r) - \frac{\ell}{2r^2}>0\right\} \subset \left\{r\geq 0; \ st. \ e+\frac{|f|_{L^1}}{r} - \frac{\ell}{2r^2}>0\right\},$$
and hence using from \fref{mine2} that $|f|_{L^1}^2+2e\ell>0$ for $e>e_{\phi,\ell}$:
$$
\left\{r\geq 0; \ st. \ e+\frac{|f|_{L^1}}{r} - \frac{\ell}{2r^2}>0\right\} \subset \left[\frac{\ell}{|f|_{L^1}+\sqrt{|f|_{L^1}^2+2e\ell}} ,\frac{\ell}{|f|_{L^1}-\sqrt{|f|_{L^1}^2+2e\ell}} \right].
$$
We then use the definitions \fref{r1} and \fref{r2} to get
$$
0<\frac {\ell}{|f|_{L^1}+\sqrt{|f|_{L^1}^2+2e\ell}}\leq r_1(\phi,e,\ell) < r_2(\phi,e,\ell) \leq \frac {\ell}{|f|_{L^1}-\sqrt{|f|_{L^1}^2+2e\ell}},
$$ which implies \fref{majr2-simple}.

Let us now prove the continuity of the function $\ell \mapsto e_{\phi,\ell}$ on $\RR_+^*$. Let $0<\ell_1<\ell_2$ be fixed. From the definitions \fref{psiel} and \fref{ephil}, for all $\ell\in[\ell_1,\ell_2]$ we have 
$$e_{\phi,\ell}\leq e_{\phi,\ell_2}$$
thus, applying \fref{majr2-simple} with $e=\frac{1}{2}e_{\phi,\ell_2}$ gives
$$\alpha_1 \leq r_0(\phi,\ell)\leq \alpha_2,\quad \mbox{with}\quad\alpha_1=\frac{\ell_1}{2|f|_{L^1}}>0,\quad \alpha_{2} = \frac{2|f|_{L^1}}{|e_{\phi,\ell_2}|}>0.$$
Hence, $(r,\ell)\mapsto \psi_{\phi,\ell}(r)$ being continuous,  the function $$\ell \in [\ell_1,\ell_2]\mapsto e_{\phi,\ell}=\min_{r\in[\alpha,\alpha_{2}]}\psi_{\phi,\ell}(r)$$ is continuous.

It remains to prove the concavity bound \fref{clvm}. Let
$$w(r)=e-\psi_{\phi,\ell}(r)-\frac{\ell}{2r^2\,r_1r_2}(r-r_1)(r_2-r),$$
then
$$(rw(r))''=\frac{1}{r}\left(-(r^2\psi_{\phi,\ell}'(r))'+\frac{\ell}{r^2}\right)=-r\rho_f(r)\leq 0,$$
where we used \fref{psi'}. Hence the function $r\mapsto rw(r)$ is concave. Since it vanishes at $r_1$ and $r_2$, we conclude that $w(r)\geq 0$ for all $r\in  [r_1,r_2]$ and \fref{clvm} is proved.

This concludes the proof of Lemma \ref{lemmapotentials}.
\end{proof}

Let us now define the Jacobian function $a_{\phi}(e,\ell)$ and examine its differentiability properties in $e$ :

\begin{lemma}[Definition and differentiability properties in $e$ of the Jacobian $a_{\phi}$]
\label{lemmadefaphi}
For $\phi=\phi_f\in \Phi_{rad}$ non zero and $\ell>0$, we define: 
\be
\label{aphi-rul}
a_{\phi}(e,\ell)=\left\{\begin{array}{ll} 
        \nu_\ell\left\{(r,u)\in(\RR_+)^2:\, \frac{u^2}{2}+\phi(r)<e \right\} \ \ \mbox{for} \ \ e<0 \  \mbox{and } \  \ell >0, \\
        +\infty \ \ \mbox{for} \ \ e\geq 0, \mbox{and}\   \ell>0,\end{array}\right. 
\ee
where $\nu_\ell$ is the measure given by \fref{nul}, equivalently: $\forall \ell>0$, $\forall e<0$,
\be
\label{defaphijo}
a_{\phi}(e,\ell)=  8 \pi^2\sqrt{2}  \int_{r_1}^{r_2}  \left( e- \psi_{\phi,\ell}(r)\right)^{1/2} dr.
\ee
Then:\\
(i) Behavior of $a_{\phi}$: $a_{\phi}(e,\ell)=0$ for $e<e_{\phi,\ell}$ and:
\be
\label{limaphi}
\forall \ell >0, \ \ a_{\phi}(e_{\phi, \ell},\ell)=0 ,  \ \ \lim_{e\to 0^-} a_{\phi}(e,\ell)=+\infty.
\ee
(ii) Uniform bounds on $a_{\phi}$: let 
\be
\label{mphi} 0< m_{\phi}:=\inf_{r\geq  0} (r+1)|\phi (r)| <+\infty,
\ee
then there holds the bounds:
\be
\label{minaphi} 
\forall e<0, \ \ a_{\phi}(e,\ell)\leq 16\pi^2  |e|^{-1/2} |f|_{L^1}, 
\ee
and 
\be
\label{estiminf}
\forall e\in (-\frac{m_{\phi}^2}{4(2m_{\phi}+\ell)},0), \qquad a_{\phi}(e,\ell)\geq \frac{4\pi^2}{3}  |e|^{-1/2} m_{\phi}.
\ee
(iii) Differentiability in $e$: the map $e \rightarrow a_{\phi}(e,\ell)$  is a $\calC^1$-diffeomorphims from $(e_{\phi, \ell},0)$ to $(0,+\infty)$ with: 
\be
\label{aphi'}
 \forall e\in ]e_{\phi,\ell},0[, \ \ \frac{\partial a_\phi}{\partial e}(e,\ell)= 4\pi^2 \sqrt{2} \int_{r_1}^{r_2} \left(e- \psi_{\phi,\ell}(r)\right)^{-1/2}dr>0.
\ee
Abusing notations, we shall denote in the sequel $a_\phi^{-1}(\cdot,\ell):(0,+\infty)\to (e_{\phi,\ell},0)$ its inverse function.
\end{lemma}

\begin{proof} {\it Step 1.} Bounds on $a_{\phi}$.\\

First compute from the definitions \fref{aphi-rul} and \fref{nul}: $\forall e<0$, $\forall \ell >0$ : 
$$
\begin{array}{ll}
\ds a_{\phi}(e,\ell)& \ds =  8\pi^2  \int_{r>0}\int_{u>0}\un_{\frac{u^2}{2}+\phi(r)<e} \un_{r^2u^2>\ell} \ (r^2u^2 - \ell)^{-1/2} r u drdu\\
 & \ds = 8\pi^2  \int_{r>0}\un_{ e- \phi(r) - \frac{\ell}{2r^2}>0}\left(\int_{u=\frac{\sqrt{\ell}}{r}}^{\sqrt{2(e-\phi(r))}}(r^2u^2 - \ell)^{-1/2}  u du\right) rdr\\
 &  \ds = 8\pi^2 \sqrt{2}\int_{r_1}^{r_2} \left(e- \phi(r) - \frac{\ell}{2r^2}\right)_+^{1/2}dr,
 \end{array}
$$
this is \fref{defaphijo} or, equivalently, \fref{defaphi}. Then $a_{\phi}(e,\ell)=0$ for $e\leq e_{\phi,\ell}$ and $a_{\phi}(e,\ell)>0$ on $(e_{\phi,\ell},0)$ from Lemma \ref{lemmapotentials}.

We now estimate $a_\phi$ from above for $e<0$ using \fref{minphi} and  \fref{majr2-simple} as follows
$$\begin{array}{ll} \ds a_{\phi}(e,\ell)&\ds = 8\pi^2 \sqrt{2} \int_{r_1(\phi,e,\ell)}^{r_2(\phi,e,\ell)} \left(e- \phi(r) - \frac{\ell}{2r^2}\right)^{1/2}dr \\[3mm]
&\ds  \leq 8\pi^2 \sqrt{2} \int_{r_1(\phi,e,\ell)}^{r_2(\phi,e,\ell)} \left(\frac{|f|_{L^1}}{r}\right)^{1/2}dr \\[3mm]
&\ds \leq 16\pi^2 \sqrt{2} |f|_{L^1}^{1/2}\left(r_2(\phi,e,\ell)^{1/2}- r_1(\phi,e,\ell)^{1/2}\right)\\[3mm]
& \ds \leq 16\pi^2  |f|_{L^1}|e| ^{-1/2},
\end{array}
$$
and \fref{minaphi} is proved.
To estimate $a_\phi(e,\ell)$ from below, first observe that \fref{mphi} follows from \fref{poisson}. We then write:
\bee
a_{\phi}(e,\ell) &\geq&  8\pi^2 \sqrt{2}   \int_{0}^{ +\infty} \left(e+\frac{m_{\phi}}{r+1}  - \frac{\ell}{2r^2}\right)_{+}^{1/2}dr\\
&\geq&  8\pi^2 \sqrt{2}   \int_{1+ \ell/m_{\phi}}^{ +\infty} \left(e+\frac{m_{\phi}}{r+1}  - \frac{\ell}{2r^2}\right)_{+}^{1/2}dr
\eee
and observe that
for $r\geq 1+ \ell/m_{\phi}$ , we have
$$
\frac{m_{\phi}}{r+1}  - \frac{\ell}{2r^2} \geq \frac{m_{\phi}}{r+1}  - \frac{\ell}{2(r^2-1)}\geq \frac{m_{\phi}}{r+1} \left(1- \frac{\ell}{2m_{\phi}(r-1)}\right) \geq \frac{m_{\phi}}{2(r+1)}.$$
Thus:
\bee
a_{\phi}(e,\ell)&\ds \geq&  8\pi^2 \sqrt{2} \int_{1+ \ell/m_{\phi}}^{+\infty}\left(e+ \frac{m_{\phi}}{2(r+1)}\right)_+^{1/2}dr \nonumber\\
&\ds  \geq& 8\pi^2 \sqrt{2}   \left(\frac{|e|}{m_{\phi}}\right)^{1/2}\int_{1+ \ell/m_{\phi}}^{ +\infty} \left(m_{\phi}- 2|e|(r+1)\right)_{+}^{1/2}dr \nonumber \\
&\ds  \geq& 8\pi^2 \frac{\sqrt{2}}{3}  |e|^{-1/2} m_{\phi} \left(1- 2|e|\frac{2m_{\phi}+\ell}{m_{\phi}^2}\right)_{+}^{3/2}.
\eee
This means that for $ |e|\leq \frac{m_{\phi}^2}{4(2m_{\phi}+\ell)}$, 
$
a_{\phi}(e,\ell)\geq \frac{4\pi^2}{3}  |e|^{-1/2} m_{\phi},
$
and \fref{limaphi} and \fref{estiminf} are proved. The continuity and the monotonicity of the application $e\mapsto a_\phi(e,\ell)$ is a consequence of \fref{majr2-simple} and of the dominated convergence theorem, since
$$\left(e- \phi(r) - \frac{\ell}{2r^2}\right)_+^{1/2} \leq  (-\phi(0))^{1/2}, \ \ \ \mbox{for all}\ \  r\in ]e_{\phi, \ell}, 0[.$$

\bs
\ni
{\it Step 2.} Differentiability of $a_{\phi}$.\\

We are now in position to prove the differentiability of the function $e\to a_\phi(e,\ell)$ which follows from  the version of Lebesgue's derivation theorem given by Lemma \ref{lebesgue}. Let us fix $\ell>0$ and write
$$a_\phi(e,\ell)=\int_0^{+\infty} g(e,r)dr$$
with
$$g(e,r)=8 \pi^2\sqrt{2} \left(e- \psi_{\phi,\ell}(r)\right)^{1/2}\un_{r_1(\phi,e,\ell)<r<r_2(\phi,e,\ell)}=8 \pi^2\sqrt{2} \left(e- \psi_{\phi,\ell}(r)\right)^{1/2}_+\,,$$ and with $\psi_{\phi,\ell}$ given by \fref{psiel}. 
Let $e_0\in ]e_{\phi,\ell},0[$ and $I=]e_0-\eps,e_0+\eps[$, with $\eps$ small enough such that $I\subset  ]e_{\phi,\ell},0[$. Let us check the assumptions of Lemma \ref{lebesgue}. By \fref{minphi}, we have
$$
0\leq g(e,r)\leq  C\left(\frac{|f|_{L^1}}{r}\right)^{1/2}\un_{r_1(\phi,e,\ell)<r<r_2(\phi,e,\ell)}\\
\leq \frac{C}{\sqrt{\ell}}|f|_{L^1}\un_{0<r<\frac{|f|_{L^1}}{2|e_0+\eps|}}
$$
where we used \fref{majr2-simple}. Hence, by standard dominated convergence, $g\in \calC^0(I,L^1(\RR_+))$. Moreover, from the $\mathcal C^1$ regularity of the boundary $r_1,r_2$ with respect to $e$ and the cancellation $e-\psi_{\phi,\ell}(r_1)=e-\psi_{\phi,\ell}(r_2)=0$, the distributional partial derivative of $g$ is given by
$$\frac{\pa g}{\pa e}(e,r)=4 \pi^2\sqrt{2} \left(e- \psi_{\phi,\ell}(r)\right)^{-1/2}\un_{r_1(\phi,e,\ell)<r<r_2(\phi,e,\ell)}\,,$$
and hence the continuity of $r_1,r_2$ with respect to $e$ implies the a.e. convergence in $r$ of this function when $e\to e_0$:
\be
\label{cvpp}
\frac{\pa g}{\pa e}(e,r)\to \frac{\pa g}{\pa e}(e_0,r)\quad \mbox{as }e\to e_0, \,\mbox{ for all } r\neq r_1(\phi,e_0,\ell), \,\,r\neq r_2(\phi,e_0,\ell).
\ee
Now from the concavity estimate \fref{clvm}:
\bea
0\leq \frac{\pa g}{\pa e}(e,r)&\leq &\frac{C}{\sqrt{\ell}}\frac{r\sqrt{r_1r_2}}{\sqrt{(r-r_1)(r_2-r)}}\un_{r_1<r<r_2}\nonumber \\
&\leq &\frac{C}{\sqrt{\ell}}\frac{|f|_{L^1}^2}{|e_0+\eps|^2}\frac{1}{\sqrt{(r-r_1)(r_2-r)}}\un_{r_1<r<r_2}=:q_{e}(r),\label{majq}
\eea
where we applied \fref{majr2-simple} and recall that $r_1=r_1(\phi,e,\ell)$, $r_2=r_2(\phi,e,\ell)$. We observe that
\be
\label{pi}
\int_{r_1}^{r_2}q_{e}(r) dr=\frac{C}{\sqrt{\ell}}\frac{|f|_{L^1}^2}{|e_0+\eps|^2}\pi,
\ee
which implies in particular that $q_e\in L^1(\RR_+)$ and $ \frac{\pa g}{\pa e}(e,r)\in L^1(I\times \RR_+)$. 
Together with \fref{cvpp}, this implies that Assumption (i) of Lemma \ref{lebesgue} is satisfied. 

Let us now check Assumption (ii). From the continuity of $r_1(\phi,e,\ell)$ and $r_2(\phi,e,\ell)$ with respect to $e$, we deduce that
$$q_e(r)\to q_{e_0}(r)\quad \mbox{as }e\to e_0, \,\mbox{ for all } r\neq r_1(\phi,e_0,\ell), \,r\neq r_2(\phi,e_0,\ell).$$
This a.e. convergence, coupled to the fact that, by \fref{pi}, the integral of the positive functions $q_e$ is independent of $e$, is enough to conclude, thanks to the Br\'ezis-Lieb theorem (see Theorem 1.9 of \cite{lieb-loss}), that $q_e$ converges to $q_{e_0}$ in $L^1(\RR_+)$ as $e\to e_0$. Assumption (ii) is then satisfied. Hence Lemma \ref{lebesgue} can be applied and $a_{\phi}(e,\ell)$ is $\calC^1$ with respect to $e$ with its derivative given by \fref{aphi'}. From  Lemma \ref{lemmapotentials}, $\frac{\partial a_\phi}{\partial e}(e,\ell)>0 $ on $ ]e_{\phi, \ell},0[$, so by \fref{limaphi} $e\mapsto a_\phi(e,\ell)$ is a   $\calC^1$ diffeomorphism from $ ]e_{\phi, \ell},0[$ to $]0, +\infty[$.\\
This concludes the proof of Lemma \ref{lemmadefaphi}.
\end{proof}


\subsection{Regularity properties in $\phi$ of $a_{\phi}$ and $a_{\phi}^{-1}$}
\label{sectnotationaphiderivability}


We continue the analysis of the Jacobian $a_{\phi}$ and claim further continuity and differentiability properties with respect to $\phi$.

\begin{lemma}[Continuity properties of $a_\phi(e,\ell)$ with respect to $\phi$]
\label{contaphi}
Let $f\in \calE_{rad}$ nonzero, let $f_n$ be a bounded sequence in $\calE_{rad}$ and denote $\phi_n=\phi_{f_n}$, $\phi=\phi_f$. Assume that $\na {\phi_{f_n}}\to \na{\phi_f}$ in $L^2(\RR^3)$ as $n\to +\infty$. Then, for all $\ell>0$ fixed, the following convergence properties hold as $n\to +\infty$:
\bea
&&e_{\phi_{n}, \ell} \to e_{\phi, \ell},\label{co1}\\
&&r_1(\phi_n,e,\ell)\to r_1(\phi,e,\ell),\,\,r_2(\phi_n,e,\ell)\to r_2(\phi,e,\ell),\quad \forall e\in]e_{\phi,\ell},0[,\label{co2}\\
&&\inf_{n}\inf_{r>0}(1+r)|\phi_n(r)|>0,\label{defm}\\
&&a_{\phi_{n}}(\cdot,\ell) \to a_{\phi}(\cdot,\ell)\quad\mbox{uniformly on any compact subset of $]-\infty,0[,$}\label{co3}\\
&&a_{\phi_{n}}^{-1}(s,\ell) \to a_{\phi}^{-1}(s,\ell)\quad\mbox{for all $s>0$}.\label{co4}
\eea
\end{lemma}

\begin{proof} {\em Step 1. Convergence of the potentials.}

\ms
\ni
As a standard consequence of interpolation and  $\na \phi_n\to \na \phi$ in $L^2$, we have that 
\be
\label{L53}
\rho_{f_n}\rightharpoonup \rho_f \quad \mbox{ weakly in }L^{5/3}(\RR^3).
\ee
 Moreover, by Sobolev embeddings and elliptic regularity, together with the spherical symmetry of $f_n$, we have:
 \be
\label{uniform}
\phi_n\to \phi \quad \mbox{in }L^\infty(\RR_+)\quad \mbox{as }n\to +\infty
\ee
and
\be
\label{uniform2}
\phi'_n\to \phi' \quad \mbox{in }L^\infty([a,b])\quad \mbox{as }n\to +\infty, \mbox{ for all }0<a<b.
\ee

\bs
\ni
{\em Step 2. Proof of \fref{co1}.}
\nopagebreak

\ms
\ni
To prove that $e_{\phi_{n}, \ell} \to e_{\phi, \ell}$, we first pass to the limit $n\to \infty$ into \fref{ephil}, and get
\be
\label{limsup}
\limsup_{n\to +\infty} \, e_{\phi_{n}, \ell} \leq e_{\phi, \ell}.
\ee
Now, we know that the infimum $e_{\phi_{n}, \ell}$ is attained at $r_n=r_0(\phi_n,\ell)$ and from
\fref{minphi} we have
\be
\label{www} e_{\phi_{n}, \ell} = \phi_n(r_n) + \frac{\ell}{2r_n^2} \geq -\frac{|f_n |_{L^1}} {r_n} + \frac{\ell}{2r_n^2}.
\ee
We observe that
\be
\label{www2}
\inf_{n}|f_n |_{L^1}>0.
\ee
Otherwise we would have, up to a subsequence, $f_n\to 0$ in $L^1(\RR^6)$, and then $\phi_f=0$. This means $\rho_f=\Delta \phi_f=0$, and then $f=0$, which contradicts the assumption $f\neq 0$.
Therefore, \fref{www} and \fref{www2} ensure that the sequence $r_n$ is bounded and bounded away from $0$ since \fref{limsup} implies  $\limsup \, e_{\phi_{n}, \ell}<0$.  Then, any subsequence $r_{n'}$ of $r_{n}$ satisfies  $r_{n'} \to r_0>0$ (up to extraction) and one can pass to the limit in $e_{\phi_{n'}, \ell} = \phi_{n'}(r_{n'}) + \frac{\ell}{2r_{n'}^2} $ to get
$$e_{\phi_{n'}, \ell}\to \phi(r_0) + \frac{\ell}{2r_0^2}\geq e_{\phi, \ell},$$
thus 
\be
\label{liminfbi}
\liminf_{n\to +\infty} \, e_{\phi_{n}, \ell} \geq e_{\phi, \ell}.
\ee
Finally \fref{limsup} and \fref{liminfbi}  imply \fref{co1}.

\bs
\ni
{\em Step 3. Proof of \fref{co2}.}

\ms
\ni
Let $e\in]e_{\phi,\ell},0[$. From Step 1, we have $e\in]e_{\phi_n,\ell},0[$ for $n$ large enough, thus $r_1(\phi_n,e,\ell)$ is well-defined. By Lemma \ref{lemmapotentials},  $r_1(\phi_n,e,\ell)$  is characterized by 
\be
\label{caract1}
e=\phi_n(r_1)+\frac{\ell}{2r_1^2}\quad \mbox{and} \quad \phi'_n(r_1)-\frac{\ell}{r_1^3}<0.
\ee
Moreover, by \fref{majr2-simple}, $r_1(\phi_n,e,\ell)$ lies in a compact interval of $\RR_+^*$. Therefore, after extraction of a subsequence, we have $r_1(\phi_n,e,\ell)\to r_*>0$. Thanks to \fref{uniform} and \fref{uniform2}, one can pass to the limit in \fref{caract1} and obtain
$$e=\phi(r_*)+\frac{\ell}{2r_*^2}\quad \mbox{and} \quad \phi'(r_*)-\frac{\ell}{r_*^3}\leq 0.$$
This is enough to conclude that $r_*=r_1(\phi,e,\ell)$. We have thus proved \fref{co2} for $r_1$. the proof of $r_2(\phi_n,e,\ell)\to r_2(\phi,e,\ell)$ is similar.

\bs
\ni
{\em Step 4. Proof of \fref{defm}.}

\ms
\ni
Let us prove \fref{defm}:
$$
\inf_n{m_{\phi_n}}=:m>0,
$$
where $m_{\phi_n}$ is defined by \fref{mphi}. Assume that $m=0$. Then there exists a sequence $r_n$ such that
$(1+r_n)\phi_n(r_n)\to 0$ as $n\to +\infty$. If $r_n$ is bounded, then, up to a subsequence, it goes to some $r_0\geq 0$ and from \fref{uniform}, we get $\phi(r_0)=0$, which is not possible from \fref{poisson} and $f\neq 0$. Hence $r_n$ is not bounded and, up to a subsequence, it goes to $+\infty$. Let $a>0$. We get from \fref{be}
$$(1+r_n)|\phi_n(r_n)|\geq r_n|\phi_n(r_n)|\geq 4\pi \int_0^{r_n}s^2\rho_{f_n}(s)ds\geq 4\pi \int_0^{a}s^2\rho_{f_n}(s)ds$$
for $n$ large enough. Passing to the limit in this inequality and using \fref{L53}, we get $\rho_f=0$, which again contradicts the assumption $f\neq 0$. We have thus proved \fref{defm}.

\bs
\ni
{\em Step 5. Proof of \fref{co3}.}

\ms
\ni
Now, we prove the uniform convergence of $a_{\phi_n}$. We observe  from \fref{majr2-simple} that
 the interval of integration in the expression  \fref{defaphijo} of $a_{\phi_n}$ is bounded. Thus, the dominated convergence theorem applies since 
$$\left( \ e-\phi_n(r) - \frac{\ell}{2r^2}\right)_+^{1/2} \leq (-\phi_n(0))^{1/2} \leq C,$$
from \fref{uniform}.
This yields $a_{\phi_n}(e,\ell) \to  a_{\phi}(e,\ell)$, for all $e<0$, $\ell >0$. 
Now using the monotonicity of the function $e\to a_{\phi_n}(e,\ell)$ at fixed $\ell$ and applying the second Dini's theorem,
we get the desired uniform convergence.

\bs
\ni
{\em Step 6. Proof of \fref{co4}.}

\nopagebreak
\ms
\ni
Let $(s,\ell)\in (\RR_+^*)^2$. Denote
$$e_n=a^{-1}_{\phi_n}(s,\ell),\qquad e_0=a^{-1}_{\phi}(s,\ell).$$
We will prove that $e_n\to e_0$. From \fref{minaphi}, we get
\be
\label{bs}|e_n|\leq C\left(\frac{ |f_n|_{L^1}}{a_{\phi_n}(e_n,\ell)}\right)^2=C\left(\frac{|f_n|_{L^1}}{s}\right)^2.
\ee

Now we claim that
\be
\label{bi}
|e_n|\geq  C\left(\frac{m}{s}\right)^2>0   \ \ \  \mbox{if} \ \ \ |e_n|\leq  \frac{m^2}{4(\ell + 2m)}.
\ee
Indeed, we first get from \fref{estiminf}
\be
\label{bibi}
|e_n|\geq  C\left(\frac{ m_{\phi_n}}{a_{\phi_n}(e_n,\ell)}\right)^2=C\left(\frac{m_{\phi_n}}{s}\right)^2>0,
\ee
provided that $|e| \leq \frac{m_{\phi_n}^2}{4(\ell + 2m_{\phi_n})}$, with $m_{\phi_n}$ defined by \fref{mphi}.   From  \fref{defm}, we have $m_{\phi_n} \geq  m>0$.   Therefore, \fref{bibi} implies \fref{bi} since the function $t\mapsto  \frac{t^2}{\ell + 2t}$ is increasing.  
 
 We then deduce from  \fref{bs} and \fref{bi} that the sequence $e_n$ belongs to a compact interval of $\RR_-^{*}$ thus, up to a subsequence, we have $e_n\to e_\infty\in \RR_-^*$ as $n\to +\infty$. Using \fref{co3}, we have $$s=a_{\phi_n}(e_n,\ell)\to a_{\phi}(e_\infty,\ell)\quad \mbox{as }n\to +\infty.$$
Hence,
$$a_{\phi}(e_\infty,\ell)=a_{\phi}(e_0,\ell)=s\in (0,\infty).$$
Since $e\mapsto a_{\phi}(e,\ell)$ is invertible from $(e_{\phi,\ell},0)$ onto $(0,\infty)$, we deduce that
$$e_0=a^{-1}_{\phi}(s,\ell)=e_\infty,$$
which means that $e_n\to e_0$ as $n\to +\infty$. The proof of \fref{co4} is complete.\\
This concludes the proof of Lemma \ref{contaphi}.
\end{proof}

Let us now examine the differentiability of $a_\phi$ and $a_\phi^{-1}$ with respect to $\phi$. To shorten the statement of the next lemma, we introduce a few notations. We consider two nonzero potentials $\phi=\phi_f\in \Phi_{rad}$ and $\widetilde \phi=\phi_{\widetilde f}\in \Phi_{rad}$ and set:
\be
\label{not1}
h=\widetilde \phi -\phi.
\ee
For all $\ell>0$ and $\lambda\in [0,1]$, we recall the notation
\be
\label{not2}
e_{\phi+\lambda h,\ell}=\inf_{r\geq 0} \left(\psi_{\phi,\ell}(r)+\lambda h(r)\right),
\ee
where $\psi_{\phi_Q,\ell} (r)$ is defined by \fref{psiel},
and denote
\be
\label{not2,5}
\Omega(\phi,\widetilde \phi,\ell)=\left\{(\lambda,e):\,\lambda\in [0,1]\mbox{ and }e\in ]e_{\phi+\lambda h,\ell},0[\right\}.
\ee
Let $s\in \RR_+^*$ and $\lambda\in [0,1]$. Recall that, by Lemma \ref{Lemmarearrangement}, there exists a unique $e\in ]e_{\phi+\lambda h,\ell},0[$, denoted by $a_{\phi+\lambda h}^{-1}(s,\ell)$, such that $a_{\phi+\lambda h}(e,\ell)=s$. Finally, we set
\be
\label{not4}
M=\max(|f|_{L^1},|\widetilde f|_{L^1}).
\ee

\begin{lemma}[Differentiability of $a_\phi(e,\ell)$ with respect to $\phi$]
\label{cont-aphi-l-phi}\mbox{}\\
Let $\ell>0$ be fixed. Consider $\phi, \widetilde \phi\in \Phi_{rad}$ both nonzero and let $h=\phi-\widetilde \phi$. Then, with the notations \fref{not1}--\fref{not4}, the following holds:\\
(i) The function
$$(\lambda,e)\mapsto a_{\phi+\lambda h}(e,\ell)$$
is a $\calC^1$ function on $\Omega(\phi,\widetilde \phi,\ell)$. Moreover, we have
\be
\label{der1}
\frac{\pa }{\pa \lambda} a_{\phi+\lambda h}(e,\ell)=-4\pi^2 \sqrt{2} \int_{r_1(\phi+\lambda h,e,\ell)}^{r_2(\phi+\lambda h,e,\ell)} \left(e- \psi_{\phi,\ell}(r)-\lambda h(r)\right)^{-1/2}h(r)dr,
\ee
with the bound:
\be
\label{der1bound}
\left|\frac{\pa a_{\phi+\lambda h}}{\pa \lambda} (e,\ell)\right|\leq C\frac{M^2}{e^2\sqrt{\ell}},\quad \forall (\lambda,e)\in \Omega(\phi,\widetilde \phi,\ell),
\ee
for some universal constant $C>0$.\\
(ii) Let $s\in \RR_+^*$. Then the function $\lambda \mapsto a_{\phi+\lambda h}^{-1}(s,\ell)$ is differentiable on $[0,1]$ and we have
\be
\label{der2}
\frac{\pa }{\pa \lambda} a_{\phi+\lambda h}^{-1}(s,\ell)=\frac{\ds \int_{r_1}^{r_2} \left(a_{\phi+\lambda h}^{-1}(s,\ell)- \psi_{\phi,\ell}(r)-\lambda h(r)\right)^{-1/2}h(r)dr}{\ds \int_{r_1}^{r_2} \left(a_{\phi+\lambda h}^{-1}(s,\ell)- \psi_{\phi,\ell}(r)-\lambda h(r)\right)^{-1/2}dr},
\ee
where $(r_i)_{i=1,2}$ shortly denotes $r_i(\phi+\lambda h,a_{\phi+\lambda h}^{-1}(s,\ell),\ell)$.
\end{lemma}
\begin{proof} Recall from Lemma \ref{contaphi} that the functions $e_{\phi+\lambda h,\ell}$, $r_1(\phi+\lambda h,e,\ell)$ and $r_2(\phi+\lambda h,e,\ell)$ are continuous functions of $\lambda$ (for fixed $e$ and $\ell$).

\bs
{\ni \em Step 1. Proof of (i).}

\ms
\ni
This proof of {\em (i)} will be done with Lemma \ref{lebesgue}, exactly in the same manner as the regularity of $a_{\phi+\lambda h}(e,\ell)$ with respect to $e$ in Lemma \ref{Lemmarearrangement}. We fix $\ell>0$ and introduce the following function
$$g(\lambda,e,r)=8\pi^2 \sqrt{2} \left(e- \psi_{\phi,\ell}(r)-\lambda h(r)\right)_+^{1/2},$$
so that
$$
a_{\phi+\lambda h}(e,\ell)=\int_{r_1(\phi+\lambda h,e,\ell)}^{r_2(\phi+\lambda h,e,\ell)} g(\lambda,e,r)dr.
$$
By \fref{minphi} and \fref{majr2-simple}, we have the following uniform bound:
$$g(\lambda,e,r)\leq C\left(\frac{M}{r_1(\phi+\lambda h,e,\ell)}\right)^{1/2}\leq C\frac{M}{\sqrt{\ell}},$$
where $M$ is defined by \fref{not4}.
Hence, one deduces from standard dominated convergence that $(\lambda,e)\mapsto a_{\phi+\lambda h}(e,\ell)$ is a $\calC^0$ function on $[0,1]\times \RR_-$ and satisfies
$$a_{\phi+\lambda h}(e,\ell)>0\Leftrightarrow (\lambda,e)\in \Omega(\phi,\widetilde \phi,\ell).$$

Let us now prove the differentiability of $a_{\phi+\lambda h}(e,\ell)$ with respect to $\lambda$. Let $\lambda_0\in [0,1]$, $e_{0}=e_{\phi+\lambda_0 h, \ell}$,  and $e \in ]e_{0},0[$ be fixed. From the continuity of $e_{\phi+\lambda h,\ell}$ with respect to $\lambda$, we have $e \in ]e_{\phi+\lambda h,\ell},0[$ for $\lambda$ in a neighborhood $I_0$ of $\lambda_0$. Hence, for $\lambda\in I_0$, the  distributional partial derivative of $g$ is given by
$$\frac{\pa g}{\pa \lambda}(\lambda,e,r)=-4 \pi^2\sqrt{2}  \left(e- \psi_{\phi,\ell}(r)-\lambda h(r)\right)^{-1/2}\un_{r_1(\phi+\lambda h,e,\ell)<r<r_2(\phi+\lambda h,e,\ell)}\,h(r).$$
Moreover, from the continuity of $r_1(\phi+\lambda h,e,\ell)$ and $r_2(\phi+\lambda h,e,\ell)$ with respect to $\lambda$, we get that
\be
\label{cvpp2}
\frac{\pa g}{\pa \lambda}(\lambda,e,r)\to \frac{\pa g}{\pa \lambda}(\lambda_0,e,r)\quad \mbox{as }\lambda\to \lambda_0,
\ee
for all  $r\neq r_1(\phi+\lambda_0 h,e,\ell)$, $r\neq r_2(\phi+\lambda_0 h,e,\ell)$. Now, we use \fref{clvm} and \fref{majr2-simple}:
\be
0\leq \frac{\pa g}{\pa \lambda}(\lambda,e,r)\leq C\frac{M^2}{e^2\sqrt{\ell}}\frac{1}{\sqrt{(r-r_1)(r_2-r)}}\un_{r_1<r<r_2}=:q_{\lambda,e}(r),\label{majq2}
\ee
with
\be
\label{pi2}
\int_{r_1}^{r_2}q_{\lambda,e}(r) dr=C\frac{M^2}{e^2\sqrt{\ell}}\pi.
\ee
As in Step 2 of the proof of Lemma \ref{Lemmarearrangement}, one deduces from \fref{cvpp2}, \fref{majq2}, \fref{pi2} and from the Br\'ezis-Lieb theorem that Assumptions (i) and (ii) of Lemma \ref{lebesgue} are satisfied. Hence the function $a_{\phi+\lambda h}(e,\ell)$ is differentiable with respect to $\lambda$ and its differential $\frac{\pa}{\pa \lambda }a_{\phi+\lambda h}(e,\ell)$ is given by \fref{der1}.

We now claim that $\frac{\pa}{\pa \lambda }a_{\phi+\lambda h}(e,\ell)$ is a continuous function of $(\lambda,e)$ on $\Omega(\phi,\widetilde \phi,\ell)$. Indeed, let $\lambda_0\in [0,1]$ and $e_0= e_{\phi+\lambda_0 h,\ell}$ be fixed. A direct adaptation of the proof of \fref{co2}  enables to show that
$$r_1(\phi+\lambda h,e,\ell)\to r_1(\phi+\lambda_0 h,e_0,\ell),\quad r_2(\phi+\lambda h,e,\ell)\to r_2(\phi+\lambda_0 h,e_0,\ell)$$
as $(\lambda,e)\to (\lambda_0,e_0)$.
Thus we have the following a.e. convergence:
\be
\label{cvpp3}
\frac{\pa g}{\pa \lambda}(\lambda,e,r)\to \frac{\pa g}{\pa \lambda}(\lambda_0,e_0,r)\quad \mbox{ as }(\lambda,e)\to (\lambda_0,e_0),\,\,\mbox{ for all }r\neq r_2(\phi+\lambda_0 h,e_0,\ell).
\ee
Hence, using again the domination \fref{majq2} and the fact that $q_{\lambda,e}\to q_{\lambda_0,e_0}$ in $L^1(\RR_+)$ as $(\lambda,e)\to (\lambda_0,e_0)$, we deduce from dominated convergence theorem that
$$
\int_0^{+\infty} \frac{\pa g}{\pa \lambda}(\lambda,e,r)dr \to \int_0^{+\infty} \frac{\pa g}{\pa \lambda}(\lambda_0,e_0,r)dr \quad \mbox{ as }(\lambda,e)\to (\lambda_0,e_0).
$$
In other words, $\frac{\pa}{\pa \lambda }a_{\phi+\lambda h}(e,\ell)$ is a continuous function of $(\lambda,e)$. Similarly,  $\frac{\pa}{\pa e}a_{\phi+\lambda h}(e,\ell)$ is a continuous function of $(\lambda,e)$. Therefore, the function $(\lambda,e)\mapsto a_{\phi+\lambda h}(e,\ell)$ is $\calC^1$ on $\Omega(\phi,\widetilde \phi,\ell)$. Since the bound \fref{der1bound} stems directly from \fref{majq2} and \fref{pi2}, the proof of Item {\em (i)} of Lemma \ref{cont-aphi-l-phi} is complete.

\bs
{\ni \em Step 2. Differentiability of $a^{-1}_{\phi+\lambda h}(s,\ell)$.}

\ms
\ni
Let $(s,\ell)\in (\RR_+^*)^2$.  From Lemma \ref{contaphi}, we already know that the function $\lambda \mapsto a^{-1}_{\phi+\lambda h}(s,\ell)$ is continuous. Let $\lambda_0\in[0,1]$ and consider a sequence $\lambda_n\in [0,1]$ such that $\lambda_n \to \lambda_0$ as $n\to +\infty$. We write
\be
\label{aqw}
\frac{a^{-1}_{\phi+\lambda_n h}(s,\ell)-a^{-1}_{\phi+\lambda_0 h}(s,\ell)}{\lambda_n-\lambda_0} =A_1(\lambda_n)\,A_2(\lambda_n),
\ee
where we have set
$$A_1(\lambda_n)=\frac{a^{-1}_{\phi+\lambda_n h}(s,\ell)-a^{-1}_{\phi+\lambda_0 h}(s,\ell)}{a_{\phi+\lambda_0 h}(a^{-1}_{\phi+\lambda_n h}(s,\ell),\ell)-a_{\phi+\lambda_0 h}(a^{-1}_{\phi+\lambda_0 h}(s,\ell),\ell)}$$
and
\bee
A_2(\lambda_n)&=&\frac{a_{\phi+\lambda_0 h}(a^{-1}_{\phi+\lambda_n h}(s,\ell),\ell)-a_{\phi+\lambda_0 h}(a^{-1}_{\phi+\lambda_0 h}(s,\ell),\ell)}{\lambda_n-\lambda_0}\\
&=&\frac{a_{\phi+\lambda_0 h}(a^{-1}_{\phi+\lambda_n h}(s,\ell),\ell)-a_{\phi+\lambda_n h}(a^{-1}_{\phi+\lambda_n h}(s,\ell),\ell)}{\lambda_n-\lambda_0}
\eee
and where we simply used that $a_{\phi+\lambda_0 h}(a^{-1}_{\phi+\lambda_0 h}(s,\ell),\ell)=s=a_{\phi+\lambda_n h}(a^{-1}_{\phi+\lambda_n h}(s,\ell),\ell)$. Let us examine separately the convergence of the two factors $A_1$ and $A_2$ in \fref{aqw}. From the continuity of $\lambda \mapsto a^{-1}_{\phi+\lambda h}(s,\ell)$, we have:
\be
\label{convpoint}
\lim_{n\to +\infty}a^{-1}_{\phi+\lambda_n h}(s,\ell)=a^{-1}_{\phi+\lambda_0 h}(s,\ell),
\ee
hence
\bea
\nonumber \lim_{n\to +\infty}A_1(\lambda_n)&=&\frac{1}{\ds \frac{\pa a_{\phi+\lambda_0 h}}{\pa e}(a_{\phi+\lambda_0 h}^{-1}(s,\ell),\ell)}\\
\label{aqw1}&=&\frac{1}{\ds 4\pi^2\sqrt{2}\int_{r_1}^{r_2} \left(a_{\phi+\lambda_0 h}^{-1}(s,\ell)- \psi_{\phi,\ell}(r)-\lambda_0 h(r)\right)^{-1/2}dr}.
\eea
\sloppy
Let us now examine the convergence of the term $A_2(\lambda_n)$, that we rewrite as follows:
$$A_2(\lambda_n)=-\frac{1}{\lambda_n-\lambda_0}\int_{\lambda_0}^{\lambda_n}\frac{\pa A}{\pa \lambda}(\mu,e_n)d\mu,$$
where we have denoted
$$A(\lambda,e)= a_{\phi+\lambda h}(e,\ell),\qquad e_n=a^{-1}_{\phi+\lambda_n h}(s,\ell).$$
By \fref{convpoint}, we have $e_n\to e_0=a^{-1}_{\phi+\lambda_0 h}(s,\ell)$. Moreover, from Step 1, we know that the function $(\lambda,e)\mapsto \frac{\pa A}{\pa \lambda}(\lambda,e)$ is continuous at $(\lambda_0,e_0)$. Hence, we have
$$\lim_{\mu \to \lambda_0,\,n\to \infty}\frac{\pa A}{\pa \lambda}(\mu,e_n)=\frac{\pa A}{\pa \lambda}(\lambda_0,e_0)=\left(\frac{\pa a_{\phi+\lambda h}}{\pa \lambda}\right)_{\lambda=\lambda_0}(a^{-1}_{\phi+\lambda_0 h}(s,\ell),\ell),$$
from which we deduce that
\bea
\nonumber \lim_{n\to +\infty}A_2(\lambda_n)&=&-\left(\frac{\pa a_{\phi+\lambda h}}{\pa \lambda}\right)_{\lambda=\lambda_0}(a^{-1}_{\phi+\lambda_0 h}(s,\ell),\ell)\\
\nonumber \label{aqw2}&=&4\pi^2\sqrt{2} \int_{r_1}^{r_2} \left(a_{\phi+\lambda_0 h}^{-1}(s,\ell)-\psi_{\phi,\ell}(r)-\lambda_0 h(r)\right)^{-1/2}h(r)dr.
\eea
where we used \fref{der1}. Finally, \fref{aqw1} and \fref{aqw2} give \fref{der2}. This concludes the proof of Lemma \ref{cont-aphi-l-phi}.
\end{proof}


\subsection{Rearrangement with respect to a given microscopic energy}
\label{rearrangement}


In this section, we introduce the Schwarz symmetrization of a function $f$ with respect to a given microscopic energy $e=\frac{|v|^2}{2}+\phi(x)$ at given momentum $\ell>0$.

We start by defining  a suitable rearrangement of $f\in \calE_{rad}$ at given momentum $\ell>0$  which preserves the generalized Casimir functionals \fref{loi2}. We proceed similarly like for the usual 
Schwarz symmetrization, see \cite{lieb-loss,kavian,ATL}. Let us recall the definition \fref {defmuf}, \fref{defmuf-rul} of the distribution function of $f$ at given $\ell>0$: 
\bee
\mu_f(s, \ell) & = &  \nu_\ell\left\{(r,u)\in \Omega_\ell,\ \ f(r,u,\ell)>s \right\}\\
& = & 4 \pi^2   \int_{r=0}^{+\infty} \int_{u=-\infty}^{+\infty}\un_{f(r,u,\ell)>s}  (r^2u^2 - \ell )^{-1/2} r |u|  \un_{r^2u^2>\ell}\, dr du \,.
\eee
We then have the following elementary lemma:

\begin{lemma}[Properties of $\mu_f$]
\label{defprop}
 Let  $f\in L^1\cap L^{\infty}(\RR^6)$, nonnegative and  spherically symmetric, and let $\mu_f(s,\ell)$ be the distribution function of $f$ at given $\ell$ as defined by \fref{defmuf-rul}. Then there exists a set $A$ with $|A|_{\RR^+}=0$ such that 
 \be
 \label{hoeohe}
 \forall \ell\in \RR^+\backslash A, \ \ \forall s>0, \ \ \mu_f(s,\ell)<+\infty,
 \ee
 \be
 \label{estlinftiy}
\forall \ell\in \RR^+\backslash A, \ \  \forall s\geq |f|_{L^{\infty}}, \ \ \mu_f(s,\ell)=0.
\ee
 Moreover, $\forall \ell\in \RR^+\backslash A$, the map $s\to \mu_f(s,\ell)$ is right continuous on $\RR^*_+$. 
 \end{lemma}
 
 We may now introduce the generalized Schwarz symmetrization:
 
\begin{proposition}[Schwarz symmetrization at fixed $\ell>0$]
\label{defpropschwartz}
 Let  $f\in L^1\cap L^{\infty}(\RR^6)$, nonnegative and  spherically symmetric, let $\mu_f(t,\ell)$ given by \fref{defmuf-rul} and let $A$ be the zero measure set given by Lemma \pref{defprop}.  We define the  Schwarz symmetrization $f^*(\cdot, \ell)$ of $f$ at fixed $\ell$ as being the pseudo inverse of $\mu_f (\cdot , \ell)$:
 \be
 \label{f*1}
 \forall t\geq 0, \ \   \forall  \ell\in \RR^*_+\backslash A,  \ \ f^*(t,\ell)=\left\{\begin{array}{ll}\sup\{s\geq 0: \ \ \mu_f(s,\ell)>t\} \ \ \mbox{for} \ \ t<\mu_f(0,\ell)\\[2mm]
 0 \ \ \mbox{for} \ \ t\geq \mu_f(0,\ell)\\
 \end{array}
 \right .
 \ee
with $\mu_f$ given by \fref{defmuf-rul}. 
Then $f^*(\cdot, \ell)$ is a  nonincreasing function on $[0,\infty)$ and
\be
\label{muf*}
\forall t\geq 0, \ \  \forall  \ell\in \RR^*_+\backslash A, \ \ \mu_f(t,\ell)=\left|\{s>0; f^*(s,\ell)>t\right\}|_{\RR^+}
\ee
 In particular
\be
\label{f*Lp}
|f^*|_{L^p(\RR_+\times \RR_+)}= |f|_{L^p(\RR^6)}, \ \ \forall \ p \in [1,+\infty]\,.
\ee
Moreover, there holds the contractivity relation: 
\be
\label{contra}
|f^*-g^*|_{L^1}\leq|f-g|_{L^1}.
\ee
\end{proposition}

Lemma \ref{defprop} and Proposition \ref{defpropschwartz} can be derived from standard arguments by adapting for example the arguments in \cite{Mossino}, this is left to the reader.\\

Given $f\in \mathcal E_{rad}$ and $\phi\in \Phi_{rad}$, we now define the rearrangement of $f$ with respect to the microscopic energy $\frac{|v|^2}{2}+\phi(x)$.
\sloppy
\begin{proposition}[Symmetric rearrangement with respect to a given microscopic energy]
\label{Lemmarearrangement}  \mbox{} \\
Let  $f\in\calE_{rad}$ and $\phi\in\Phi_{rad}$ non zero. Let $f^*$  be its symmetric rearrangement defined by \fref{f*1}. We define the rearrangement  $f^{*\phi}$ of $f$ with 
respect to the microscopic energy $\frac{|v|^2}{2}+\phi(x)$  by:
\be
\label{symgenphi}
f^{*\phi}(x,v)=f^*\left(a_{\phi}\left(\frac{|v|^2}{2}+\phi(x), |x \times v|^2\right), |x \times v|^2\right)\un_{\frac{|v|^2}{2}+\phi(x)<0}\,,
\ee
where $a_\phi$ is defined in Lemma \pref{lemmadefaphi}.
Then 
\be
\label{equimeasur}
f^{*\phi}\in  \Eq(f)
\ee
where $\Eq(f) $ is defined  by \fref{defeqf}, and in particular:
\be
\label{fphiLp}
|f^{*\phi}|_{L^p}= |f|_{L^p}, \ \ \forall p \in [1,+\infty]\,.
\ee
Moreover, $f^{*\phi}\in \mathcal E_{rad}$ with estimate:
\be
\label{f11}
||v|^2f^{*\phi}|_{L^1}\leq C|\na \phi|_{L^2}^{4/3}\,|f|_{L^1}^{7/9}\,|f|_{L^\infty}^{2/9}\,.
\ee
\end{proposition}

Before proving this proposition, we give a corollary which says that nonincreasing steady states are invariant through the above rearrangement.

\begin{corollary}[Identification of $Q^{*\phi_Q}$]
\label{coroll}
Let $Q$ satisfy Assumption (A) and let $Q^{*\phi_Q}$ be defined according to \fref{symgenphi}. Then $Q^{*\phi_Q}=Q$ and we have
\be
\label{qstarphiq}
Q^*\left(a_{\phi_Q}(e,\ell), \ell\right)=F(e,\ell), \qquad \forall \ell>0 ,\quad  \forall e\in [e_{\phi_Q,\ell},0[,
\ee
where $e_{\phi_Q,\ell}$ is defined in Lemma \pref{lemmapotentials}. In particular, for all $\ell>0$, $Q^*(\cdot,\ell)$ is  a $\calC^1$ function 
on $]0,\mu_{Q}(0,\ell)[$, where $\mu_{Q}$ is defined by \fref{defmuf}.
\end{corollary}

\bs
\ni
{\em Proof of Corollary \pref{coroll}.}
\ms
\ni
Let $\ell > 0$ be fixed and recall the function $F$ defined in Assumption (A). Assume first that $F(e,\ell)=0$ for all $e <0$.  From definition \fref{defmuf}  we have
$\mu_Q(s,\ell)=\nu_\ell\left\{(r,u)\,:\,\,F\left(\frac{|u|^2}{2}+\phi_Q(r),\ell\right)>s\right\}=0$ for all $s\geq 0$. This implies  
  from \fref{f*1} that  $Q^*(\cdot,\ell)=0,$ and then identity \fref{qstarphiq} is satisfied.

Assume now that  $F(\cdot,\ell)$ is not zero on $\RR_{-}^*$ and let
\be
\label{e0(ell)}
e_0(\ell)=\sup\left\{e<0\,:\,\,F(e,\ell)>0\right\}.
\ee
By Assumption (A), we have $e_0(\ell)\leq e_0<0$ and the function $e\mapsto F(e,\ell)$ is continuous, strictly decreasing on $]-\infty,e_0(\ell)]$ and vanishes for $e\geq e_0(\ell)$.
As $F$ is nonnegative, we have from \fref{defmuf}:
$$\mu_Q(F(e,\ell),\ell)=\nu_\ell\left\{(r,u)\,:\,\,F\left(\frac{|u|^2}{2}+\phi_Q(r),\ell\right)>F(e,\ell)\right\}, \quad \forall e\in \RR,$$
and, $F(\cdot,\ell)$ being strictly decreasing on $]-\infty,e_0(\ell)]$, this identity implies
\bea
\mu_Q(F(e,\ell),\ell)&=& \nu_\ell\left\{(r,u)\,:\,\,\frac{|u|^2}{2}+\phi_Q(r)<e\right\}, \quad  \forall e\leq e_{0}(\ell)\nonumber\\
&=&a_{\phi_{Q}}(e,\ell), \quad  \forall e\leq e_{0}(\ell).\label{muQF}
\eea
Assume  that $\mu_Q(0,\ell)=0$,  then $\mu_Q(\cdot ,\ell)=0$ since it is a nonincreasing function.  Hence,  from definition \fref{f*1}  we get 
$Q^*(\cdot,\ell)=0$. Now, we write \fref{muQF} for $e=e_{0}(\ell)$ and deduce from the structure of $a_{\phi_{Q}}$ that $e_{0}(\ell) \leq e_{\phi_{Q},\ell}.$ This means that $F(e,\ell)=0$ for $e\in  [e_{\phi_Q,\ell},0[$, and identity \fref{qstarphiq} is satisfied.

We now assume $\mu_Q(0,\ell)>0$, which implies from \fref{muQF}  that $e_{0}(\ell) > e_{\phi_{Q},\ell}$. We know that $a_{\phi_{Q}}(\cdot,\ell)$ (resp. $F(\cdot,\ell)$) is continuous and one-to-one from  $[e_{\phi_Q,\ell},e_{0}(\ell)]$ to  $[0,a_{\phi_{Q}}(e_{0}(\ell),\ell)]$ (resp. $[0,F(e_{\phi_Q,\ell},\ell)]$). Hence, identity \fref{muQF} ensures that $\mu_{Q}(\cdot,\ell)$ is  invertible from $[0, F(e_{\phi_{Q},\ell},\ell)]$ to $[0,a_{\phi_{Q}}(e_{0}(\ell),\ell)]$ and $Q^*$ (which is by definition its pseudoinverse)  is its inverse in this case.  Therefore, \fref{muQF} implies
$$Q^*\left( a_{\phi_{Q}}(e,\ell), \ell\right) = F(e,\ell), \qquad  \forall e\in  [e_{\phi_Q,\ell},e_{0}(\ell)].$$
Now \fref{muQF}  implies that $a_{\phi_{Q}}(e,\ell)\geq a_{\phi_{Q}}(e_{0}(\ell),\ell)=\mu_Q(0,\ell)$ for $e\in  [e_{0}(\ell),0 [$, which together with the definition of $Q^*$ ensure
that both terms in \fref{qstarphiq} vanish for $e\in  [e_{0}(\ell),0 [$. This ends the proof of \fref{qstarphiq}. Finally, using   \fref{qstarphiq}, we conclude that the stated 
$\calC^1$ regularity of $Q^*$ on $]0,a_{\phi_{Q}}(e_{0}(\ell),\ell)[$   is an immediate consequence of the $\calC^1$ regularity and the non vanishing derivatives of $F$ and $a_{\phi_{Q}}$  
on $]e_{\phi_Q,\ell},e_{0}(\ell)[$. 

To end the proof of Corollary \ref{coroll} , it remains to identify $Q$ and $Q^{*\phi_Q}$ for a.e. $x,v$. Let $(x,v)\in \RR^6$ such that $\ell =|x\times v|^2>0$ and let $e(x,v)=\frac{|v|^2}{2}+\phi_Q(r)\geq \psi_{\phi_Q,\ell}(r)\geq e_{\phi_Q,\ell}$, where we used that $|v|^2\geq \frac{\ell}{r^2}$.
If $e(x,v)<0$, then \fref{qstarphiq} gives directly $Q(x,v)=F(e(x,v),\ell)=Q^{*\phi_Q}(x,v)$, by Assumption (A) and \fref{symgenphi}.
If $e(x,v)\geq 0$, then  we have $Q(x,v)=F(e(x,v),\ell)=Q^{*\phi_Q}(x,v)=0$, using again \fref{symgenphi}. 
This concludes the proof of Corollary \ref{coroll}.
\qed

\bs
\ni
{\em Proof of Proposition \pref{Lemmarearrangement}}

\ms
\ni
We first notice that the formula \fref{symgenphi} is well-defined for a.e. $(x,v)\in \RR^6$ by Proposition \ref{defpropschwartz}. Indeed, from \fref{changevariables} we have that
$$\left|\left\{(x,v)\in\RR^6\,:\,\,|x\times v|^2\in A\right\}\right|_{\RR^6}=0,$$
where $A$ is the measure zero exceptional set given in Lemma \ref{defprop}.

\bs
\ni
{\it Step 1.} The change of variables formula.

\ms
\ni
The equimeasurability of $f$ and $f^{*\phi}$ relies on the following elementary change of variables formula: let two nonnegative functions $\alpha\in {\calC}^0(\RR)\cap L^{\infty}(\RR)$, $\beta\in L^1(\RR^+\times \RR^+)$, then $\forall \ell >0$,
\be
\label{change2}
\begin{array}{ll}
\ds 2\int_0^{+\infty}  \int_0^{+\infty} \alpha\left(\frac{u^2}{2}+\phi(r)\right ) & \ds  \beta\left(a_\phi\left(\frac{u^2}{2}+\phi(r), \ell \right),\ell\right)\un_{\frac{u^2}{2}+\phi(r)<0}d\nu_{\ell} \\ & \ds =\int_0^{+\infty} \alpha(a_{\phi}^{-1}(s,\ell))\beta(s,\ell)ds.
\end{array}
\ee
where $\nu_{\ell}$ is given by \fref{nul}.
This implies in particular from \fref{changevariables}:
\be
\label{change2bis}
\begin{array}{ll}
\ds \int_{\RR^6}\alpha\left(\frac{|v|^2}{2}+\phi(x)\right) & \ds \beta\left(a_\phi\left(\frac{|v|^2}{2}+\phi(x), |x\times v|^2\right),|x\times v|^2\right)dxdv \\ &  \ds =\int_0^{+\infty} \int_0^{+\infty}\alpha(a_{\phi}^{-1}(s,\ell))\beta(s,\ell) dsd\ell.
\end{array}
\ee
Let us prove \fref{change2}.  We first perform the change of variable on the integral on $u$ in the lhs  of \fref{change2},  $e=\frac{u^2}{2}+\phi(r)$, to obtain
\bee
& & 2\int_{r,u \geq 0} \alpha\left(\frac{u^2}{2}+\phi(r)\right )  \ds  \beta\left(a_\phi\left(\frac{u^2}{2}+\phi(r), \ell \right),\ell\right) \un_{\frac{u^2}{2}+\phi(r)<0}d\nu_{\ell}\\
& = & 4\pi^2 \sqrt{2} \int_{r=0}^{+\infty}\int _{e=-\infty}^{0}\alpha(e )  \beta\left(a_\phi(e, \ell),\ell\right)\left(e-\phi(r) - \frac{\ell}{2r^2} \right)^{-1/2}\un_{e-\phi(r) - \frac{\ell}{2r^2}>0}de dr ,
\eee
Now from Fubini and \fref{aphi'}, we get
\bee
&&2\int_{r,u \geq 0} \alpha\left(\frac{u^2}{2}+\phi(r)\right )  \ds  \beta\left(a_\phi\left(\frac{u^2}{2}+\phi(r), \ell \right),\ell\right)\un_{\frac{u^2}{2}+\phi(r)<0}d\nu_{\ell} \\
&&\qquad \qquad = \int _{e_{\phi,\ell}}^{0}\alpha(e )  \beta\left(a_\phi(e, \ell),\ell\right)\frac{\partial a_\phi}{\partial e}(e,\ell) de .
\eee 
From Lemma \ref{lemmadefaphi}, for any $\ell >0$,  the map $e\mapsto a_\phi(e, \ell)$  is a  $\calC^1$-diffeomorphism  from $]e_{\phi,\ell}, 0[$ to $] 0,+\infty[$, and we may therefore perform the change of variable 
 $s=a_\phi(e, \ell)$, which together with \fref{limaphi} yields \fref{change2}.

\bs
\ni
{\em Step 2. Equimeasurability and proof of \fref{f11}.} 

\ms
\ni 
Let now $\phi\in \Phi_{rad}$,   $f\in\calE_{rad}$ and $f^{*\phi}$ given by (\ref{symgenphi}).  We first prove that $f^{*\phi}\in \Eq(f)$ according to the definition \fref{defeqf}. 
For all $t \geq 0$ and $\ell >0$, we use the definition of $\mu_f(\cdot,\ell)$ \fref{defmuf-rul}, of $f^{*\phi}$  \fref{symgenphi}  and the formula  \fref{change2} with $\alpha=1$ and $\beta(s,\ell)= \un_{f^*(s,\ell)>t} $ to get:
$$\mu_{f^{*\phi}}(t,\ell)=2\int_{0}^{+\infty}\int_{0}^{+\infty} \un_{f^{*\phi}(r,u,\ell)>t} d\nu_{\ell}=  
\int_0^{+\infty} \un_{f^*(s,\ell)>t} ds$$ and hence from \fref{muf*}:
$$\forall t\geq 0, \ \ a.e.\ \ \ell>0, \ \ \mu_{f^{*\phi}}(t,\ell)=\mu_f(t,\ell),$$  which implies the equimeasurability of $f$
and $f^{*\phi}$ according to the definition \fref{defeqf}.\\
It remains  to control the kinetic energy of $f^{*\phi}$ according to \fref{f11}. Indeed:
\bee
 ||v|^2f^{*\phi}|_{L^1}& = & 2\int\left(\frac{|v|^2}{2}+\phi\right)f^{*\phi}(x,v)dxdv+2\int\na \phi(x)\cdot \na \phi_{f^{*\phi}}dx\\
& \leq & 2\int\na \phi(x)\cdot \na \phi_{f^{*\phi}}dx\lesssim |\na \phi|_{L^2}\,||v|^2f^{*\phi}|_{L^1}^{1/4}\,|f^{*\phi}|_{L^1}^{7/12}\,|f^{*\phi}|_{L^\infty}^{1/6}\eee
where we used (\ref{symgenphi}) and the interpolation inequality \fref{interpolation}. This together with a straightforward localization argument concludes the proof of (\ref{f11}).\\
This concludes the proof of Proposition \pref{Lemmarearrangement}.\\
\qed

Let us conclude this section with an elementary Lemma which will be useful in the sequel.

\begin{lemma}[Pseudo inverse of $f^*(a_{\phi}(\cdot,\ell),\ell)$]
\label{lemmapseudoinverse}
Let $f\in\calE_{rad}$  and $\phi \in \Phi_{rad}$ be  given nonzero functions, and let $\ell>0$ such that $f^*(0,\ell)>0$.  The function $e \mapsto f^*(a_\phi(e,\ell),\ell)$ is nonincreasing from
$[e_{\phi,\ell}, 0[$ to $[0,f^*(0,\ell)]$.
  We define its pseudo inverse, which we denote (with abuse of notation)  $s\rightarrow (f^*\circ a_{\phi})^{-1}(s, \ell)$, as follows:
\be
\label{defpseduoinverse}
(f^*\circ a_{\phi})^{-1}(s,\ell)=\sup\{e\in [e_{\phi,\ell}, 0[: \ \ f^*( a_{\phi}(e,\ell),\ell)>s\},
\ee
for all $s\in ]0,f^*(0,\ell)[$.
Then $s\rightarrow (f^*\circ a_{\phi})^{-1}(s, \ell)$ is a nonincreasing function  and $\forall (x,v)\in (\RR^3)^2$ such that $|x\times v|^2=\ell$,  $\forall s\in]0,f^*(0,\ell)[$,
\be
\label{implicationone}
 f^{*\phi}(x,v)>s \Longrightarrow \frac{|v|^2}{2}+\phi (x)\leq (f^*\circ a_{\phi})^{-1}(s,\ell),
\ee
\be
\label{implicationtwo}
 f^{*\phi}(x,v)\leq s \Longrightarrow \frac{|v|^2}{2}+\phi (x)\geq (f^*\circ a_{\phi})^{-1}(s,\ell).
 \ee
\end{lemma}

\begin{proof}  
Let $\ell >0$ and $s\in (0,f^*(0,\ell))$, then $f^*( a_{\phi}(e_{\phi,\ell},\ell),\ell)=f^*(0,\ell)>s$ and hence
$$
\label{cohoeh}
\{e\in [e_{\phi,\ell},0): \ \ f^*( a_{\phi}(e,\ell),\ell)>s\} \ \ \mbox{is not empty}.
$$
This means that $(f^*\circ a_{\phi})^{-1}(s,\ell)$ is well defined  for $s\in (0,f^*(0,\ell))$.  The monotonicity of $(f^*\circ a_{\phi})^{-1}$ follows from the monotonicity of $f^*$ and $a_{\phi}$. 

Let now $(x,v) \in \RR^6$ be  such that $|x\times v|^2=\ell >0$. Assume $f^{*\phi}(x,v)> s$,  then $f^*(a_{\phi}(\frac{|v|^2}{2}+\phi (x), \ell),  \ell)> s$ and thus $\frac{|v|^2}{2}+\phi (x)< 0$ . Thus we have either $\frac{|v|^2}{2}+\phi (x)< e_{\phi,\ell}$, and in this case  (\ref{implicationone}) is trivial, or $\frac{|v|^2}{2}+\phi (x)\in [e_{\phi,\ell},0)$, and this implies  $\frac{|v|^2}{2}+\phi (x)\leq (f^*\circ a_{\phi})^{-1}(s,\ell))$ from the definition (\ref{defpseduoinverse}). Thus (\ref{implicationone}) is  then proved. 
Assume now $f^{*\phi}(x,v)\leq s$. If $\frac{|v|^2}{2}+\phi (x)\geq 0$ then (\ref{implicationtwo}) is trivial, 
otherwise $f^{*\phi}(x,v)\leq s<f^*(0,\ell)$ implies that $\frac{|v|^2}{2}+\phi (x)\in (e_{\phi,\ell},0)$. Thus for all $e\in\{e\in [e_{\phi,\ell},0): \ \ f^*(a_{\phi}(e,\ell),  \ell)>s\}$ which is a non empty set, $\frac{|v|^2}{2}+\phi (x)\geq e$, and (\ref{implicationtwo}) follows.
\end{proof}


\section{Nonlinear stability of the Vlasov Poisson system}
\label{sectionstar}


This section is devoted to the proof of the main results of this paper. We first exhibit the key monotonicity formula involving the generalized symmetric rearrangement  with respect to the Poisson field \fref{symgenphi}, Proposition \ref{propclemonotonie}, which allows us to reduce the study of the minimization problem of Theorem \ref{theo1} to the one of an unconstrained minimization problem on the Poisson field only. The study of this new problem, that is the proof of Proposition \ref{propA}, is postponed to section \ref{sectionj}, and immediately yields Theorem \ref{theo1}. We then show how to extract compactness from minimizing sequences to prove Theorem \ref{theo2} which now implies Theorem \ref{theo3} from standard arguments.


\subsection{The monotonicity formula}


Given $f\in\calE_{rad}$, we will note to ease notation:
\be
\label{deffhat}
\widehat{f}=f^{*\phi_f}, 
\ee
and recall from Proposition \ref{Lemmarearrangement} that:
\be
\label{propftilde}
\widehat{f}\in\calE_{rad} \cap \Eq(f).
\ee
We introduce the functional of $\phi\in \Phi_{rad}$:
\be
\label{defiphi}
 \mathcal J _{f^*}(\phi) =\calH(f^{*\phi})+\frac{1}{2}|\na \phi-\na \phi_{f^{*\phi}}|^2,
\ee
and claim the following monotonicity formula which is a fundamental key for our analysis -see also \cite{aly} for related statements-:

\begin{proposition}[Monotonicity of the Hamiltonian under the $f^{*\phi_f}$ rearrangement]
\label{propclemonotonie}
Let $f\in\calE_{rad}$, non zero, and $\widehat{f}$ given by \fref{deffhat}, then:
\be
\label{keymononicity}
\calH(f)\geq  \mathcal J _{f^*}(\phi_f)\geq \calH(\widehat{f}).
\ee
Moreover, $\calH(f)= \calH(\widehat{f})$ if and only if $f=\widehat{f}$.
\end{proposition}

\begin{proof}

Let $f,g\in \calE_{rad}$, then:
\bee
\calH(f)& = & \frac{1}{2}\int_{\RR^6}|v|^2f-\frac{1}{2}\int_{\RR^3}|\nabla\phi_f|^2\\
& = & \int_{\RR^6}\left(\frac{|v|^2}{2}+\phi_f\right)(f-g)+\frac{1}{2}\int_{\RR^6}|v|^2 g+\int_{\RR^3}\phi_fg+\frac{1}{2}\int|\nabla\phi_f|^2\\
& = & \calH(g)+\int_{\RR^6}\left(\frac{|v|^2}{2}+\phi_f\right)(f-g)+\frac{1}{2}\int|\nabla \phi_g|^2-\int\nabla \phi_f\cdot\nabla\phi_g+\frac{1}{2}\int|\nabla\phi_f|^2
\eee
and hence the general formula: $\forall f,g\in\calE_{rad}$,
\be
\label{identity}
\calH(f)=\calH(g)+\frac{1}{2}|\na\phi_f-\na\phi_{g}|_{L^2}^2+\int_{\RR^6}\left(\frac{|v|^2}{2}+\phi_f(x)\right)(f-g)\,dxdv.
\ee
We apply this formula with $g=\widehat{f}=f^{*\phi_f}$ and rewrite the result using (\ref{defiphi}):
$$\calH(f)=\mathcal J _{f^*}(\phi_f)+\int_{\RR^6}\left(\frac{|v|^2}{2}+\phi_f(x)\right)(f-\widehat{f})\,dxdv.$$ We now claim:
\be
\label{identity2}
 \int_{\RR^6}\left(\frac{|v|^2}{2}+\phi_f(x)\right)(f-\widehat{f})\,dxdv\geq 0,
\ee
with equality if  and only if $f=\widehat{f}$, which immediately implies (\ref{keymononicity}).

The proof of (\ref{identity2}) is reminiscent from the standard inequality for symmetric rearrangement known as the "bathtub" principle $$\int_{\RR^6}|x|f^*\geq\int_{\RR^6}|x|f,$$ see \cite{lieb-loss}. Indeed, use $f(x,v)=\int_{t=0}^{+\infty}\un_{t<f(x,v)} dt$ and Fubini to derive:
\bea
\label{mesureurto}
\nonumber & & \int_{\RR^6}\left(\frac{|v|^2}{2}+\phi_f\right)(f-\widehat f)\,dxdv  =  \int_ {t=0}^{+\infty}dt\int_{\RR^6}\left(\frac{|v|^2}{2}+\phi_f\right)\left(\un_{t<f(x,v)} -\un_{t<\widehat f(x,v)}\right)dxdv\\ 
\nonumber & = &  \int_{t=0}^{+\infty}dt\int_{\RR^6}\left(\frac{|v|^2}{2}+\phi_f\right)\left(\un_{\widehat f(x,v)\leq t<f(x,v)}-\un_{f(x,v)\leq t<\widehat f(x,v)}\right)dxdv\\
& = & \int_{\ell=0}^{\infty}d\ell \int_{t=0}^{f^*(0,\ell)}dt \left(\int_{S_{1,\ell}(t)}  - \int_{S_{2,\ell}(t)}\right)  \left(\frac{u^2}{2}+\phi_f(r)\right)d\nu_{\ell}
\eea
where $d\nu_{\ell}$ is given by \fref{nul}, and
$$S_{1,\ell}(t)=\{(r,u)\in \Omega_\ell,  \widehat f(r,u,\ell)\leq t<f(r,u,\ell)\}, $$
$$ S_{2,\ell}(t)=\{(r,u)\in \Omega_\ell,   f(r,u,\ell)\leq t < \widehat f(r,u,\ell)\},$$
We now use (\ref{implicationone}) in  Lemma \ref{lemmapseudoinverse} to obtain: $\forall t\in (0,f^*(0,\ell))$,
$$\ds  \int_{S_{2,\ell}(t)}\left(\frac{u^2}{2}+\phi_f(r)\right)d\nu_{\ell}\leq (f^{*}\circ a_{\phi})^{-1}(t, \ell)  \nu_{\ell}(S_{2,\ell}(t))
$$
where we recall that
$$\nu_{\ell}(S_{2,\ell}(t))= 4\pi^2  \int_{S_{2,\ell}(t)}\un_{r^2u^2>\ell} (r^2u^2-\ell)^{-1/2} r|u|drdu.$$
We then observe from $\widehat{f}\in \Eq(f)$ that: 
$$\mbox{for a.e.} \ \ t>0, \ \ \nu_{\ell}(S_{1,\ell}(t))=\nu_{\ell}(S_{2,\ell}(t)),
$$ 
and deduce
$$\int_{S_{2,\ell}(t)} \left(\frac{u^2}{2}+\phi_f(r)\right)d\nu_{\ell}  \leq (f^{*}\circ a_{\phi})^{-1}(t, \ell)\int_{S_{1,\ell}(t)}d\nu_{\ell}.
$$
 Injecting this into (\ref{mesureurto}) and using  (\ref{implicationtwo}) yields:
$$\begin{array}{l} \ds \int_{\RR^6}\left(\frac{|v|^2}{2}+\phi_f\right)(f-\widehat f)\,dxdv  \geq \\
\ \ \ \ \ \ \ \ \ \ \ \ \ \ds \int_{\ell=0}^{\infty}d\ell \int_{t=0}^{f^*(0,\ell)}dt \int_{S_{1,\ell}(t)} \left(\frac{u^2}{2}+\phi_f(r)- (f^{*}\circ a_{\phi})^{-1}(t, \ell)\right)d\nu_{\ell}
\geq 0
\end{array}$$
and the analogous inequality for $S_{2,\ell}(t)$:

$$\begin{array}{l} \ds\int_{\RR^6}\left(\frac{|v|^2}{2}+\phi_f\right)(f-\widehat f)\,dxdv  \geq \\ \ds
\ \ \ \ \ \ \ \ \ \ \ \ \ \int_{\ell=0}^{\infty}d\ell \int_{t=0}^{f^*(0,\ell)}dt \int_{S_{2,\ell}(t)} \left(  (f^{*}\circ a_{\phi})^{-1}(t, \ell) -  \frac{u^2}{2}-\phi_f(r)\right)d\nu_{\ell}\geq 0. \end{array}
$$
 Moreover, assume that $ \int_{\RR^6}\left(\frac{|v|^2}{2}+\phi_f(x)\right)(f-\widehat{f})\,dxdv=0$. Recalling that $\nu_\ell(S_{1,\ell}(t))=\nu_\ell(S_{2,\ell}(t))=0$ for $t>f^*(0,\ell)$, the above two chains of equalities imply that for a.e $t,\ell>0$, either $\nu_\ell( S_{1,\ell}(t))=\nu_\ell(S_{2,\ell}(t))=0$ or $\nu_\ell( S_{1,\ell}(t))=\nu_\ell(S_{2,\ell}(t))>0$ with: 
 $$ \frac{u_1^2}{2}+\phi_f(r_1)=(f^{*}\circ a_{\phi})^{-1}(t,\ell)=\frac{u_2^2}{2}+\phi_f(r_2),$$
 for a.e $(r_1,u_1)\in S_{1,\ell}(t)$,  a.e $(r_2,u_2)\in S_{2,\ell}(t)$, which contradicts the fact that $\widehat f(r_1,u_1,\ell)\leq t<\widehat f(r_2,u_2,\ell)$. We conclude that a.e $t,\ell>0$, $\nu_\ell(S_{1,\ell}(t))=\nu_\ell(S_{2,\ell}(t))=0$ which implies $f=\widehat f$. This concludes the proof of (\ref{identity2}) and of Proposition \ref{propclemonotonie}.
\end{proof}


\subsection{Reduction to a variational problem on $\phi$ and proof of Theorem \ref{theo1}}


We now claim the following local coercivity property of the functional of $\phi$ given by (\ref{defiphi}). To ease notations, we let for $\phi\in \Phi_{rad}$:
\be
\label{defiphiq}
  \mathcal J (\phi)= \mathcal J _{Q^*}(\phi) =\calH(Q^{*\phi})+\frac{1}{2}\int_{\RR^3}|\na \phi-\na \phi_{Q^{*\phi}}|^2
\ee

\begin{proposition}[$\phi_Q$ is a local strict minimizer of $\mathcal J$]
\label{propA}
\mbox{}\\
There exist a constant $C_0>0$ such that the following holds. For all $R>0$, there exists $\delta_0(R)\in ]0,\frac{1}{2}|\na \phi_Q|_{L^2}]$ such that, for all $f\in \calE_{rad}$ satisfying
$$|f-Q|_{\calE}\leq R,\qquad |\na\phi_f-\na\phi_Q|_{L^2}\leq \delta_0(R),$$
we have
\be
\label{coerc}
\mathcal J(\phi_f)-\mathcal J(\phi_Q)\geq C_0|\na \phi_f-\na \phi_Q|^2_{L^2}\,.
\ee
\end{proposition}

The proof of this Proposition essentially relies on Antonov's coercivity property and is postponed to section \ref{sectionj}. Theorem \ref{theo1} is now a straightforward 
consequence of Propositions \ref{propclemonotonie} and \ref{propA}.

\bs
\ni
{\em Proof of Theorem \pref{theo1}.}

\ms
\ni
Let $R>0$ and $f\in \calE_{rad}\cap \Eq(Q)$ satisfying \fref{a-theo1}, where $\delta_0(R)$ is as in Proposition \ref{propA}. In particular, note that
$$|\na \phi_f-\na\phi_Q|_{L^2}\leq \frac{1}{2}|\na\phi_Q|_{L^2}$$
implies that $\phi_f\neq 0$ and $f\neq 0$.
 Then the monotonicity property (\ref{keymononicity}), $f^*=Q^*$ and (\ref{defiphi}) yield:
\be
\label{hohfeofh}
\calH(f)-\calH(Q)\geq  \mathcal J _{f^*}(\phi_f)-\calH(Q)= \mathcal J(\phi_f)-\calH(Q).
\ee 
On the other hand, recall from Corollary \ref{coroll} that our assumption on the ground state $Q$ ensures $$\widehat{Q}=Q^{*\phi_Q}=Q \ \ \mbox{and thus} \ \ \calH(Q)=\mathcal J(\phi_Q).$$ Injecting this together with (\ref{coerc}) into (\ref{hohfeofh}) yields: 
\be
\label{oveonoveo}
\calH(f)-\calH(Q)\geq \mathcal J(\phi_f)-\mathcal J(\phi_Q)\geq C_0|\na \phi_f-\na \phi_Q|^2_{L^2},
\ee this is (\ref{minoration}). If in addition $\calH(f)=\calH(Q)$, then $\phi_f=\phi_Q$ and hence using $f^*=Q^*$:
$$\mathcal H(f^{*\phi_f})=\mathcal H(Q^{*\phi_f})=\mathcal H(Q^{*\phi_Q})=\mathcal H(Q)=\calH (f).$$ We thus are in the case of equality of Proposition \ref{propclemonotonie} from which: $$f=f^{*\phi_f}=f^{*\phi_Q}=Q^{*\phi_Q}=Q.$$ This concludes the proof of Theorem \ref{theo1}.


\subsection{Compactness of minimizing sequences}
\label{sectionthmmain}

We are now in position to prove Theorem \ref{theo2}.

\bs
\ni
{\em Proof of Theorem \pref{theo2}.}

\ms
\ni
The key to extract compactness is the monotonicity formula \fref{oveonoveo} which yields a lower bound on the Hamiltonian involving the Poisson field $\phi_f$ only, while standard Sobolev embeddings ensure that $\phi_f$ enjoys nice compactness properties in the radial setting.

\bs
\ni
{\em Step 1. Weak convergence in $L^p$, $p>1$}.

\ms
\ni
Let
\be
\label{RRR}
R=|Q|_{\calE}+C(1+|\na \phi_Q|_{L^2})^{4/3}\,|Q|_{L^1}^{7/9}\,|Q|_{L^\infty}^{2/9}+|Q|_{L^1}+|Q|_{L^\infty},
\ee
where $C$ is the constant in the interpolation inequality \fref{f11}.

Let $f_n\in \calE_{rad}$ be a sequence satisfying \fref{condfn1}, \fref{condfn2}, where $\delta$ will be fixed further, satisfying in particular
\be
\label{DDD}
\delta\leq \min\left(1,\frac{1}{2}|\na\phi_Q|_{L^2}\right).
\ee
Observe that \fref{condfn1} and \fref{DDD} imply $\phi_{f_n}\neq 0$.
The sequence $f_n^*$ is bounded in $L^1$ by \fref{condfn2}, so $f_n$ is itself bounded in $L^1$. Moreover, from $\calH(f_n)<C$, the $L^\infty$ bound of $f_n$ and the interpolation inequality \fref{interpolation}, $|v|^2f_n$ is uniformly bounded in $L^1$. Hence $f_n$ is bounded in the energy space $\calE_{rad}$. We then get:
\be
\label{Lpfaible}
f_n \rightharpoonup f\
\in \calE_{rad} \ \  \mbox{ in }L^p \quad \mbox{for all } 1< p<+\infty,
\ee
up to a subsequence. Moreover, by a standard consequence of interpolation, Sobolev embeddings and elliptic regularity, we have
\be
\label{champ1}
|\na \phi_{f_n}- \na \phi_f|_{L^2}\to 0 \ \ \mbox{and} \ \ |\phi_{f_n} - \phi_{f}|_{L^\infty}\to 0 \ \ \mbox{as} \ \ n\to +\infty.
\ee
 From assumptions (\ref{condfn1}) and \fref{condfn2}:
\be
\label{fhfheoheoh}
|\nabla \phi_f-\nabla \phi_Q|_{L^2}\leq \delta.
\ee
In particular, $\phi_{f}\neq 0$,  since $\delta<|\nabla \phi_Q|_{L^2}$ from \fref{DDD}. Hence, by Proposition \ref{Lemmarearrangement}, we  have
\be
\label{Q*phi-eq-Q}
Q^{*\phi_{f}} \in \Eq(Q).
\ee

\bs
\ni
{\em Step 2.  Strong convergence in $\calE$ of the sequence $Q^{*\phi_{f_n}}$}.

\ms
\ni
We now aim at extracting a preliminary compactness from $f_n$. Let 
\be
\label{tilde}
\widetilde f_n=Q^{*\phi_{f_n}},\qquad \widetilde f=Q^{*\phi_{f}},
\ee
and observe that  $\widetilde f_n$  is in fact a function of $\phi_{f_n}$. We then claim that the strong convergence (\ref{champ1}) automatically implies some strong compactness in $\calE$ for $\widetilde f_n$:
\be
\label{convfntilde}
(1+|v|^2)\widetilde f_n\to (1+|v|^2)\widetilde f \mbox{ in }L^1(\RR^6).
\ee
We claim also that there exists $\delta_1(R)$ such that, for $0< \delta\leq \delta_1(R)$ we have
\be
\label{fhjoigo}
|\nabla \phi_{\widetilde f}-\nabla \phi_Q|_{L^2}\leq \frac{\delta_0(R)}{2},
\ee
where $R$ is defined by \fref{RRR} and $\delta_0(R)$ is defined in Theorem \ref{theo1}. We are now ready to fix the constant $\delta$ of Theorem \ref{theo2} as follows:
$$\delta=\min\left(1,\frac{1}{2}|\na\phi_Q|_{L^2},\delta_1(R)\right).$$

\bs
\ni
{\it Proof of (\ref{convfntilde}), (\ref{fhjoigo})}. We first claim the a.e convergence:
\be
\label{aeconegrence}
\widetilde f_n\to \widetilde f \ \ \mbox{as} \ \ n\to +\infty  \ \ \mbox{for a.e} \ \ (x,v)\in \RR^6.
\ee
Indeed,  let  $(x,v)\in \RR^6$ such that $|x \times v|^2=\ell >0$. If $e=\frac{|v|^2}{2}+\phi_{f}(x)<0$, then from (\ref{champ1}), $\frac{|v|^2}{2}+\phi_{f_{n}}(x)<e/2$ for $n$ large enough  and  
\be
\label{bb}
\left|\frac{|v|^2}{2}+\phi_{f_{n}}(x)\right|=-\frac{|v|^2}{2}-\phi_{f_{n}}(x) \leq -\phi_{f_{n}}(0)\leq C.
\ee 
 We now recall  from Lemma \ref{contaphi} that for all $\ell>0$: 
\be
\label{convaphit}
a_{\phi_{f_n}}(e,\ell)\to a_{\phi_f}(e,\ell),
\ee 
uniformly with respect to $e$ lying in a compact subset of $]-\infty,0[$. Therefore, from $\frac{|v|^2}{2}+\phi_n(x)<e/2<0$ and from \fref{bb},
\bee
a_{\phi_{f_n}}\left(\frac{|v|^2}{2}+\phi_{f_n}(x), |x\times v|^2\right) \to  a_{\phi}\left(\frac{|v|^2}{2}+\phi_{f}(x), |x\times v|^2\right)\mbox{ as }n\to +\infty.
\eee
Since, by Corollary \ref{coroll}, Lemma \ref{lemmadefaphi} and Assumption (A), the function $Q^*(\cdot,\ell)$ is continuous, this implies $Q^{*\phi_{f_n}}(x,v)\to Q^{*\phi}(x,v)$. Similarly, $\frac{|v|^2}{2}+\phi_{f}(x)>0$ implies $\frac{|v|^2}{2}+\phi_n(x)>0$ for $n$ large enough and thus $Q^{*\phi_{f_n}}(x,v)= Q^{*\phi_{f}}(x,v)=0$. Hence $Q^{*\phi_{f_n}}\to  Q^{*\phi_{f}}$  a.e in $\RR^6$ and (\ref{aeconegrence}) is proved.

Now recall from Proposition \ref{Lemmarearrangement} and from $\phi_{f_n}\neq 0$, $\phi_f\neq 0$, that   $\widetilde f_n\in \Eq(Q)$ and $\widetilde f \in \Eq(Q)$ so that 
$$\forall n\geq 1, \ \ \int_{\RR^6} \widetilde f_n=\int_{\RR^6} Q=\int_{\RR^6}\widetilde f.$$ 
 The almost everywhere convergence of  $\widetilde f_n=Q^{*\phi_{f_{n}}}$ to $\widetilde f$ and the  fact that
 $|\widetilde f_n|_{L^1} = |\widetilde f|_{L^1}$ allows us to apply the Br\'ezis-Lieb Lemma (see \cite{lieb-loss}, Theorem 1.9) and get the strong $L^1$ convergence 
\be
\label{hohfohfoe}
\widetilde f_n\to \widetilde f \ \ \mbox{in} \ \ L^1 \ \ \mbox{as} \ \ n\to +\infty.
\ee 
It remains to prove the strong convergence of the kinetic energy. Let us decompose
$$\widetilde f_n=\un_{|v|^2\leq R}\widetilde f_n+\un_{|v|^2> R}\widetilde f_n=g_{n,R}+h_{n,R}.$$
The $L^1$ convergence \fref{hohfohfoe} implies: $\forall R>0$, $|v|^2 g_{n,R}\to |v|^2\un_{|v|^2\leq R}\widetilde f_n $ in $L^1$. Consider the other term. We recall that $\widetilde f_n=Q^{*\phi_{f_n}}$ is supported in the set $\frac{|v|^2}{2}+\phi_{f_n}(x)< 0$. Hence, by interpolation,
$$\begin{array}{ll}
\ds ||v|^2h_{n,R}|_{L^1} =\int |v|^2 h_{n,R}(x,v)\,dxdv&\ds \leq -2\int \phi_{f_n}(x)h_{n,R}(x,v)\,dxdv\\[4mm]
&\ds \lesssim |\na \phi_{f_n}|_{L^2}\,||v|^2h_{n,R}|_{L^1}^{1/4}\,|h_{n,R}|_{L^1}^{7/12}\,|Q|_{L^\infty}^{1/6}\,,
\end{array}
$$
which yields 
$$||v|^2h_{n,R}|_{L^1}\leq C\,|h_{n,R}|_{L^1}^{7/9}.$$
By writing
$$|h_{n,R}|_{L^1}\leq |Q^{*\phi_{f_n}}-Q^{*\phi_f}|_{L^1}+\int_{|v|^2>R}Q^{*\phi_{f}}(x,v)\,dxdv,$$
we obtain that $||v|^2h_{n,R}|_{L^1}$ converges to 0 when $R\to +\infty$ and $n\to +\infty$ independently. This together with the convergence of $|v|^2g_{n,R}$ concludes the proof of (\ref{convfntilde}).\

We now turn to the proof of (\ref{fhjoigo}) and claim that it follows directly from (\ref{fhfheoheoh}) and the definition $\widetilde f=Q^{*\phi_f}$. Indeed, arguing by contradiction, we extract a subsequence $\nabla \phi_{n}\to \nabla \phi_{Q}$ in $L^2$ and $\widetilde g_n=Q^{*\phi_{n}}$ such that $|\nabla \phi_{\widetilde g_n}-\nabla \phi_Q|_{L^2}\geq \frac{\delta_0(R)}{2}$. From (\ref{f11}), $\widetilde g_n$ is a bounded sequence in $\calE_{rad}$ and then the same proof like for (\ref{convfntilde}) yields $$(1+|v|^2)\widetilde g_n \to (1+|v|^2)Q^{*\phi_Q}=(1+|v|^2)Q\ \  \ \mbox{in} \ \ L^1$$ and hence $\nabla \phi_{\widetilde g_n}\to \nabla \phi_{Q}$ in $L^2$, a contradiction. This concludes the proof of (\ref{fhjoigo}).

\bs
\ni
{\em Step 3. Identification of the limit}.

\ms
\ni
Following (\ref{deffhat}), we let:
\be
\label{hat}
\widehat f_n={f_n}^{*\phi_{f_n}}.
\ee
We now claim that the variational characterization of $Q$ given by Theorem \ref{theo1} and the monotonicity of Proposition \ref{propclemonotonie} allow us to identify the limit:
\be
\label{identifiucaitonfirst}
\widetilde{f}=Q \ \ \mbox{and} \ \ \phi_{\tilde{f}}=\phi_f=\phi_Q,
\ee
and to obtain the additional convergence:
\be
\label{keyconvergence}
\int_{\RR^6}\left(\frac{|v|^2}{2}+\phi_{f_n}(x)\right)(f_n-\widehat f_n)\,dxdv \to 0\ \ \mbox{as} \ \ n\to +\infty.
\ee

\bs
\ni
{\it Proof of (\ref{identifiucaitonfirst}), (\ref{keyconvergence})}. First observe from (\ref{convfntilde}), $|\widetilde f_n|_{L^{\infty}}=|Q^*|_{L^{\infty}}$ and (\ref{interpolation}) that:
\be
\label{champ2}
\calH(\widetilde f_n)\to \calH(\widetilde f),\qquad 
\na \phi_{\widetilde f_n}\to \na \phi_{\widetilde f}\quad \mbox{ in }L^2.
\ee
From (\ref{change2bis}), there holds:
\bee
|\widehat f_n-\widetilde f_n|_{L^1}& = & \int_{\RR^6}|f_n^*-Q^*|\left(a_{\phi_{f_n}}\left(\frac{|v|^2}{2}+\phi_{f_n},| x\times v|^2\right),| x\times v|^2\right)\un_{\frac{|v|^2}{2}+\phi_{f_n}(x)<0}dxdv\\
&=&\int_0^{+\infty}\int_0^{+\infty}|f_n^*-Q^*|(s,\ell)\,ds d\ell\\
& \to & 0 \ \mbox{as} \ \ n\to +\infty
\eee
from the assumption \fref{condfn2}. Together with \fref{convfntilde}, this yields:
\be
\label{convfnhat}
\widehat f_n \to \widetilde f \mbox{ in }L^1(\RR^6).
\ee
We now invoke the identity \fref{identity} with $g=\widetilde f_n$ to derive:
\be
\label{maj1}
\begin{array}{l}
  \ds \frac{1}{2}|\na\phi_{f_n}-\na\phi_{\widetilde f_n}|^2+\calH(\widetilde f_n)-\calH(Q)+\int_{\RR^6}\left(\frac{|v|^2}{2}+\phi_{f_n}(x)\right)(f_n-\widehat f_n)\,dxdv \\[4mm]
\ds \qquad \qquad = \calH(f_n)-\calH(Q)+\int_{\RR^6}\left(\frac{|v|^2}{2}+\phi_{f_n}(x)\right)(\widetilde f_n-\widehat f_n)\,dxdv.
\end{array}
\ee
Let us examine the various terms of this identity. From (\ref{champ2}) and (\ref{fhjoigo}), 
\be
\label{klm}
|\nabla \phi_{\widetilde f_n}-\nabla \phi_Q|_{L^2}<\delta_0(R)
\ee
 for $n$ large enough, where $\delta_0(R)$ is defined in Theorem \ref{theo1}. Moreover, from the definition \fref{tilde} and from Proposition \ref{Lemmarearrangement}, we have the following estimates on $\widetilde f_n$ and $\widetilde f$:
$$
 |\widetilde f_n|_{L^1}=|Q|_{L^1},\quad |\widetilde f_n|_{L^\infty}=|Q|_{L^\infty},\quad 
 ||v|^2\widetilde f_n|_{L^1}\leq C|\na \phi_{f_n}|_{L^2}^{4/3}\,|Q|_{L^1}^{7/9}\,|Q|_{L^\infty}^{2/9},
 $$
 $$
 |\widetilde f|_{L^1}=|Q|_{L^1},\quad |\widetilde f|_{L^\infty}=|Q|_{L^\infty},\quad 
 ||v|^2\widetilde f|_{L^1}\leq C|\na \phi_{f}|_{L^2}^{4/3}\,|Q|_{L^1}^{7/9}\,|Q|_{L^\infty}^{2/9},
 $$
where $C$ is the constant in the interpolation inequality \fref{f11}. Since, by \fref{DDD} and \fref{condfn1}, we have
$$|\na \phi_{f_n}|_{L^2}\leq 1+ |\na \phi_Q|_{L^2},$$
we deduce from \fref{champ1} and \fref{RRR} that
\be
\label{klm2}
|\widetilde f_n-Q|_{\calE}\leq |Q|_{\calE}+|\widetilde f_n|_{\calE}\leq R,\qquad |\widetilde f-Q|_{\calE}\leq |Q|_{\calE}+|\widetilde f|_{\calE}\leq R.
\ee
Therefore, from \fref{klm}, \fref{klm2} and $\widetilde f_n\in \Eq(Q)$, the variational characterization of $Q$ given by Theorem \ref{theo1} ensures:
$$
\calH(\widetilde f_n)-\calH(Q)\geq 0.
$$
Next, from (\ref{identity2}):
$$\int_{\RR^6}\left(\frac{|v|^2}{2}+\phi_{f_n}(x)\right)(f_n-\widehat f_n)\,dxdv\geq 0.$$
We now claim that:
\be
\label{diffftilfhat}
\int_{\RR^6}\left(\frac{|v|^2}{2}+\phi_{f_n}\right)(\widetilde f_n-\widehat f_n)\,dxdv\to 0 \mbox{ as }n\to +\infty.
\ee
Indeed, from (\ref{change2bis}):
$$\begin{array}{l}
\ds \int_{\RR^6}\left(\frac{|v|^2}{2}+\phi_{f_n}(x)\right)(\widetilde f_n-\widehat f_n)\,dxdv=\\[4mm]
\ds \qquad =\int_{\RR^6}\left(\frac{|v|^2}{2}+\phi_{f_n}(x)\right)(Q^*-f_n^*)\left(a_{\phi_{f_n}}\left(\frac{|v|^2}{2}+\phi_{f_n}(x), |x\times v|^2\right), |x\times v|^2\right)\\[4mm]
\ds \qquad \qquad \qquad \qquad \qquad \qquad\qquad \qquad \qquad\qquad  \qquad \qquad\qquad\times \un_{\frac{|v|^2}{2}+\phi_{f_n}<0}\,dxdv\\[4mm]
\ds \qquad =\int_0^{+\infty}a_{\phi_{f_n},\ell}^{-1}(s)(Q^*-f_n^*)(s,\ell)\,ds d\ell,
\end{array}
$$
where we recall that $a_{\phi_{f_n},\ell}^{-1}$ is the inverse of the diffeomorphism $e\to a_\phi(e,\ell)$.
From Lemma \ref{lemmadefaphi},  \fref{mine2}, and  (\ref{champ1}), we have $|a_{\phi_{f_n},\ell}^{-1}(s)|\leq -e_{\phi_{f_n},\ell}\leq -\phi_{f_n}(0)\leq C$, and hence the $L^1$ convergence of $f_n^*$ to $Q^*$ yields (\ref{diffftilfhat}).

Finally, since \fref{condfn2} gives $\underset{\footnotesize n\to +\infty}{\mbox{\rm lim sup}}\,\calH(f_n)\leq \calH(Q)$, we deduce that all the nonnegative quantities in the left-hand side of \fref{maj1} converge to 0:
\be
\label{champ3}
\na \phi_{f_n}- \na \phi_{\widetilde f_n}\to 0\mbox{ in }L^2,
\ee
\be
\label{champ4}
\calH(\widetilde f_n)\to \calH(Q),
\ee
and (\ref{keyconvergence}) holds, and moreover:
\be
\label{necessarliy}
\mathcal H(f_n)\to \mathcal H(Q).
\ee
Hence, \fref{champ2} and \fref{champ4} imply $\calH(\widetilde f)=\calH(Q)$ and thus (\ref{fhjoigo}), \fref{klm2} and Theorem \ref{theo1} ensure:
$$\widetilde f=Q\quad \mbox{and}\quad \phi_{\widetilde f}=\phi_Q=\phi_f\,,$$
where we used \fref{champ1}, \fref{champ2} and \fref{champ3} for the last identity. This concludes the proof of (\ref{identifiucaitonfirst}), (\ref{keyconvergence}).

\bs
\ni
{\em Step 4. Strong convergence in $L^1$ of $f_n$ to $Q$ }.

\ms
\ni
We first note that \fref{convfnhat} and (\ref{identifiucaitonfirst}) imply 
$$\widehat f_n\to Q \ \ \mbox{in} \ \ L^1.$$
We then claim that the extra gain (\ref{keyconvergence}) -which is again a consequence of the monotonicity (\ref{keymononicity})- allows us to identify $Q$ as the limit of $f_n$. We indeed claim that
(\ref{keyconvergence}) and $|\widehat f_n- Q |_{L^1} \to 0$ imply 
\be
\label{Qfnplus}
\int_{\RR^6}(Q-f_n)_+ dxdv  \to 0\mbox{ as }n\to +\infty.
\ee
Let us assume (\ref{Qfnplus}) and conclude the proof of Theorem \ref{theo2}. 
We first claim that  (\ref{Qfnplus}), together with $\widehat f_n\to Q $ in $L^1$,  imply
\be
\label{fnQplus}
\int_{\RR^6}(f_n-Q)_+ dxdv  \to 0\mbox{ as }n\to +\infty.
\ee
Indeed we  observe that for all $g,h \in L^1(\RR^6)$ with $g\geq 0$, $h\geq 0$, we have
\be
\label{magicidentity}
\int_{0}^{+\infty}  \mbox{meas}\left\{g\leq t<h\right\} dt = \int_{\RR^6}(h-g)_+ dxdv \leq | g- h |_{L^1},
\ee
and thus:
\bee
\int_{\RR^6}(f_n-Q)_+ dxdv &\leq & \int_{\RR^6}(f_n-\widehat f_n)_+ dxdv +\int_{\RR^6}(\widehat f_n-Q)_+ dxdv \\
&\leq &  \int_{0}^{+\infty}\mbox{meas}\left\{\widehat f_n\leq t<f_n\right\} dt+|\widehat f_n- Q |_{L^1}  \\
&&  =\int_{0}^{+\infty}\mbox{meas}\left\{f_n\leq t<\widehat f_n\right\} dt+|\widehat f_n- Q |_{L^1}  \\
&&  =\int_{\RR^6}(\widehat f_n-f_n)_+ dxdv+|\widehat f_n- Q |_{L^1}  \\
&\leq&  \int_{\RR^6}(Q-f_n)_+ dxdv+ \int_{\RR^6}(\widehat f_n-Q)_+ +|\widehat f_n- Q |_{L^1}  \\
&\leq&  \int_{\RR^6}(Q-f_n)_+ dxdv+ 2|\widehat f_n- Q |_{L^1}  \\
\eee
where we repeatedly used \fref{magicidentity} and the fact that $\widehat f_n \in \Eq(f_n)$ implies
$$\forall t>0,\ \  \mbox{meas}\left\{\widehat f_n\leq t<f_n\right\}=\mbox{meas}\left\{f_n\leq t<\widehat f_n\right\}.$$
As $\widehat f_n\to Q $ in $L^1$, we then conclude that (\ref{Qfnplus}) implies (\ref{fnQplus}). 
Finally adding (\ref{Qfnplus}) and  (\ref{fnQplus}) gives
$$| f_n- Q |_{L^1} \to 0\ \mbox{as} \ \ n\to +\infty.$$

Furthermore,  \fref{necessarliy} and the strong convergence $\na \phi_{f_n}\to \na \phi_Q$ in $L^2$ imply: $$||v|^2f_n|_{L^1}\to ||v|^2Q|_{L^1}\mbox{ as }n\to +\infty,$$
Together with the a.e. convergence of $f_n$, this yields the strong convergence of $|v|^2f_n$ to $|v|^2Q$. Note that the uniqueness of the limit now implies the convergence of all the sequence $f_n$ which completes the proof of (\ref{convfn}).

\ss

\noindent {\it Proof of (\ref{Qfnplus})}. It is a consequence of (\ref{keyconvergence}) and of the convergence:
\be
\label{convL1-1}\widehat f_n\to Q \ \ \mbox{in} \ \ L^1.
\ee 
We first claim that  (\ref{keyconvergence}) remains true if one replaces $\phi_{f_{n}}$ by $\phi_{Q}$:
\be
\label{monotonie0}
 \ds  0\leq T_n := \int_{\RR^6}  \left(\frac{|v|^2}{2}+\phi_{Q}(x)\right) \left( f_{n}- f_{n}^{*\phi_{Q}}\right)dxdv   \  \to 0, \quad\mbox{as}\ \  n\to+\infty.
 \ee
 The fact that $T_{n}\geq 0$ can be proved exactly in the same way as for \fref{identity2}, since $f_{n}$ and $f_{n}^{*\phi_{Q}}$ are equimeasurable. Let us now prove that $T_{n} \to 0$.
  We observe from assumption \fref{condfn2} that $|f^{*\phi_{Q}}_{n}|_{L^1}= |f^*_{n}|_{L^1}=|f_{n}|_{L^1} \to |Q|_{L^1}$,  and that
  $$ f_{n}^*(s,\ell) \to Q^*(s,\ell), \ \ \mbox{for}\ \ ae \ \  s>0, \ \ell >0.$$
  This implies that
  $$f^{*\phi_{Q}}_{n} (x,v)\to Q(x,v),\ \  \mbox{for}\ \ ae\ \   (x,v)\in \RR^6.$$
  As a consequence of the Br\'ezis-Lieb Lemma (see \cite{lieb-loss}, Theorem 1.9), we then get
  \be
  \label{convL1-2}
  |f^{*\phi_{Q}}_{n}- Q|_{L^1} \ \ \to 0.
  \ee
  We now write
  
  $$\begin{array}{l}
  \ds  T_{n}- \int_{\RR^6}\left(\frac{|v|^2}{2}+\phi_{f_{n}}(x)\right) \left( f_{n}- f_{n}^{*\phi_{f_{n}}}\right)dxdv=\\
 \ds  \int_{\RR^6} \left(\phi_{Q}-\phi_{f_{n}}\right) \left(  f_{n}- f_{n}^{*\phi_{Q}}\right)+ \int_{\RR^6} \left(\frac{|v|^2}{2}+\phi_{f_{n}}(x)\right) \left( f_{n}^{*\phi_{f_{n}}}-f_{n}^{*\phi_{Q}}\right)dxdv\\
\ds  \leq  \int_{\RR^6} \left|\phi_{f_{n}}- \phi_{Q}\right| \left| f_{n}- f_{n}^{*\phi_{Q}}\right|+  \int_{\RR^6} \left(-\frac{|v|^2}{2}-\phi_{f_{n}}(x)\right)_{+} \left| f_{n}^{*\phi_{f_{n}}}- f_{n}^{*\phi_{Q}}\right|dxdv\\
\ds  \leq |\phi_{f_{n}}- \phi_{Q}|_{L^{\infty}}| f_{n}- f_{n}^{*\phi_{Q}}|_{L^{1}}  +|\phi_{f_{n}}(0)|| f_{n}^{*\phi_{f_{n}}}- f_{n}^{*\phi_{Q}}|_{L^{1}} \ \ \ \to 0,
 \end{array}$$
 where we have used the definition \fref{symgenphi} of $f^{*\phi}$, the uniform convergence of the potential $\phi_{f_{n}}$, the boundedness of $f_{n}$ and $ f_{n}^{*\phi_{Q}}$ in the energy space, and the $L^1$ convergences \fref{convL1-1} and \fref{convL1-2} . Using $T_{n}\geq 0$ and the convergence
 \fref{keyconvergence}, we finally deduce that $T_{n}\to 0$, and \fref{monotonie0} is proved.
 
 Arguing as in the proof of \fref{identity2}, we write \fref{monotonie0} in the following equivalent form
 \be
\label{monotonie-1}
 \ds T_{n}=\int_{t=0}^{+\infty} dt  \left(\int_{S_{1}^n(t)}  \left(\frac{|v|^2}{2}+\phi_{Q}(x)\right)dxdv   - \int_{S_{2}^n(t)}  \left(\frac{|v|^2}{2}+\phi_{Q}(x)\right)dxdv\right) \ \ \  \to 0,
 \ee
 where
$$S^n_{1}(t)=\{(x,v)\in \RR^6,  f_n^{*\phi_{Q}}(x,v)\leq t<f_n(x,v)\}, $$
$$ S^n_{2}(t)=\{(x,v)\in \RR^6,   f_n(x,v)\leq t< f_n^{*\phi_{Q}}((x,v)\}.$$
Now from \fref{implicationtwo}, we have

$$  \frac{|v|^2}{2}+\phi_{Q}(x)\geq (f^*_{n}\circ a_{\phi_{Q}})^{-1}(t, |x\times v|^2),$$
if $(x,v) \in  S^n_{1}(t)$. Thus

\be
\label{inegalite1}T_{n}\geq \int_{t=0}^{+\infty} dt  \left(\int_{S_{1}^n(t)}  (f^*_{n}\circ a_{\phi_{Q}})^{-1}(t, |x\times v|^2)dxdv   - \int_{S_{2}^n(t)}  \left(\frac{|v|^2}{2}+\phi_{Q}(x)\right)dxdv\right).
\ee
As a consequence of the equimeasurability  of $f_n^{*\phi_{Q}}$ and $f_{n}$, we claim that
 \be
 \label{equi-L}
 \left(\int_{S_{2}^n(t)}  - \int_{ S_{1}^n(t)}\right) (f_{n}^*\circ a_{\phi_{Q}})^{-1}(t, |x\times v|^2) dx dv \, = \, 0.
 \ee
 Indeed, we first use the change of variables 
 $$
 r=|x|,\quad u=|v|,\quad \ell=|x\times v|^2,
 $$  to get
\bee
\int_{S_{1}^n(t)}(f_{n}^*\circ a_{\phi_{Q}})^{-1}(t, |x\times v|^2) dx dv & =& \int_{\ell=0}^{\infty}\int_{S^n_{1,\ell}(t)} 
 (f_{n}^*\circ a_{\phi_Q})^{-1}(t, \ell)  d \nu_\ell(r,u)  d\ell,\\  
 &=&
 \int_{\ell=0}^{\infty}(f_{n}^*\circ a_{\phi_Q})^{-1}(t, \ell) \nu_\ell (S^n_{1,\ell})(t) d\ell,
\eee
and  the same identity holds for $S^n_{2}(t)$, 
where $\nu_\ell$ is given  by \fref{nul}, and
$$S^n_{1,\ell}(t)=\{(r,u)\in \Omega_\ell,\,  f_n^{*\phi_{Q}}(r,u,\ell)\leq t<f_n(r,u,\ell)\}, $$
$$ S^n_{2,\ell}(t)=\{(r,u)\in \Omega_\ell,\,   f_n(r,u,\ell)\leq t< f_n^{*\phi_{Q}}(r,u,\ell)\}.$$ 
Since $f_n^{*\phi_{Q}}\in \Eq(f_n)$, we have:
$$ \nu_\ell(S_{1,\ell}^n(t))=\nu_\ell (S_{2,\ell}^n(t))= 4\pi^2 \int_{S_{2,\ell}^n(t)}  \un_{r^2u^2>\ell} (r^2u^2-\ell)^{-1/2} r|u|drdu.$$
This implies \fref{equi-L} and then \fref{inegalite1} gives:
\be
\label{inegalite2}T_{n}\geq \int_{t=0}^{+\infty} dt  \int_{S_{2}^n(t)}  \left[(f^*_{n}\circ a_{\phi_{Q}})^{-1}(t, |x\times v|^2)   -   \left(\frac{|v|^2}{2}+\phi_{Q}(x)\right)\right]dxdv.
\ee
Now from \fref{implicationone}, we have
$$ (f^*_{n}\circ a_{\phi_{Q}})^{-1}(t, |x\times v|^2)\geq \frac{|v|^2}{2}+\phi_{Q}(x),$$
for $(x,v) \in S_{2}^n(t)$. Thus, from \fref{monotonie0} and \fref{inegalite2}, we get

\be
\label{keyaeconv1}
A_{n}=\left[(f^*_{n}\circ a_{\phi_{Q}})^{-1}(t, |x\times v|^2)   -   \left(\frac{|v|^2}{2}+\phi_{Q}(x)\right)\right]\un_{S_{2}^n(t)} (x,v)\to 0, 
\ee
as $n\to + \infty$, for almost every $(t,x,v)\in \RR_{+}\times \RR^3\times \RR^3$. We now claim that this implies
\be
\label{keyaeconv2}
B_{n}=\left[(Q^*\circ a_{\phi_{Q}})^{-1}(t, |x\times v|^2)   -   \left(\frac{|v|^2}{2}+\phi_{Q}(x)\right)\right]\un_{\overline S_{2}^n(t)} (x,v)\to 0, 
\ee
as $n\to + \infty$, for almost every $(t,x,v)\in \RR_{+}\times \RR^3\times \RR^3$, where
 $$\overline S^n_{2}(t)=\{(x,v)\in \RR^6),  f_n(x,v)\leq t< Q(x,v)\}. $$
 To prove \fref{keyaeconv2}, we write 
 $$S_{2}^n= \left(S_{2}^n\backslash \overline S_{2}^n\right) \cup \left(S_{2}^n\cap  \overline S_{2}^n\right), \ \ \ \ \ \  \overline S_{2}^n= \left(\overline S_{2}^n\backslash S_{2}^n\right) \cup \left(S_{2}^n\cap  \overline S_{2}^n\right),$$
 and get
 \be
 \begin{array}{lll}
 \label{decomposition}
 A_{n}-B_{n} &= &\ds \left(\frac{|v|^2}{2}+\phi_{Q}(x) - (Q^*\circ a_{\phi_{Q}})^{-1}(t, |x\times v|^2) \right)\un_{\overline S_{2}^n(t)\backslash S_{2}^n(t)} \\
  & +& \ds \left((f^*_{n}\circ a_{\phi_{Q}})^{-1}(t, |x\times v|^2)-\frac{|v|^2}{2}-\phi_{Q}(x) \right) \un_{S_{2}^n(t)\backslash \overline S_{2}^n(t)} \\
 &+&\ds  \left[ (f^*_{n}\circ a_{\phi_{Q}})^{-1}(t, |x\times v|^2)- (Q^*\circ a_{\phi_{Q}})^{-1}(t, |x\times v|^2) \right]\un_{S_{2}^n(t)\cap  \overline S_{2}^n(t)}
 \end{array}
 \ee
 We shall now examine the behavior of each of  these terms  when $n\to \infty$.
 We first observe from \fref{magicidentity} and \fref{convL1-2} that  
 $$\int_0^{+\infty} \mbox{meas}(S_{2}^n(t)\backslash \overline S_{2}^n(t)) dt \leq |f_{n}^{*\phi_{Q}}-Q|_{L^1} \ \ \to 0,$$
 which implies (up to a subsequence extraction)
 $$\un_{S_{2}^n(t)\backslash \overline S_{2}^n(t)} \longrightarrow 0,  \ \ \ \mbox{for ae} \ (t,x,v)\in \RR_{+}\times \RR^3\times \RR^3.$$
 Using in addition the estimate
 $$\left|(f_{n}^*\circ a_{\phi_{Q}})^{(-1)}(t, |x\times v|^2)\right| \leq |e_{\phi_Q, |x\times v|^2}| \leq |\phi_Q(0)|,$$ 
 we deduce that the first two terms of the decomposition \fref{decomposition} go to $0$ when $n$ goes to infinity, for almost every $(t,x,v)\in \RR_{+}\times \RR^3\times \RR^3.$
 We now treat the third term and show that
 \be
 \label{liminf} q_{0}=\liminf_{n\to \infty} \left[ (f^*_{n}\circ a_{\phi_{Q}})^{-1}(t, |x\times v|^2)- (Q^*\circ a_{\phi_{Q}})^{-1}(t, |x\times v|^2) \right]\un_{S_{2}^n(t)\cap  \overline S_{2}^n(t)} \geq 0,
 \ee
 for almost every $(t,x,v)$. To prove \fref{liminf}, one may  assume that $ \un_{S_{2}^n(t)\cap  \overline S_{2}^n(t)} (x,v)=1$ for $n$ large enough, $(t,x,v)$ being fixed, otherwise $q_{0}=0$ and \fref{liminf} is proved. Let us also recall from standard argument that the strong $L^1$ convergence \fref{condfn2} together with the monotonicity of $f_n^*$ in $e$ and the {\it continuity} of $Q^*$ in $e$ ensure:
$$
 a.e.  \ \ \ell>0 , \ \ \forall e\in(e_{\phi_Q,\ell},0), \ \ f_n^*(a_{\phi_Q}(e,\ell),\ell)\to Q^*(a_{\phi_Q}(e,\ell),\ell).
 $$
Hence, from \fref{changevariables}, we deduce that for a.e. $(x,v)\in \RR^6$, we have
 \be
 \label{vohoehohve}
\forall e\in(e_{\phi_Q,\ell},0), \ \ f_n^*(a_{\phi_Q}(e,\ell),\ell)\to Q^*(a_{\phi_Q}(e,\ell),\ell),\quad \mbox{where }\ell=|x\times v|^2>0.
 \ee
 Let then  $(t,x,v)$ being fixed such that $ \un_{S_{2}^n(t)\cap  \overline S_{2}^n(t)} (x,v)=1$ for $n$ large enough and \fref{vohoehohve} holds. From $$Q(x,v)=Q^*\left(a_{\phi_Q,\ell}\left(\frac{|v|^2}{2}+\phi_Q(x),\ell\right),\ell\right)>t$$ and from the continuity of $Q^*(\cdot,\ell)$, we deduce that
 \be
 \label{:::}
 (Q^*\circ a_{\phi_Q})^{-1}(t,\ell)=\sup\{e\in ]e_{\phi_Q,\ell}, 0[: \ \ Q^*( a_{\phi_Q}(e,\ell),\ell)>t\}.
 \ee
 Take now any $e$ such that 
 \be
 \label{econ}e_{\phi_{Q},\ell}<e<0, \ \ \ \mbox{and}\ \ \  Q^*(a_{\phi_{Q}}(e,\ell),\ell) >t,
 \ee
 then from  \fref{vohoehohve}:
 $$f_{n}^*(a_{\phi_{Q}}(e,\ell),\ell) > t,$$
 for $n$ large enough. Using the definition of the pseudoinverse given in Lemma \ref{lemmapseudoinverse}, we then get 
 $e\leq (f^*_{n}\circ a_{\phi_{Q}})^{-1}(t, |x\times v|^2)$ for $n$ large enough, and hence
 $$e\leq \liminf_{n\to \infty}  (f^*_{n}\circ a_{\phi_{Q}})^{-1}(t, |x\times v|^2).$$
 Since this equality  holds for all $e$ satisfying \fref{econ},  we conclude from \fref{:::} that
 $$\liminf_{n\to \infty}  (f^*_{n}\circ a_{\phi_{Q}})^{-1}(t, |x\times v|^2)\geq (Q^*\circ a_{\phi_{Q}})^{-1}(t, |x\times v|^2),$$
 and \fref{liminf} is proved.
 
 We now turn to  the decomposition \fref{decomposition} and get from \fref{liminf}
 $$\liminf (A_{n}-B_{n})\geq 0,\quad \mbox{for a.e. }(t,x,v).$$
 Finally, observing that $B_{n}\geq 0$ and using \fref{keyaeconv1},  we conclude that
 \fref{keyaeconv2} holds true:
\be
\label{ae1}
\left((Q^*\circ a_{\phi_{Q}})^{(-1)}(t, |x\times v|^2)- \frac{|v|^2}{2}-\phi_Q(x)\right) \un_{\{f_n\leq t<Q\}} \  \ \to  0 ,
\ee
for ae $(t,x,v)$ in $\RR^*_+\times\RR^6$.

Observe now that 
$$t<Q(x,v) \ \ \mbox{implies} \ \ Q(x,v)=F\left(\frac{|v|^2}{2}+\phi_Q(x), |x\times v|^2\right)>t,$$
By Assumption (A) and Corollary \ref{coroll}, $e\to F(e, |x\times v|^2)$ is continuous and strictly decreasing with respect to $e=\frac{|v|^2}{2}+\phi_Q(x)$ for $(x,v)\in \{Q>0\}$, and thus:
$$t<Q(x,v)\ \ \mbox{implies} \ \ (Q^*\circ a_{\phi_{Q}})^{(-1)}(t, |x\times v|^2)- \frac{|v|^2}{2}-\phi_Q(x) >0.$$
We then deduce from \fref{ae1} that
$$\un_{\{f_n\leq t<Q\}} \  \ \to  0, \ \ \mbox{as} \ \ n\to \infty,$$
for ae $(t,x,v)\in\RR^*_+\times\RR^6$. 
Now from $\un_{\{f_n\leq t<Q\}} \leq \un_{\{t<Q\}} $ and
$$\int_0^{\infty} \int_{\RR^6} \un_{\{t<Q\}} dx dv dt = |Q|_{L^1} <+\infty.$$
we may apply the dominated convergence theorem to conclude:
$$\int_0^{\infty} \int_{\RR^6} \un_{\{f_n\leq t<Q\}} dx dv dt \to  0 \ \ \mbox{as} \ \ n\to \infty .$$
Injecting this into \fref{magicidentity} yields \fref{Qfnplus} .\\
This concludes the proof of Theorem \ref{theo2}.


\subsection{Nonlinear stability of $Q$}


We now turn to the proof of Theorem \ref{theo3} which is a direct consequence of Theorem \ref{theo2} and the known regularity of strong solutions to the Vlasov-Poisson system.

\bs
\ni
{\em Proof of Theorem \pref{theo3}.} Let $f_0\in \calE_{rad}\cap \calC^1_c$ and let $f(t)\in \calE_{rad}$ be the corresponding global strong solution to \fref{vp}. Then from the properties of the flow of the Vlasov-Poisson system \fref{vp}, there holds:
\be
\label{cont1}
\na \phi_f\in \calC([0,+\infty),L^2(\RR^3))
\ee
and
\be
\label{cont2}
\forall t\geq 0,\qquad f(t)\in\Eq(f_0),\qquad \calH(f(t))= \calH(f_0).
\ee
Note that $f(t)\in\Eq(f_0)$ means the equality of the symmetric rearrangements at given $\ell$:  for all $t\geq 0$, $f(t)^*=f_0^*$. Recall also from the property of contractivity of the symmetric rearrangement \fref{contra} that:
\be
\label{lieb}
|f_0^*-Q^*|_{L^1}\leq |f_0-Q|_{L^1}\,.
\ee
Remark finally that the following inequality can be proved by interpolation: for all $f\in \calE_{rad}$ satisfying $\calH(f)\leq \calH(Q)+1$, $|f|_{L^1}\leq |Q|_{L^1}+1$ and $|f|_{L^\infty}\leq M$, we have
\be
\label{se}
|\na\phi_f-\na\phi_Q|_{L^2}\leq C_Q\,|f-Q|_{L^1}^{7/12}\,.
\ee
Let us fix $\eps_0>0$ such that
\be
\label{eps0}
C_Q\,\eps_0^{7/12}=\frac{\delta}{2}\,,
\ee
where $C_Q$ is the constant in \fref{se} and $\delta$ is as in Theorem \ref{theo2}.

 An equivalent reformulation of Theorem \ref{theo2} is the following: for all $\eps>0$, there exists $\eta$ satisfying
\be
\label{defeta}
0<\eta\leq \min (\eps_0,1)
\ee
such that the following holds true: if $f\in \calE_{rad}$ is such that
\be
\label{e1}
|f^*-Q^*|_{L^1}<\eta,\quad |f|_{L^\infty}<M,\quad \calH(f)<\calH(Q)+\eta
\ee
and
\be
\label{e2}
|\na \phi_f-\na\phi_Q|_{L^2}<\delta,
\ee
then we have
\be
\label{e3}
|f-Q|_{L^1}<\min(\eps,\eps_0),\quad |f|_{L^\infty}<M,\quad ||v|^2(f-Q)|_{L^1}<\eps.
\ee
Let $f_0\in \calE_{rad}\cap\calC^1_c$ satisfying \fref{eta}. From \fref{cont2} and \fref{lieb}, we first deduce that the corresponding solution $f(t)$ of \fref{vp} satisfies \fref{e1} for all $t\geq 0$. Hence, if we prove that
\be
\label{w}
\forall t\geq 0,\qquad |\na \phi_f(t)-\na\phi_Q|_{L^2}<\delta,
\ee
then it will imply that $f(t)$ satisfies \fref{e3} for all $t\geq 0$, which is nothing but \fref{eta2}.

Therefore, it remains to prove \fref{w}. By \fref{eta} and \fref{defeta}, the initial data $f(0)=f_0$ satisfies
$$|f(0)-Q|_{L^1}<\min(\eps_0,1),\quad |f(0)|_{L^\infty}< M \quad \mbox{and}\quad \calH(f(0))< \calH(Q)+\eta,$$
thus \fref{se} and \fref{eps0} imply that $|\na\phi_f(0)-\na\phi_Q|_{L^2}<\delta/2$. 
Now, assume that \fref{w} is not true. Then there exists $t_1>0$ such that $|\na\phi_{f}(t_1)-\na\phi_Q|_{L^2}\geq \delta$ and, by the continuity property \fref{cont1}, there exists $t_2>0$ such that
\be
\label{delta/3}|\na\phi_f(t_2)-\na\phi_Q|_{L^2}=2\delta/3\,.
\ee
Hence, since the function $f(t_2)$ satisfies \fref{e1} and \fref{e2}, it satisfies \fref{e3}. Therefore, we have
$$|f(t_2)-Q|_{L^1}<\eps_0,\qquad |f(t_2)|_{L^\infty}\leq M \quad \mbox{and}\quad \calH(f(t_2))\leq \calH(Q)+\eta,$$
thus \fref{se} and \fref{eps0} imply that $|\na\phi_f(t_2)-\na\phi_Q|_{L^2}\leq \delta/2$, which contradicts \fref{delta/3}. The proof of Theorem \pref{theo3} is complete.
\qed


\section{Study of the reduced functional $\mathcal J$}
\label{sectionj}


This section is devoted to the proof of Proposition \ref{propA} which requires a detailed study of the reduced functional $\mathcal J$ defined by (\ref{defiphiq}). In particular, we aim at proving that $\mathcal J$ is twice differentiable at $\phi_Q$, that $\phi_Q$ is a critical point and that the Hessian at $\phi_Q$ is definite positive. Remarkably enough, this last fact holds because the Hessian of $\mathcal J$ is deeply connected to the structure of the linearized transport operator close to $\phi_Q$ which is {\it explicit in the radial setting}. The coercivity of the corresponding  Hartree-Fock exchange operator, \cite{lynden},  then follows from Antonov's celebrated coercivity property \cite{A1,A2}.


\subsection{Antonov's coercivity and the structure of the linearized transport operator}
\label{sectantonovcoercive}


Let us start by describing some properties of the linearized transport operator generated by $\phi_Q$ which is deeply connected to the structure of the Hessian of $\mathcal J$, see Proposition \ref{lemA1}. These properties rely on standard abstract functional analysis 
results and a remarkable coercivity property due to Antonov, \cite{A1,A2}.\\

Let $F'_e$ denote the partial derivative with respect to $e$ of the function $F(e,\ell)=Q^*\left(a_{\phi_Q}(e,\ell),\ell\right)$ defined in Assumption (A). Abusing notations, we will note when no confusion is possible:
\be
\label{deffphiprime}
F'_e=F'_e(x,v)=\frac{\partial F} {\partial e}\left(\frac{|v|^2}{2}+\phi_Q(x), |x\times v|^2\right).
\ee
Denote
\be
\label{omega}
\Omega=\left\{(x,v)\in\RR^6\,:\quad Q(x,v)>0\right\}.
\ee
At any $(x,v)\in \Omega$, we have $(\frac{|v|^2}{2}+\phi_Q(x), |x\times v|^2)\in \calO$, where $\calO$ is defined in Assumption (A), hence $F'_e(x,v)<0$. Moreover, the function $(x,v)\mapsto F'_e(x,v)$ is continuous on $\Omega$.

We now consider the $L^2$ weighted Hilbert space:
$$L^{2,r}_{|F'_e|}=\left\{f\in L^{1}_{loc}(\Omega)\mbox{ spherically symmetric with }\int_\Omega \frac{f^2}{|F'_e|}\,dxdv<+\infty\right\}$$ and introduce an orthogonal decomposition:
$$L^{2,r}_{|F'_e|}=L^{2,even}_{|F'_e|}\oplus L^{2,odd}_{|F'_e|}$$ where
$$L^{2,even}_{|F'_e|}=\left\{f\in L^{2,r}_{|F'_e|} \ \ \mbox{with} \ \ f(x,-v)=f(x,v)\right\},$$
$$L^{2,odd}_{|F'_e|}=\left\{f\in L^{2,r}_{|F'_e|} \ \ \mbox{with} \ \ f(x,-v)=-f(x,v)\right\}.$$
We then consider the unbounded transport operator:
$$\calT f=v\cdot\nabla_xf-\nabla \phi_Q\cdot\nabla_vf, \ \ D(\calT)=\left\{f\in L^{2,r}_{|F'_e|}, \ \ \calT f\in L^{2,r}_{|F'_e|}\right\}.$$ Note that $\calC^{\infty}_c(\Omega)\subset D(\calT)$ is dense in $L^{2,r}_{|F'_e|}$ and hence $D(\calT)$ is dense in $L^{2,r}_{|F'_e|}$.
We claim the following properties of $\calT$:

\begin{proposition}[Properties of $\calT$]
\label{propotransport}
(i) {\em Structure of the kernel}: $i\calT$ is a self adjoint operator with kernel:
\be
\label{kernelt}
N(\calT)=\left\{f\in L^{2,r}_{|F'_e|} \ \ \mbox{of the form} \ \ f(x,v)=\tilde{f}\left(\frac{|v|^2}{2}+\phi_Q(x),|x\times v|^2\right)\right\}.
\ee
(ii) {\em Coercivity of the Antonov functional}: The Antonov functional
\be
\label{antonovpositive}
\mathscr A(g,g):=\int \frac{g^2}{|F'_e|}\,dxdv-|\na \phi_g|_{L^2}^2
\ee
is continuous on $L^{2,r}_{|F'_e|}$. Moreover, 
\be
\label{keyestimate}
\forall \xi\in D(\calT)\cap L^{2,odd}_{|F'_e|}, \qquad  \mathscr A(\calT\xi,\calT\xi)\geq \int_\Omega \frac{(\xi)^2}{|F'_e|}\,\frac{\phi_Q'(r)}{r}\,dxdv.
 \ee
 (iii) Let $g\in \left[N(\calT)\right]^{\perp}\cap L^{2,even}_{|F'_e|}$. Then $A(g,g)\geq 0$ and we have $A(g,g)=0$ if and only if $g=0$.
 \end{proposition}

 \begin{proof} {\em Step 1: Description of the kernel}
 
  \ms
\ni
Property {\em (i)} relies on the integration of the characteristic equations associated with $\mathcal Tf=0$ and is a standard consequence of the integrability of Newton's equation with central force field in radial symmetry. The proof follows similarly like for the proof of Jean's theorem in \cite{Batt}, see also \cite{GuoLin}.

\bs
\ni
{\em Step 2. Proof of (ii)} 

 \ms
\ni  
Let $g\in \calC^0_c(\Omega)$. We integrate by parts to get: $$\begin{array}{ll}
\ds |\na \phi_g|_{L^2}^2=-\int g(x,v)\phi_g(x)\,dxdv&\ds \leq |g|_{L^{2,r}_{|F'_e|}}\left|\int_{\RR^6}(\phi_g(x))^2\,F'_e\,dxdv\right|^{1/2}\\[4mm]
&\ds \lesssim |g|_{L^{2,r}_{|F'_e|}}|\na \phi_g|_{L^2}\,,
\end{array}
$$
where we used \fref{houmtsouk} proved in the Appendix. The density of $\calC^0_c(\Omega)$ into $L^{2,r}_{|F'_e|}$ allows us to extend this estimate:
$$\forall  g\in L^{2,r}_{|F'_e|}, \ \ |\nabla\phi_g|_{L^2}\lesssim |g|_{L^{2,r}_{|F'_e|}},$$ and the continuity of (\ref{antonovpositive}) onto $L^{2,r}_{|F'_e|}$ follows.

Antonov's coercivity property is now the following claim: $\forall \xi\in  \calC^{\infty}_c(\Omega)\cap L^{2,odd}_{|F'_e|}$,
\be
\label{coercivityproeprty}
 \mathscr A(\calT\xi,\calT\xi)\geq \int_\Omega \frac{(\xi)^2}{|F'_e|}\,\frac{\phi_Q'(r)}{r}\,dxdv.
 \ee
 In the case where the function $F$ depends only on $e=|v|^2/2 +\phi(x)$,  a proof of this inequality can be found in  \cite{kandrup1,GR4,PA,SDLP}.  In our context $F$ depends on $e=|v|^2/2 +\phi(x)$ and $\ell= |x\times v|^2$, and for a sake of clarity
 and completeness,  we give a proof of this inequality in Appendix \ref{antonov} which is a simple extension of the proof in \cite{GR4}.

Let us extend this estimate to all $\xi\in D(\calT)\cap L^{2,odd}_{|F'_e|}$ using standard regularization arguments. Let $\xi\in D(\calT)\cap L^{2,odd}_{|F'_e|}$ and assume first that $\mbox{Supp}(\xi)\subset \Omega$. From the continuity of $F'_e$, we deduce that $F'_e(x,v)\leq \delta <0$ for all $(x,v)\in \mbox{Supp}(\xi)$. Let a mollifying sequence $\zeta_n(x,v)=\frac{1}{n^6}\zeta(\frac{|x|}{n},\frac{|v|}{n})\in \calC_c^{\infty}(\RR^6)$ with $\zeta\geq 0$, then first from standard regularization arguments: 
$$\zeta_n\star \xi\to \xi, \ \ \zeta_n\star(\calT \xi)\to \calT\xi \ \ \mbox{in} \ \ L^2_{|F'_e|}\ \ \mbox{as} \ \ n\to +\infty,
$$ 
and
$$
\calT(\zeta_n\star\xi)\to \calT\xi \ \ \mbox{in} \ \ L^2_{|F'_e|}\ \ \mbox{as} \ \ n\to +\infty.
$$
Antonov's coercivity property applied to $\zeta_n\star\xi\in \calC_c^{\infty}(\Omega)\cap L^{2,odd}_{|F'_e|}$, the continuity of $\mathscr A$ on $L^2_{|F'_e|}$ and the boundedness of $\frac{\phi'_Q(r)}{r}$ yield the claim.  Consider now a general $\xi\in D(\calT)\cap L^{2,odd}_{|F'_e|}$. We let $\widetilde \chi_n$ a $\calC^\infty$ function such that
\be
\label{chitilde}
\left\{\begin{array}{l}
\widetilde \chi(s)=0 \quad\mbox{for } s\leq \frac 1 {2n},\\[2mm]
\widetilde \chi\mbox{ increasing on }\left[\frac1{2n},\frac 1 n\right],\\[2mm]
\widetilde \chi(s)=1 \quad\mbox{for } s\geq \frac 1 n,
\end{array}\right.
\ee
and we set
\be
\label{chi}
\chi_n(x,v)=\widetilde \chi_n(Q(x,v)).
\ee
Then $\chi_n$ is a $\calC^1$ function with a compact support in $\Omega$, satisfying $\calT \chi_n=0$. Therefore  $\chi_n\xi\in L^{2,odd}_{|F'_e|}$, has compact support in $\Omega$ and $$\calT(\chi_n\xi)=\chi_n\calT\xi\to \calT\xi \ \ \mbox{in} \ \ L^2_{|F'_e|},$$ and hence the previous step and the continuity of $\mathscr A$ on $L^2_{|F'_e|}$ yield \fref{keyestimate}.

 \bs
\ni
{\em Step 3. Proof of (iii)} 

\bs
\ni  
We first observe that the transport operator exchanges parity in $v$: $\forall \xi\in D(\calT)$, 
$$
\left(\xi\in L^{2,odd}_{|F'_e|} \Longrightarrow  \calT \xi\in  L^{2,even}_{|F'_e|}\right), \ \ \left(\xi\in L^{2,even}_{|F'_e|}\Longrightarrow \calT \xi\in  L^{2,odd}_{|F'_e|}\right).
$$
This implies: 
\be
\label{nbv}
 \overline{R(\calT\mbox{}_{|_{\mbox{\footnotesize $L$}^{2,odd}_{|F'_e|}}})}=\overline{R(\calT)}\cap L^{2,even}_{|F'_e|}.
 \ee
On the other hand, $i\calT$ being self-adjoint, there holds --see Cor. II.17, p. 28 in \cite{Brezis}--: 
\be
\label{nbv2}
\overline{R(\calT)}=N(\calT)^{\perp}.
\ee

Let $g\in \left[N(\calT)\right]^{\perp}\cap L^{2,even}_{|F'_e|}$. From \fref{nbv} and \fref{nbv2}, we infer the existence of a sequence $\xi_n\in D(\calT)\cap L^{2,odd}_{|F'_e|}$ such that
\be
\label{txin}
\calT \xi_n\to g \ \ \mbox{in} \ \ L^{2,r}_{|F'_e|}
\ee
as $n\to +\infty$.
Hence, from the continuity of the Antonov functional on $L^{2,r}_{|F'_e|}$, we have
\be
\label{nbv25}
\mathscr A(\calT \xi_n,\calT \xi_n)\to \mathscr A(g,g).
\ee
Moreover, by \fref{keyestimate}, we have
\be
\label{nbv3}
\mathscr A(\calT \xi_n,\calT \xi_n)\geq \int_\Omega \frac{(\xi_n)^2}{|F'_e|}\,\frac{\phi_Q'(r)}{r}\,dxdv\geq 0.
\ee
Thus \fref{nbv25} and \fref{nbv3} imply $\mathscr A(g,g)\geq 0$.

Assume now that $\mathscr A(g,g)=0$. Then \fref{nbv25} and \fref{nbv3} imply that
\be
\label{zidane}\int_\Omega \frac{(\xi_n)^2}{|F'_e|}\,\frac{\phi_Q'(r)}{r}\,dxdv\to 0.
\ee
as $n\to +\infty$. 
Solving the Poisson equation in radial coordinates yields: 
$$
r^2\phi_Q'(r)=4\pi\int_0^r\rho_Q(s)s^2ds.
$$
Denote $r_0=\inf_{(x,v)\in \Omega}|x|$. From the definition \fref{omega} of $\Omega$ and the continuity of $Q$, we have a sequence $r_j\to r_0$, $r_j>r_0$, such that $\rho_Q(r_j)>0$. Hence, for all $r>r_0$, we have $r^2\phi_Q'(r)\geq r_j^2\phi_Q'(r_j)>0$, for $j$ large enough. Thus, the function $\frac{\phi_Q'(r)}{|F'_e|r}$ is continuous and strictly positive on $\Omega$ and \fref{zidane} implies that 
$$\xi_n\to 0 \ \ \mbox{in} \ \ L^{2}_{loc}(\Omega).$$
Therefore, $\calT \xi_n\rightharpoonup 0$ in the distribution sense $\calD'(\Omega)$ and, by \fref{txin}, $g=0$.  This concludes the proof of Proposition \ref{propotransport}.

\end{proof}

A standard consequence of the explicit description of the kernel of $\calT$ given by (\ref{kernelt}) is that we can compute the projection onto $N(\calT)$ --see \cite{GuoLin} for related statements--. For later use, we introduce the following homogeneous Sobolev space
$$\X=\left\{h\in \dot H^1(\RR^3)\quad \mbox{s.t. $h$ is radially symmetric}\right\}.$$
\begin{lemma}[Projection onto the kernel of $\calT$]
\label{lemmaprojection}
Let
\be
\label{D}
{\mathcal D}=\left\{(e,\ell)\in \RR_-^*\times \RR_+^*:\,e>e_{\phi_Q,\ell}\right\},
\ee
where $e_{\phi_Q,\ell}$ is defined by \fref{ephil}.
Given $h\in \X$, we define the projection operator:
\be
\label{pikernel}
{\mathcal P} h(x,v)=\frac{\ds \int_{r_1}^{r_2}\left(e(x,v)-\phi_Q(r)-\frac{\ell(x,v)}{2r^2}\right)^{-1/2}h(r)dr}{\ds \int_{r_1}^{r_2}\left(e(x,v)-\phi_Q(r)-\frac{\ell(x,v)}{2r^2}\right)^{-1/2}dr}\,\un_{(e(x,v),\ell(x,v))\in \mathcal D}\,,
\ee
where $r_1=r_1(\phi_Q,e(x,v),\ell(x,v))$, $r_2=r_2(\phi_Q,e(x,v),\ell(x,v))$ are defined by \fref{r1}, \fref{r2}, and where
$$e(x,v)=\frac{|v|^2}{2}+\phi_Q(x),\qquad \ell(x,v)=|x\times v|^2.$$
Then:
\be
\label{estun}
hF'_e\in L^{2,r}_{|F'_e|},\qquad  ({\mathcal P} h)F'_e\in L^{2,r}_{|F'_e|}
\ee
and
\be
\label{projectionkernel}
({\mathcal P} h)|F'_e|\in N(\calT), \ \ (h-{\mathcal P} h) F'_e\in \left[N(\calT)\right]^{\perp}\cap L^{2,even}_{|F'_e|}
\ee
with $F'_e$ given by \fref{deffphiprime}.
\end{lemma}

The proof is given in Appendix \ref{appendixlemma}.


\subsection{Differentiability of $\mathcal J$}
\label{sectionjdiff}


Our aim in this section is to prove the differentiability of $\mathcal J$ at $\phi_Q$ and to compute the first two derivatives. We shall in particular exhibit an intimate link between the Hessian of $\mathcal J$ and the projection operator (\ref{pikernel}).
\begin{proposition}[Differentiability of $\mathcal J$]
\label{lemA1}
The functional $\mathcal J$ defined by \fref{defiphiq} on $\Phi_{rad}$ satisfies the following properties.\\
(i) {\em Differentiability of $\mathcal J$}. Let $\phi =\phi_f\in \Phi_{rad}$ and $\widetilde \phi=\phi_{\widetilde f}\in \Phi_{rad}$, both nonzero. Then, the functional $$\lambda\mapsto \mathcal J(\phi+\lambda (\widetilde \phi-\phi))$$ is twice differentiable on $[0,1]$.\\
(ii) {\em Taylor expansion of $\mathcal J$ near $\phi_Q$}. Let $R>0$ and
\be
\label{hyptaylor}
f\in B_R:=\left\{g\in\calE_{rad}\mbox{ such that }|g-Q|_{\calE}< R\right\}.
\ee
Then we have the following Taylor expansion near $\phi_Q$:
\be
\label{phiQderiv2}
\mathcal J(\phi_f)-\mathcal J(\phi_Q)=\frac{1}{2}D^2\mathcal J(\phi_Q)(\phi_f-\phi_Q,\phi_f-\phi_Q)+\eps_R(\phi_f)\,|\na \phi_f-\na \phi_Q|_{L^2}^2
\ee
where
$$\eps_R(\phi_f) \to 0 \quad \mbox{as }|\na \phi_f-\na \phi_Q|_{L^2}\to 0\quad \mbox{with }  f\in B_R,$$
and where the second derivative of $\mathcal J$ in the direction $h$ is given by
\be
\label{quadratic}
D^2\mathcal J(\phi_Q)(h,h) =\int_{\RR^3}|\na h|^2\,dx+\int_{\RR^6}h(x)(h(x)-{\mathcal P} h(e,\ell)) F'_e(e,\ell)dxdv
\ee
with ${\mathcal P} h$ given by \fref{pikernel} and $e=\frac{|v|^2}{2}+\phi_Q(x)$, $\ell=|x\times v|^2$.
\end{proposition}
\begin{proof}
Let us decompose $\mathcal J$ into a kinetic part and a potential part:
\be
\label{JJ}
\mathcal J (\phi)= \mathcal J _{Q^*}(\phi) =\calH(Q^{*\phi})+\frac{1}{2}|\na \phi-\na \phi_{Q^{*\phi}}|^2=\frac{1}{2}\int|\nabla \phi|^2dx+ \mathcal J_0(\phi)
\ee with 
\be
\label{A0}
\begin{array}{ll}
\ds \mathcal J_0(\phi)&\ds =\int_{\RR^6}\left(\frac{|v|^2}{2}+\phi(x)\right)Q^*\left(a_\phi\left(\frac{|v|^2}{2}+\phi(x),|x\times v|^2\right),|x\times v|^2\right)dxdv\\[4mm]
&\ds=\int_{\RR^6}\left(\frac{|v|^2}{2}+\phi(x)\right)Q^{*\phi}(x,v)\,dxdv.
\end{array}\ee
Observe that \fref{A0} seems to suggest that two derivatives of $\mathcal J_0$ should involve two derivatives of $Q^*$ and $a_{\phi}$ which are not available in particular from the $\sqrt{\cdot}$ regularity only of the integral \fref{aphi-rul} defining $a_{\phi}$. We claim that is in fact not the case and that suitable integration by parts and change of variables and a careful track of the dependence on $(e,\phi,\ell)$ of the various estimates on $a_{\phi}$ and its derivatives given by Lemmas \ref{lemmadefaphi}, \ref{contaphi}, \ref{cont-aphi-l-phi} will yield the result.\\

\bs
\ni
{\em Step 1. Bounds for the support of $Q^*$}. 

\ms
\ni
In Corollary \ref{coroll}, we have identified the function $Q^*$:
\be
\label{Qstarident}
Q^*(s,\ell)=F\left(a_{\phi_Q}^{-1}(s,\ell),\ell\right),\qquad \forall \ell> 0,\quad \forall s\geq 0.
\ee
Recall that, by Assumption (A), for all $\ell\geq 0$ the function $e\mapsto F(e,\ell)$ is nonincreasing. Let us define
\be
\label{L}
L=\left\{\ell> 0\,:\quad F(e_{\phi_Q,\ell},\ell)>0\right\},
\ee
where $e_{\phi_Q,\ell}$ is defined in Lemma \ref{lemmapotentials}. By Lemma \ref{lemmapotentials} {\em (i)} and by the continuity of $F$, the function $\ell\mapsto F(e_{\phi_Q,\ell},\ell)$ is continuous on $\RR_+^*$, thus $L$ is an open set.

If $\ell\in \RR_+^*\setminus L$, then $a_{\phi_Q}^{-1}(s,\ell)\geq e_{\phi_Q,\ell}$ implies 
$$F(a_{\phi_Q}^{-1}(s,\ell),\ell)\leq F(e_{\phi_Q,\ell},\ell)=0$$ for all $s\geq 0$, thus
\be
\label{ss1}
\forall \ell\in \RR_+^*\setminus L,\qquad Q^*(\cdot,\ell)=0.
\ee
In particular, since $Q=Q^{*\phi_Q}$ is not zero, the measure of $L$ cannot be zero.

Let now $\ell\in L$ and let $$s_0(\ell)=a_{\phi_Q}(e_0(\ell),\ell),$$ where we recall the definition \fref{e0(ell)} of $e_0(\ell)$.
From Assumption (A), Lemma \ref{lemmadefaphi} and \fref{Qstarident}, we infer that the function $Q^*(\cdot,\ell)$ is continuous on $\RR_+$, that its support is $[0,s_0(\ell)]$ and that this function is strictly decreasing and $\calC^1$ on $]0,s_0(\ell)[$. Furthermore, from \fref{minaphi}, we deduce that
\be
\label{ss2}
\forall \ell\in L,\quad 0<s_0(\ell)\leq s_0:=16\pi^2 |Q|_{L^1}|e_0|^{-1/2}.
\ee

Finally, let us prove that the set $L$ is bounded. From Assumption (A),  $(x,v)\mapsto Q(x,v)$ is compactly supported, thus there exist $r_0,u_0>0$ such that $Q(x,v)=0$ for all $(x,v)$ such that $|x|\geq r_0$ or $|v|\geq u_0$.
Hence, we have $Q(x,v)=0$ for all $(x,v)$ such that $|x\times v|^2\geq r_0u_0$ and then, by definition of $Q^*$, $Q^*(\cdot,\ell)=0$ for all $\ell\geq \ell_0:=r_0^2u_0^2$. Therefore, we have
\be
\label{boundL}
L\subset ]0,\ell_0[.
\ee

\bs
\ni
{\em Step 2. First derivative of $\mathcal J_0$}. 

\ms
\ni
We first transform the expression \fref{A0} of $\mathcal J_0$. Using the change of variable \fref{change2bis} and the bounds \fref{ss1}, \fref{ss2} for the support of $Q^*$, we get 
\be
\label{J0new}
\forall \phi \in \Phi_{rad}\backslash \{0\}, \qquad \mathcal J_0(\phi) =\int_{\ell\in L}\int_0^{s_0(\ell)} a_{\phi}^{-1}(s,\ell)Q^*(s,\ell)d\ell ds,
\ee
where we recall that $a_{\phi}^{-1}(\cdot,\ell)$ is defined as the inverse function of $e\mapsto a_{\phi}(\cdot,\ell)$ at given $\phi\in \Phi_{rad}\backslash  \{0\},$ and $\ell>0$.

Let $\phi$ and $\widetilde \phi$ as in Proposition \ref{lemA1} {\em (i)} and $h=\widetilde \phi -\phi$. Let us differentiate the following function with respect to $\lambda \in [0,1]$:
\be
\label{J0newlambda}
\mathcal J_0(\phi+\lambda h) =\int_{L}\int_0^{s_0(\ell)} a_{\phi+\lambda h}^{-1}(s,\ell)Q^*(s,\ell)d\ell ds.
\ee
Let
$$
g(\lambda,s,\ell)=a_{\phi+\lambda h}^{-1}(s,\ell)Q^*(s,\ell).
$$
According to \fref{der2}, we have
$$\frac{\pa g}{\pa \lambda}(\lambda,s,\ell)=\frac{\ds \int_{r_1}^{r_2} \left(a_{\phi+\lambda h}^{-1}(s,\ell)-\psi_{\phi,\ell}(r)-\lambda h(r)\right)^{-1/2}h(r)dr}{\ds \int_{r_1}^{r_2} \left(a_{\phi+\lambda h}^{-1}(s,\ell)- \psi_{\phi,\ell}(r)-\lambda h(r)\right)^{-1/2}dr}Q^*(s,\ell),
$$
where $r_i$, $i=1,2$, shortly denotes $r_i(\phi+\lambda h,a_{\phi+\lambda h}^{-1}(s,\ell),\ell)$ defined by \fref{r1}, \fref{r2}, and $\psi_{\phi,\ell}(r)$ is defined by \fref{psiel}. Therefore,
$$0\leq \frac{\pa g}{\pa \lambda}(\lambda,s,\ell)\leq |h|_{L^\infty}\,Q^*(s,\ell)\in L^1(\RR_+,\RR_+)$$
and we deduce from dominated convergence that $\mathcal J_0$ is differentiable at $\phi$ in the direction $h$ with:
$$
D\mathcal J_0(\phi)(h)=\int_{L}\int_{s=0}^{s_0(\ell)} \frac{\ds \int_{r_1}^{r_2} \left(a_{\phi}^{-1}(s,\ell)- \psi_{\phi,\ell}(r)\right)^{-1/2}h(r)dr}{\ds \int_{r_1}^{r_2} \left(a_{\phi}^{-1}(s,\ell)- \psi_{\phi,\ell}(r)\right)^{-1/2}dr}Q^*(s,\ell)\,d\ell ds.
$$
Using the change of variable $s\mapsto e=a_{\phi}^{-1}(s,\ell)$ and \fref{aphi'}, we now get the following equivalent expression:
\be
\label{derivpremiere}
D\mathcal J_0(\phi)(h)=4\pi^2\sqrt{2}\int_{L}\int_{e_{\phi,\ell}}^0 \int_{r_1}^{r_2} Q^*\left(a_{\phi}(e,\ell),\ell\right)\left(e- \psi_{\phi,\ell}(r)\right)^{-1/2}h(r)drded\ell.
\ee

\bs
\ni
{\em Step 3. Second derivative of $\mathcal J_0$}.

\ms
\ni
Let us now compute the second derivative of $\mathcal J_0(\phi+\lambda h)$ with respect to $\lambda$. First, we write the first derivative in a more convenient  form. Let
\bee
\calD_{\phi,\ell}&=&\left\{(r,e)\in \RR_+^*\times \RR_-^*\quad s.t.\quad e-\phi(r)-\frac{\ell}{2r^2}>0\right\}\\
&=&\left\{(r,e)\in \RR_+^*\times ]e_{\phi,\ell},0[\quad s.t.\quad r_1(\phi,e,\ell)<r<r_2(\phi,e,\ell)\right\}.
\eee
An integration  by parts gives
$$
\begin{array}{l}
\ds \frac{\pa}{\pa \lambda}\mathcal J_0(\phi+\lambda h)=\\
\ds =8\pi^2\sqrt{2}\int_{L}\int_{e=e_{\phi+\lambda h,\ell}}^0 Q^*\left(a_{\phi+\lambda h}(e,\ell),\ell\right)\frac{\pa}{\pa e}\left[\int_{r_1(\phi+\lambda h,e,\ell)}^{r_2(\phi+\lambda h,e,\ell)}\left(e- \psi_{\phi,\ell}(r)-\lambda h(r)\right)^{1/2}h(r)dr\right]ded\ell\\
\ds =-8\pi^2\sqrt{2}\int_{L}\int_{\calD_{\phi+\lambda h,\ell}}\frac{\pa Q^*}{\pa s}\left(a_{\phi+\lambda h}(e,\ell),\ell\right)\frac{\pa a_{\phi+\lambda h}}{\pa e}(e,\ell)\left(e- \psi_{\phi,\ell}(r)-\lambda h(r)\right)^{1/2}h(r)drded\ell
\end{array}
$$
where the boundary terms of the integration by parts vanish. Now, we perform the change of variable $e\mapsto s=a_{\phi+\lambda h}(e,\ell)$ and get
\be
\label{mm1}
 \frac{\pa}{\pa \lambda}\mathcal J_0(\phi+\lambda h)=-8\pi^2\sqrt{2}\int_{L}\int_0^{s_0(\ell)}\int_{r=0}^{+\infty}G(\lambda,s,\ell,r)drdsd\ell,
\ee
with
$$G(\lambda,s,\ell,r)= \frac{\pa Q^*}{\pa s}(s,\ell)\left(a_{\phi+\lambda h}^{-1}(s,\ell)- \psi_{\phi,\ell}(r)-\lambda h(r)\right)_+^{1/2}h(r).
$$
We have
\be
\label{qq1}
\begin{array}{l}
\ds \frac{\pa G}{\pa \lambda}=\frac{1}{2} \frac{\pa Q^*}{\pa s}(s,\ell)\left(\frac{\pa a_{\phi+\lambda h}^{-1}}{\pa \lambda}(s,\ell)-h(r) \right)h(r)\\
\ds \qquad \qquad \qquad \qquad \times\left(a_{\phi+\lambda h}^{-1}(s,\ell)- \psi_{\phi,\ell}(r)-\lambda h(r)\right)^{-1/2}\un_{r_1<r<r_2}\,.
\end{array}
\ee
From \fref{der2}:
\be
\label{qq2}
\left|\frac{\pa a_{\phi+\lambda h}^{-1}}{\pa \lambda}(s,\ell)-h(r)\right|\leq 2|h|_{L^\infty}.
\ee
Moreover, applying \fref{clvm} to the potential $\phi+\lambda h$ gives
\be
\label{qq3}
\left(a_{\phi+\lambda h}^{-1}(s,\ell)- \psi_{\phi,\ell}(r)-\lambda h(r)\right)^{-1/2}\leq \frac{r\sqrt{2r_1r_2}}{\sqrt{(r-r_1)(r_2-r)}}\frac{1}{\sqrt{\ell}}\ee
for $r\in ]r_1,r_2[$. Inserting \fref{qq2} and \fref{qq3} in \fref{qq1} yields
\be
\label{mm}
\left|\frac{\pa G}{\pa \lambda}\right|\leq C|h|_{L^\infty}|rh|_{L^\infty}\left(-\frac{\pa Q^*}{\pa s}(s,\ell)\right)\frac{\sqrt{r_1r_2}}{\sqrt{(r-r_1)(r_2-r)}}\frac{1}{\sqrt{\ell}},
\ee
for $r\in ]r_1,r_2[$, where we recall that $\frac{\pa Q^*}{\pa s}< 0$ for $\ell\in L$ and  $0<s<s_0(\ell)$.

In order to apply Lebesgue's derivation Lemma \ref{lebesgue}, one has to bound the right-hand side of \fref{mm} by suitable $L^1$ functions. To this aim, we need estimates for
$$
r_1=r_1(\phi+\lambda h,a_{\phi+\lambda h}^{-1}(s,\ell),\ell),\quad r_2=r_2(\phi+\lambda h,a_{\phi+\lambda h}^{-1}(s,\ell),\ell).
$$
We claim that, for all $\phi=\phi_f\in \Phi_{rad}$,
\be
\label{ss3}
\forall \ell\leq \ell_0,\quad \forall s\leq s_0(\ell),\qquad r_2(\phi,a_{\phi}^{-1}(s,\ell),\ell)\leq C\,|f|_{L^1}\frac{m_\phi+\ell_0+s_0^2}{m_{\phi}^2},
\ee
where $m_\phi$ is defined by \fref{mphi} and where $C$ is a universal constant. Let us assume \fref{ss3} and conclude the computation of the second derivative. Let
$$m=\min(m_\phi,m_{\widetilde \phi})>0,\quad M=\max(|f|_{L^1},|\widetilde f|_{L^1}).$$
For $\lambda \in [0,1]$, we have
$$m_{\phi+\lambda h}=\inf_{r>0}(r+1)((1-\lambda)|\phi(r)|+\lambda |\widetilde\phi(r)|)\geq (1-\lambda) m_\phi+\lambda m_{\widetilde \phi}\geq m>0.$$
From \fref{boundL}, \fref{mm} and \fref{ss3}, we get
\be
\label{mm2}
\left|\frac{\pa G}{\pa \lambda}\right|\leq C|h|_{L^\infty}|rh|_{L^\infty}\,q_\lambda(s,\ell,r),
\ee
for $\ell\in L$, $s\leq s_0(\ell)$, with
$$0\leq q_\lambda(s,\ell,r)=-\frac{\pa Q^*}{\pa s}(s,\ell)\frac{\un_{r_1< r< r_2}}{\sqrt{(r-r_1)(r_2-r)}}\frac{1}{\sqrt{\ell}},$$
and where the constant $C$ only depends on $m$, $M$, $\ell_0$ and $s_0$.
By Lemma \ref{contaphi}  and the continuity of $r_i(\phi+\lambda h,e,\ell)$, $i=1,2$ with respect to $\lambda$, we have for all $\lambda_0\in [0,1]$
$$q_\lambda(s,\ell,r)\to q_{\lambda_0}(s,\ell,r)\quad \mbox{as }\lambda\to \lambda_0\quad \mbox{for a.e. }r,s,\ell.$$
Moreover,
$$\int_{L}\int_0^{s_0(\ell)}\int_{0}^{+\infty}q_\lambda(s,\ell,r)drdsd\ell =\pi \int_LQ^*(0,\ell)\frac{d\ell}{\sqrt{\ell}}.$$
Since $Q^*\in L^\infty(\RR_+\times\RR_+)$ and from \fref{boundL}, this integral is finite and its value is independent of $\lambda$. Now we invoke the Br\'ezis-Lieb Lemma to conclude:
$$q_\lambda(s,\ell,r)\to q_{\lambda_0}(s,\ell,r)\quad \mbox{as }\lambda\to \lambda_0\quad \mbox{in } L^1([0,s_0]\times [0,\ell_0]\times \RR_+).$$
We conclude from \fref{mm1}, \fref{mm2} and Lemma \ref{lebesgue} that $\lambda\mapsto \mathcal J_0(\phi+\lambda h)$ is twice differentiable on $[0,1]$. In particular, $\mathcal J_0$ is twice differentiable at $\phi$ in the direction $h$ with:
\bee
&&D^2\mathcal J_0(\phi)(h,h)=-4\pi^2\sqrt{2}\int_{L}\int_0^{s_0(\ell)}\int_{r_1}^{r_2}
 \frac{\pa Q^*}{\pa s}(s,\ell)\left(\frac{\pa a_{\phi}^{-1}}{\pa \lambda}(s,\ell)-h(r) \right)h(r)\times \\
&&\qquad \qquad \qquad \qquad \times\left(a_{\phi}^{-1}(s,\ell)- \psi_{\phi,\ell}(r)\right)^{-1/2}drdsd\ell.
\eee
By \fref{der2}, this expression can be simplified into
\be
\label{derivseconde}
\begin{array}{l}
\ds D^2\mathcal J_0(\phi)(h,h)=\\[4mm]
\ds=4\pi^2\sqrt{2}\int_{L}\int_0^{s_0(\ell)}\frac{\pa Q^*}{\pa s}(s,\ell)\int_{r_1}^{r_2}\left(a_{\phi}^{-1}(s,\ell)- \psi_{\phi,\ell}(r)\right)^{-1/2}(h(r))^2drdsd\ell\\[4mm]
\ds-4\pi^2\sqrt{2}\int_{L}\int_0^{s_0(\ell)} \frac{\pa Q^*}{\pa s}(s,\ell)\frac{\ds \left(\int_{r_1}^{r_2} \left(a_{\phi}^{-1}(s,\ell)- \psi_{\phi,\ell}(r)\right)^{-1/2}h(r)dr\right)^2}{\ds \int_{r_1}^{r_2} \left(a_{\phi}^{-1}(s,\ell)- \psi_{\phi,\ell}(r)\right)^{-1/2}dr}dsd\ell.
\end{array}
\ee
{\it Proof of the claim \fref{ss3}}. In order to show that $r_1$ and $r_2$ are not allowed to go to infinity, we shall use \fref{majr2-simple}. To this aim, we first need to show that $e=a_{\phi}^{-1}(s,\ell)$ is not allowed to go to zero when $\ell\leq \ell_0$ and $s\leq s_0(\ell)$. From \fref{estiminf}, we deduce that, for $\ell<\ell_0$, $s\leq s_0(\ell)$,
\be
\label{ss3inter}
|e|^{1/2}\geq \min\left(\frac{m_{\phi}}{2\sqrt{2m_{\phi}+\ell}},\frac{4\pi^2}{3}\frac{m_{\phi}}{s}\right)\geq C\min\left(\frac{m_{\phi}}{\sqrt{m_{\phi}+\ell_0}},\frac{m_{\phi}}{s_0}\right).
\ee
Finally \fref{ss3} can be deduced from \fref{majr2-simple} and \fref{ss3inter}.

\bs
\ni
{\em Step 4. Identification of the first and second derivatives of $\mathcal J$ at $\phi_Q$}.

\ms
\ni
Let $f\in \calE_{rad}$ and $h=\phi_f-\phi_Q$. We claim that
\be
\label{first}
D\mathcal J(\phi_Q)(h)=0.
\ee
In order to prove this claim, we first remark from \fref{JJ} that
\be
\label{first2}
D\mathcal J(\phi_Q)(h)=D\mathcal J_0(\phi_Q)(h)+\int_{\RR^3}\na\phi_Q\cdot \na h\,dx.
\ee
Moreover, by \fref{derivpremiere}, we have
$$
\begin{array}{l}
\ds D\mathcal J_0(\phi_Q)(h)= \\
\ds\, =4\pi^2\sqrt{2}\int_{L}\int_{e_{\phi_Q,\ell}}^0 \int_{0}^{+\infty} F(e,\ell)\left(e- \psi_{\phi_Q,\ell}(r)\right)^{-1/2}h(r)\un_{e- \psi_{\phi_Q,\ell}(r)>0}drded\ell.
\end{array}
$$
where we used \fref{qstarphiq}. Applying the change of variable $e\mapsto u=\sqrt{2(e-\phi_Q(r))}$, it comes
\bee
D\mathcal J_0(\phi_Q)(h)&=&
2\int_{0}^{+\infty}\int_{0}^{+\infty}\int_{0}^{+\infty}F\left(\frac{u^2}{2}+\phi_Q(r),\ell\right)\,h(r)\,\un_{ru>\ell}\,d\nu_\ell d\ell\\
&=&\int_{\RR^6}Q(x,v)h(x)\,dxdv,
\eee
where we used \fref{changevariables}, Assumption (A), and recall that $h$ is radially symmetric. Hence, from the Poisson equation, we deduce after an integration by parts that
$$D\mathcal J_0(\phi_Q)(h)=-\int_{\RR^3}\na\phi_Q\cdot \na h\,dx,
$$
which together with \fref{first2} implies \fref{first}.

Let us now identify the right second derivative of $\mathcal J$ at $\phi_Q$. We have
\be
\label{seconde}
D^2\mathcal J(\phi_Q)(h,h)=D^2\mathcal J_0(\phi_Q)(h,h)+\int_{\RR^3}|\na h|^2\,dx
\ee
and, by \fref{derivseconde},
$$
\begin{array}{l}
\ds D^2\mathcal J_0(\phi_Q)(h,h)=\\[4mm]
\ds=4\pi^2\sqrt{2}\int_{L}\int_0^{s_0(\ell)}\frac{\pa Q^*}{\pa s}(s,\ell)\int_{r_1}^{r_2}\left(a_{\phi_Q}^{-1}(s,\ell)- \psi_{\phi_Q,\ell}(r)\right)^{-1/2}(h(r))^2drdsd\ell\\[4mm]
\ds-4\pi^2\sqrt{2}\int_{L}\int_0^{s_0(\ell)} \frac{\pa Q^*}{\pa s}(s,\ell)\frac{\ds \left(\int_{r_1}^{r_2} \left(a_{\phi_Q}^{-1}(s,\ell)- \psi_{\phi_Q,\ell}(r)\right)^{-1/2}h(r)dr\right)^2}{\ds \int_{r_1}^{r_2} \left(a_{\phi_Q}^{-1}(s,\ell)- \psi_{\phi_Q,\ell}(r)\right)^{-1/2}dr}dsd\ell.
\end{array}
$$
Using first the change of variable $s\mapsto e=a_{\phi_Q}^{-1}(s,\ell)$, \fref{aphi'} and \fref{qstarphiq}, we get
$$
\begin{array}{l}
\ds D^2\mathcal J_0(\phi_Q)(h,h)=\\[4mm]
\ds=4\pi^2\sqrt{2}\int_{L}\int_{e_{\phi_Q,\ell}}^{e_0(\ell)} F'_e(e,\ell) \int_{r_1}^{r_2}\left(e- \psi_{\phi_Q,\ell}(r)\right)^{-1/2}(h(r))^2drded\ell\\[4mm]
\ds-4\pi^2\sqrt{2}\int_{L}\int_{e_{\phi_Q,\ell}}^{e_0(\ell)}F'_e(e,\ell) \frac{\ds \left(\int_{r_1}^{r_2} \left(e-\psi_{\phi_Q,\ell}(r)\right)^{-1/2}h(r)dr\right)^2}{\ds \int_{r_1}^{r_2} \left(e- \psi_{\phi_Q,\ell}(r)\right)^{-1/2}dr}ded\ell.
\end{array}
$$
We next apply the change of variable $e\mapsto u=\sqrt{2(e-\phi_Q(x))}$ and use \fref{changevariables} to get:
$$
D^2\mathcal J_0(\phi_Q)(h,h)=\int_{\RR^6}F'_e(e,\ell)(h(x))^2\,dxdv-\int_{\RR^6}F'_e(e,\ell)h(x){\mathcal P} h(e,\ell) dxdv,
$$
where we used the definition \fref{pikernel} and where we shortly denoted
$$e=\frac{|v|^2}{2}+\phi_Q(x),\qquad \ell=|x\times v|^2.$$
This together with \fref{seconde} concludes the proof of \fref{quadratic}.

\bs
\ni
{\em Step 5. Proof of the Taylor expansion \fref{phiQderiv2}}.

\ms
\ni
We are now ready to prove the Taylor expansion \fref{phiQderiv2}. We first deduce from \fref{first} and from the fact that $\mathcal J(\phi_Q+\lambda h)$ twice differentiable with respect to $\lambda$ that
$$
\mathcal J(\phi_Q+h)-\mathcal J(\phi_Q)=\int_0^1 (1-\lambda)\frac{\pa^2 }{\pa \lambda^2}\mathcal J(\phi_Q+\lambda h)\,d\lambda.
$$
Hence, for $h\neq 0$,
\be
\label{qsd}
\begin{array}{l}
\ds \mathcal J(\phi_Q+h)-\mathcal J(\phi_Q)-\frac{1}{2}D^2\mathcal J(\phi_Q)(h,h)=\\[4mm]
\ds \qquad =\int_0^1 (1-\lambda)\left(D^2\mathcal J(\phi_Q+\lambda h)-D^2\mathcal J(\phi_Q)\right)(h,h)\,d\lambda
\\[4mm]
\ds \qquad =|\na h|_{L^2}^2\int_0^1 (1-\lambda)\left(D^2\mathcal J_0(\phi_Q+\lambda h)-D^2\mathcal J_0(\phi_Q)\right)\left(\frac{h}{|\na h|_{L^2}},\frac{h}{|\na h|_{L^2}}\right)\,d\lambda.
\end{array}
\ee
We now claim the following continuity property:
\be
\label{conti2}
\sup_{\lambda\in [0,1]}\sup_{|\na \widetilde h|_{L^2}=1}\left|\left(D^2\mathcal J_0(\phi_Q+\lambda (\phi_f-\phi_Q)-D^2\mathcal J_0(\phi_Q)\right)(\widetilde h,\widetilde h)\right|\to 0
\ee
as $|\na \phi_f-\na \phi_Q|_{L^2}\to 0$, $f$ satisfying \fref{hyptaylor}.

Assume \fref{conti2}. Then:
$$\int_0^1 (1-\lambda)\left(D^2\mathcal J_0(\phi_Q+\lambda h)-D^2\mathcal J_0(\phi_Q)\right)\left(\frac{h}{|\na h|_{L^2}},\frac{h}{|\na h|_{L^2}}\right)\,d\lambda\to 0$$
and  \fref{qsd} now yields \fref{phiQderiv2}.

\bs
\ni
{\em Proof of \fref{conti2}}. We argue by contradiction and consider $\eps>0$, $f_n$ satisfying \fref{hyptaylor}, $\widetilde h_n$ and $\lambda_n\in [0,1]$  such that
\be
\label{contr1}
|\na \phi_{f_n}-\na \phi_Q|_{L^2}<\frac{1}{n}\,,\qquad |\na \widetilde h_n|_{L^2}=1\,,
\ee
and
\be
\label{contr2}
\left|D^2\mathcal J_0(\phi_Q+ \lambda_n (\phi_{f_n}-\phi_Q))(\widetilde h_n,\widetilde h_n)-D^2\mathcal J_0(\phi_Q)(\widetilde h_n,\widetilde h_n)\right|>\eps.\ee
We denote $h_n=\lambda_n (\phi_{f_n}-\phi_Q)$. Recall from \fref{derivseconde}:
\be
\label{derivseconde2}
\begin{array}{l}
\ds D^2\mathcal J_0(\phi_Q+h_n)(\widetilde h_n,\widetilde h_n)=\\[4mm]
\ds=4\pi^2\sqrt{2}\int_{L}\int_0^{s_0(\ell)}\frac{\pa Q^*}{\pa s}(s,\ell)\int_{r_1^n}^{r_2^n}\left(a_{\phi_Q+h_n}^{-1}(s,\ell)- \psi_{\phi_Q,\ell}(r)-h_n(r)\right)^{-1/2}(\widetilde h_n(r))^2drdsd\ell\\[4mm]
\ds-4\pi^2\sqrt{2}\int_L\int_0^{s_0(\ell)} \frac{\pa Q^*}{\pa s}(s,\ell)\frac{\ds \left(\int_{r_1^n}^{r_2^n} \left(a_{\phi_Q+h_n}^{-1}(s,\ell)- \psi_{\phi_Q,\ell}(r)-h_n(r)\right)^{-1/2}\widetilde h_n(r)dr\right)^2}{\ds \int_{r_1^n}^{r_2^n} \left(a_{\phi_Q+h_n}^{-1}(s,\ell)- \psi_{\phi_Q,\ell}(r)-h_n(r)\right)^{-1/2}dr}dsd\ell,
\end{array}
\ee
 where we have denoted, for $i=1,2$,
$$r_i^n=r_i\left(\phi_Q+h_n,a_{\phi_Q+h_n}^{-1}(s,\ell),\ell\right).$$
By \fref{contr1} and standard radial Sobolev embeddings, the sequence of radially symmetric functions $\widetilde h_n$ is compact in $L^\infty([a,b])$ for all $0<a<b$. By diagonal extraction, we deduce the pointwise convergence of $\widetilde h_n$ (up to a subsequence) to a function $h$:
\be
\label{pointwise}
\forall r\in \RR_+^*\qquad \widetilde h_n(r)\to\widetilde h(r)\quad \mbox{as }n\to +\infty.
\ee
Moreover, 
\be
\label{cs}
r^{1/2}|\widetilde h_n(r)|\leq \left(\int_r^{+\infty}s^2(\widetilde h_n'(s))^2\,ds\right)^{1/2}\leq |\na \widetilde h_n|_{L^2}=1,
\ee
thus, in particular, $r^{1/2}\widetilde h$ belongs to $L^\infty(\RR_+)$.

Let us analyze the convergence of \fref{derivseconde2}. In a first step, recalling \fref{ss1} and \fref{ss2}, we fix $\ell\in L$ and $s\in]0,s_0(\ell)]$ and set
$$e_n=a_{\phi_Q+h_n}^{-1}(s,\ell),\qquad e_\infty=a_{\phi_Q}^{-1}(s,\ell)<0.$$
From \fref{contr1}, the uniform bound of $f_n$ in $\calE_{rad}$ and Lemma \ref{contaphi}, we have:
\be
\label{clpm}
e_n\to e_\infty\quad \mbox{as }n\to +\infty.
\ee
For $k=0,1$ or $2$, we introduce the functions
\be
\label{not55}
g_k(n,s,\ell,r)=\left(e_n- \psi_{\phi_Q,\ell}(r)-h_n(r)\right)^{-1/2}(\widetilde h_n(r))^k\,\un_{\psi_\ell(r)+h_n(r)<e_n}
\ee
and
$$g_k(\infty,s,\ell,r)=\left(e_\infty- \psi_{\phi_Q,\ell}(r)\right)^{-1/2}(\widetilde h(r))^k\,\un_{\psi_\ell(r)<e_\infty}.$$
We claim: $\forall \ell\in L$, $\forall s\in]0,s_0(\ell)]$,
\be
\label{clclaim}
\int_{0}^{+\infty}g_k(n,s,\ell,r)\,dr \to \int_{0}^{+\infty}g_k(\infty,s,\ell,r)\,dr \quad \mbox{for }k=0,1\mbox{ or }2.
\ee
Indeed, from \fref{pointwise}, \fref{clpm}, \fref{contr1} and the bound of $f_n$ which imply $|h_n|_{L^\infty}\to 0$, we first deduce that, for all $s>0$, $\ell>0$, the function $g_k(n,s,\ell,r)$ converges pointwise in $r\in \RR_+^*$ to the function $g_k(\infty,s,\ell,r)$, for $k=0,1$ or $2$, as $n\to +\infty$. Moreover, by applying \fref{clvm} to the function $\phi_Q+h_n$, we get
\bea
0\leq g_k(n,s,\ell,r)&\leq& \frac{r|\widetilde h_n(r)|^k\sqrt{r_1^nr_2^n}\,\un_{r_1^n<r<r_2^n}}{\sqrt{\ell(r-r_1^n)(r_2^n-r)}}\nonumber\\
&\leq&|r^{1/2}\widetilde h_n(r)|_{L^\infty}^k\frac{(r_2^n)^{2-k/2}\,\un_{r_1^n<r<r_2^n}}{\sqrt{\ell(r-r_1^n)(r_2^n-r)}}\nonumber\\
&\leq&\frac{(r_2^n)^{2-k/2}\,\un_{r_1^n<r<r_2^n}}{\sqrt{\ell(r-r_1^n)(r_2^n-r)}},\label{majgk2}
\eea
where we used \fref{cs}. Now we observe from the uniform bound of $f_n$ in $\calE$, from \fref{contr1} and from \fref{defm} that we have
$$\sup_n|Q+\lambda_n (f_n-Q)|_{L^1}<\infty \quad \mbox{and}\quad 0<\inf_n\inf_{r\geq  0} (r+1)|\phi_Q (r)+h_n(r)|<+\infty.$$
Injecting the bound \fref{ss3} on $r_2^n$ into \fref{majgk2} gives:
\be
0\leq g_k(n,s,\ell,r)\leq K\,\frac{\un_{r_1^n<r<r_2^n}}{\sqrt{\ell(r-r_1^n)(r_2^n-r)}}.\label{majgk}
\ee
Denote $r_i^\infty:=r_i(\phi_Q,e_\infty,\ell)$ for $i=1,2$ and
$$q_n(r)=\frac{\un_{r_1^n<r<r_2^n}}{\sqrt{\ell(r-r_1^n)(r_2^n-r)}}\quad \mbox{and}\quad q_\infty(r):=\frac{\un_{r_1^\infty<r<r_2^\infty}}{\sqrt{\ell(r-r_1^\infty)(r_2^\infty-r)}}.$$
By Lemma \ref{contaphi}, we have $r_i^n\to r_i^\infty$ as $n\to +\infty$, for $i=1,2$.
Therefore, the function $q_n$ converges pointwise to $q_\infty$. Since
\be
\label{intpi}
\int_0^{+\infty}q_n(r)\,dr=\int_0^{+\infty}q_\infty(r)\,dr=\frac{\pi}{\sqrt{\ell}},
\ee
we deduce from the Br\'ezis-Lieb lemma that $q_n\to q_\infty$ in $L^1(\RR_+)$ as $n\to +\infty$. Finally, by applying the generalized dominated convergence theorem, we obtain \fref{clclaim}.

Now, we remark that, with the notation \fref{not55}, \fref{derivseconde2} reads
$$
D^2\mathcal J_0(\phi_Q+h_n)(\widetilde h_n,\widetilde h_n)=-4\pi^2\sqrt{2}\int_L\int_0^{s_0(\ell)}G(n,s,\ell)dsd\ell
$$
and
$$
D^2\mathcal J_0(\phi_Q)(\widetilde h,\widetilde h)=-4\pi^2\sqrt{2}\int_L\int_0^{s_0(\ell)}G(\infty,s,\ell)dsd\ell,
$$
with
$$
G(n,s,\ell)=- \frac{\pa Q^*}{\pa s}(s,\ell)\left[\int_{0}^{\infty}g_2(n,s,\ell,r)dr-\frac{\ds \left(\int_{0}^{\infty}g_1(n,s,\ell,r)dr\right)^2}{\ds \int_{0}^{\infty}g_0(n,s,\ell,r)dr}\right]
$$
and
$$
G(\infty,s,\ell)=- \frac{\pa Q^*}{\pa s}(s,\ell)\left[\int_{0}^{\infty}g_2(\infty,s,\ell,r)dr-\frac{\ds \left(\int_{0}^{\infty}g_1(\infty,s,\ell,r)dr\right)^2}{\ds \int_{0}^{\infty}g_0(\infty,s,\ell,r)dr}\right].
$$
From \fref{clclaim}, we get that, for all $\ell\in L$ and $s\in]0,s_0(\ell)[$,
\be
\label{clclaim2}
G(n,s,\ell)\to G(\infty,s,\ell) \quad \mbox{as }n\to +\infty.
\ee
Moreover, by the Cauchy-Schwarz inequality, we have
\be
\label{cs2}
0\leq \left(\int_{0}^{\infty}g_1(n,s,\ell,r)dr\right)^2\leq \left(\int_{0}^{\infty}g_0(n,s,\ell,r)dr\right)\left(\int_{0}^{\infty}g_2(n,s,\ell,r)dr\right).
\ee
Therefore, \fref{majgk}, \fref{intpi} and \fref{cs2} yield the estimate
\be
\label{clclaim3}
0\leq G(n,s,\ell)\leq -\frac{\pi K}{\sqrt{\ell}}\frac{\pa Q^*}{\pa s}(s,\ell)\,.
\ee
Remark that the function in the right-hand side of \fref{clclaim3} belongs to $L^1$ since it is nonnegative and
$$\int_L\int_0^{s_0(\ell)}-\frac{\pa Q^*}{\pa s}(s,\ell)\,ds\frac{d\ell}{\sqrt{\ell}}=\int_LQ^*(0,\ell)\frac{d\ell}{\sqrt{\ell}}<+\infty,$$
where we used \fref{boundL} and $|Q^*(0,\ell)|\leq |Q|_{L^\infty}$.
Finally, we deduce from \fref{clclaim2}, \fref{clclaim3} and from dominated convergence that
\be
\label{cbovbovb}
D^2\mathcal J_0(\phi_Q+h_n)(\widetilde h_n,\widetilde h_n)\to D^2\mathcal J_0(\phi_Q)(\widetilde h,\widetilde h)
\ee
as $n\to +\infty$. 
This contradicts \fref{contr2} and concludes the proof of \fref{conti2}.\\
This concludes the proof of Proposition \ref{lemA1}.
\end{proof}

\begin{remark}
\label{compact}  Note that a similar argument like for the proof of \fref{cbovbovb} gives that for all bounded sequence $\widetilde h_n\in\X$, after extraction of a subsequence, 
\be
\label{comp}
D^2\mathcal J_0(\phi_Q)(\widetilde h_n,\widetilde h_n)\to D^2\mathcal J_0(\phi_Q)(\widetilde h,\widetilde h)
\ee
as $n\to +\infty$. Indeed, we never used for the proof of \fref{cbovbovb} the fact that $\tilde{h}_n$ is a Poisson field. The important consequence is then that the quadratic form $D^2\mathcal J_0(\phi_Q)$ is {\it compact} on $\X$. 
\end{remark}


\subsection{Proof of Proposition \ref{propA}}
\label{sectcoercive}


 We are now in position to conclude the proof of Proposition \ref{propA}. The coercivity property \fref{coerc} will appear as a consequence of the fact that the Hessian \fref{phiQderiv2} can be connected to the Antonov functional \fref{antonovpositive} via the projection operator \fref{pikernel}, and the key step is then Antonov's coercivity property \fref{keyestimate}.

\begin{proof}
 
 \bs
\ni
{\em Step 1. Strict positivity of the Hessian}

\bs
\ni
 We claim that: 
 \be
\label{positivity}
\forall h\in \X\setminus \{0\},\qquad D^2\mathcal J(\phi_Q)(h,h)>0.
\ee
Indeed, let $h\in \dot{H}^1_{r}\setminus \{0\}$ and consider the projection ${\mathcal P} h$ given by  \fref{pikernel}. From \fref{estun}, the functions $hF'_e$ and $({\mathcal P} h)F'_e$ belong to $L^{2,r}_{|F'_e|}$ and hence $g=(h-{\mathcal P}h)F'_e\in L^{2,r}_{|F'_e|}$. By the orthogonality property \fref{projectionkernel}, we have
$$\begin{array}{ll}
\ds \int \left(h-{\mathcal P} h\right)^2 F'_e\,dxdv&\ds =-\int\left(h-{\mathcal P}h\right)^2 F'_e\,dxdv+2\int h\left(h-{\mathcal P} h\right) F'_e\,dxdv\\[3mm]
&\ds =\int \frac{g^2}{|F'_e|}\,dxdv-2\int \na h\cdot \na \phi_g\,dx,
\end{array}$$
where we used the Poisson equation. We may thus rewrite the Hessian \fref{quadratic}:
\bee
D^2\mathcal J(\phi_Q)(h,h) & = & \int|\nabla h|^2dx+\int_{\RR^6}\left(h-{\mathcal P} h\right)^2 F'_e\,dxdv\\
& = & \int \frac{g^2}{|F'_e|}\,dxdv-|\na \phi_g|_{L^2}^2+|\na h - \na \phi_g|^2_{L^2}\\
&= &\mathscr A(g,g)+|\na h - \na \phi_g|^2_{L^2}.
\eee
Now, from \fref{projectionkernel} and Proposition \ref{propotransport} {\em (iii)}, we deduce that $\mathscr A(g,g)\geq 0$. Therefore, $D^2\mathcal J(\phi_Q)(h,h)$ is nonnegative. Moreover, if $D^2\mathcal J(\phi_Q)(h,h)=0$, then $\mathscr A(g,g)= |\na h - \na \phi_g|_{L^2}=0$ and using again Proposition \ref{propotransport} {\em (iii)} enables to conclude that $g=\phi_g=h=0$. This ends the proof of \fref{positivity}.

\bs
\ni
{\em Step 2. Coercivity of the Hessian and conclusion}

\ms
\ni
In Remark \ref{compact}, we have seen that the quadratic form $D^2\mathcal J_0(\phi_Q)$ is compact on $\X$. Hence from \fref{JJ}, the Fredholm alternative can be applied to the quadratic form $D^2\mathcal J(\phi_Q)$. Together with the strict positivity property \fref{positivity}, this implies the coercivity of this quadratic form:
\be
\label{coercivity}
\forall h\in \X\qquad D^2\mathcal J(\phi_Q)(h,h)\geq c|\na h|_{L^2}^2\,,
\ee
for some universal constant $c>0$.\\
We now may conclude the proof of (\ref{coerc}). Let $R>0$ be fixed. From Proposition \ref{lemA1} {\em (ii)}, there exists $\delta_0(R)$ -- chosen in $]0,\frac{1}{2}|\na\phi_Q|_{L^2}]$ -- such that, for all $f\in \calE_{rad}$ satisfying
$$|f-Q|_{\calE}\leq R,\qquad |\na\phi_f-\na\phi_Q|_{L^2}\leq \delta_0(R),$$
we have
$$\eps_R(\phi_f)\leq \frac{c}{4},$$
where $c$ is the constant in \fref{coercivity} and $\eps_R$ is defined in \fref{phiQderiv2}. Hence, for such $f$, we deduce from \fref{coercivity} and \fref{phiQderiv2} that
$$\mathcal J(\phi_Q+h)-\mathcal J(\phi_Q)\geq \frac{c}{4}|\na h|_{L^2}^2.$$
The proof of Proposition \ref{propA} is complete.

 \end{proof}

\vspace{1cm}
\titlecontents{section} 
[1.5em] 
{\vspace*{0.5em}} 
{\contentslabel[\bf]{0em}} 
{\hspace*{-2.5em}} 
{\titlerule*[0.5pc]{\rm.}\contentspage}

\begin{appendix}


\section{A dominated convergence lemma}

\bs
\begin{lemma}
\label{lebesgue}
Let $I$ be an interval of $\RR$ and let $g(\lambda,r)$ be a real-valued function in $\calC^0(I,L^1(\RR_+))$. Let 
$$G(\lambda)=\int_{\RR_+}g(\lambda,r)dr$$
and denote by $\ds \frac{\pa g}{\pa \lambda}$ the weak partial derivative of $g$ with respect to $\lambda$. Assume that $g$ satisfies the following assumptions:
\begin{itemize}
\item[(i)] $\ds \frac{\pa g}{\pa \lambda}\in L^1(I\times \RR_+)$ and for all $\lambda_0\in I$, $\ds \lim_{\lambda\to \lambda_0}\frac{\pa g}{\pa \lambda}(\lambda,r)=\frac{\pa g}{\pa \lambda}(\lambda_0,r)$ for a.e. $r$;
\item[(ii)] for all $\lambda \in I$, $\ds\left|\frac{\pa g}{\pa \lambda}(\lambda,r)\right|\leq q_\lambda(r)$ a.e., where $q_\lambda\in L^1(\RR_+) $, and for all $\lambda_0\in I$, $q_\lambda \to q_{\lambda_0}$ in $L^1(\RR_+)$ as $\lambda \to \lambda_0$.
\end{itemize}
Then, $G$ is $\calC^1$ on $I$ and
$$G'(\lambda)=\int_{\RR_+} \frac{\pa g}{\pa \lambda}(\lambda,r)dr.$$
\end{lemma}
\begin{proof}
Let $\lambda_0, \lambda\in I$. Since $g\in \calC^0(I,L^1(\RR_+))$ and $\frac{\pa g}{\pa \lambda}\in L^1(I\times \RR_+)$, we have
$$g(\lambda,r)-g(\lambda_0,r)=\int_{\lambda_0}^\lambda \frac{\pa g}{\pa \lambda}(\mu,r)d\mu.$$
Hence, by Fubini,
\be
\label{A1}
\frac{G(\lambda)-G(\lambda_0)}{\lambda-\lambda_0}=\int_{\RR_+}\frac{g(\lambda,r)-g(\lambda_0,r)}{\lambda-\lambda_0}dr= \frac{1}{\lambda-\lambda_0}\int_{\lambda_0}^\lambda \left(\int_{\RR_+}\frac{\pa g}{\pa \lambda}(\mu,r)dr\right)d\mu.
\ee
Now, we use a generalized version of the dominated convergence as stated in \cite{lieb-loss} (see Remark after Theorem 1.8) and deduce from Assumptions (i) and (ii) that
\be
\label{clwm}
\lim_{\mu\to \lambda_0}\int_{\RR_+}\frac{\pa g}{\pa \lambda}(\mu,r)dr =\int_{\RR_+}\frac{\pa g}{\pa \lambda}(\lambda_0,r)dr.
\ee
Hence, by using \fref{clwm}, we pass to the limit in \fref{A1} and obtain
$$
\lim_{\lambda\to 0}\frac{G(\lambda)-G(\lambda_0)}{\lambda-\lambda_0}=\int_{\RR_+}\frac{\pa g}{\pa \lambda}(\lambda_0,r)dr,
$$
which proves the differentiability of $G$. In fact, we observe that the same Assumptions (i) and (ii) associated to the same generalized version of the dominated convergence theorem provide the continuity of $G'$. This ends the proof of the lemma.
\end{proof}


\section{Proof of the Antonov inequality \fref{coercivityproeprty}}
\label{antonov}


Let $\Omega$ be  defined by \fref{omega} and let $\xi\in  \calC^{\infty}_c(\Omega)\cap L^{2,odd}_{|F'_e|}$. We recall that the linear transport operator $\calT$ is defined by
$$\calT \xi = v \cdot \nabla_{x}\xi - \nabla_{x}\phi_{Q}\cdot\nabla_{v}\xi.$$ 
Our aim is to prove the coercivity property \fref{coercivityproeprty}. Let $g=\calT \xi$, we have from the Poisson equation

\bee
r^2\phi_{g}'(r)= \int_{0}^r s^2\rho_{g}(s)ds &= &\frac{1}{4\pi}\int_{|x|\leq r}\rho_{g}(x) dx\\
&= &\frac{1}{4\pi} \int_{|x|\leq r} \int_{\RR^3}\left(v \cdot \nabla_{x}\xi - \nabla_{x}\phi_{Q}\cdot\nabla_{v}\xi\right) dv dx\\
&= &\frac{1}{4\pi}  \int_{|x|\leq r} \nabla_{x} \cdot \left( \int_{\RR^3} v\xi dv\right) dx\\
&= &\frac{1}{4\pi}  \int_{|x|=r} \left( \int_{\RR^3} \frac{x}{r}\cdot v \xi dv\right) d\sigma(x).
\eee
Now we observe from the spherical symmetry of $\xi$ that the quantity $\int_{\RR^3} \frac{x}{r}\cdot v \xi dv$ only depends on
$r=|x|$. Therefore
\be
r\phi_{g}'(r)=  \int_{\RR^3} (x\cdot v) \xi dv.
\ee
We then use the Cauchy-Schwarz inequality and $\mbox{Supp($\xi$)}\subset \Omega$ to estimate:
\be
\label{antCS}
r^2\phi_{g}'(r)^2\leq  \left(\int_{\RR^3} (x\cdot v)^2 |F'_{e}| dv\right) \left(\int_{\RR^3} \frac{\xi^2} {|F'_{e}|} dv \right),
\ee
where we recall that $F'_{e}$ is given by \fref{deffphiprime}.
Now we claim that
\be
\label{antid1} 
\int_{\RR^3} (x\cdot v)^2 |F'_{e}| dv = r^2 \rho_{Q}(r).
\ee
Indeed, we first pass to spherical coordinates in $v$
\be
\label{ch1}u=|v|, \ \ \ \   |x\times v|^2= r^2u^2 \sin^2 \theta, \ \ \ \theta \in [0,\pi[,  \ \ \mbox{with}\ \ r=|x|,
\ee
and  get (recall that $F'_{e} \leq 0$)
$$ \int_{\RR^3} (x\cdot v)^2 |F'_{e}| dv =  - 4\pi\int_{0}^{\pi/2} \int_{r=0}^{+\infty}  r^2u^4 \cos^2\theta \sin\theta \frac{\pa F}{\pa e}\left( \frac{u^2}{2} +\phi_Q(r), r^2u^2\sin^2\theta\right) du d\theta$$
Now we perform the change of variable 
\be
\label{ch2} (u,\theta) \rightarrow  \left (e= \frac{u^2}{2} +\phi_Q(r), \ell = r^2u^2\sin^2\theta\right)
\ee
and obtain
$$
 \int_{\RR^3} (x\cdot v)^2 |F'_{e}| dv = -2\pi\sqrt{2}\int_{\ell=0}^{+\infty} \int_{e=\phi_Q(r)+\frac{\ell}{2r^2}}^{+\infty} 
\left( e- \phi_Q(r) - \frac{\ell}{2r^2}\right)^{1/2} \frac{\pa F}{\pa e}(e,\ell) de d\ell.$$
We then integrate by parts with respect to the variable $e$:
$$ \int_{\RR^3} (x\cdot v)^2 |F'_{e}| dv = \pi\sqrt{2}\int_{\ell=0}^{+\infty} \int_{e=\phi_Q(r)+\frac{\ell}{2r^2}}^{+\infty} 
\left( e- \phi_Q(r) - \frac{\ell}{2r^2}\right)^{-1/2} F(e,\ell) \ell de d\ell.$$
Using the same changes of variables \fref{ch1} and \fref{ch2}, we get

\bee r^2 \rho_{Q}(x) &=&r^2\int_{\RR^3}  F\left( \frac{|v|^2}{2} +\phi_Q(x), |x\times v|^2\right)  dv\\
&= &\pi\sqrt{2}\int_{\ell=0}^{+\infty} \int_{e=\phi_Q(r)+\frac{\ell}{2r^2}}^{+\infty} 
\left( e- \phi_Q(r) - \frac{\ell}{2r^2}\right)^{-1/2} F(e,\ell) \ell de d\ell,
\eee
and \fref{antid1} follows.

Now, we integrate the inequality \fref{antCS} with respect to $r$ and use \fref{antid1} to get
\bee
\int_{\RR^3}|\nabla_{x}\phi_{g}|^2 dx &=  &4\pi \int_{0}^{+\infty} r^2\phi_{g}'(r)^2dr\\
&\leq& \int_{\RR^3}\int_{0}^{+\infty} 4\pi r^2 \rho_{Q}(r) \frac{\xi^2}{|F'_{e}|} dr dv= \int_{\RR^3\times \RR^3} \rho_{Q}(x) \frac{\xi^2}{|F'_{e}|} dx dv
\eee
From the definition \fref{antonovpositive} of the Antonov functional, we then deduce 
\be
\label{antonovfirst}
 \mathscr A(\calT\xi,\calT\xi)\geq \int_{\Omega} \left( (\calT \xi)^2 - \rho_{Q}(x)\xi^2\right)\frac{1}{|F'_e|}\,dxdv.
 \ee
Let now $\xi = (x\cdot v) q(x,v)$ and write from the definition of $\calT$
\bee (\calT \xi)^2 &=& \left(q\calT (x\cdot v) + (x\cdot v) \calT q\right)^2 \\
&=& (x\cdot v)^2(\calT q)^2 + (x\cdot v) \calT(x\cdot v)  \calT (q^2) + q^2 (\calT(x\cdot v))^2\\
&=& (x\cdot v)^2(\calT q)^2 + \calT \left( (x\cdot v) q^2 \calT(x\cdot v)\right) - (x\cdot v)q^2\calT \left(\calT(x\cdot v)\right) .
\eee
We observe, from the Poisson equation $\Delta \phi_{Q}=\rho_{Q}$, that
$$ \calT \left(\calT(x\cdot v)\right) =  -(x\cdot v) \Delta \phi_{Q} - v\cdot \nabla_{x}\phi_{Q} = -(x\cdot v)\left( \rho_{Q}(r)+ \frac{\phi_{Q}'(r)}{r}\right).$$
Thus
\be
(\calT \xi)^2 - \rho_{Q}(r) \xi^2 = (x\cdot v)^2(\calT q)^2 +  \calT \left( (x\cdot v) q^2 \calT(x\cdot v)\right)  + \xi^2\frac{\phi_{Q}'(r)}{r}.
\ee
We now insert this expression into \fref{antonovfirst} and directly get the desired Antonov's inequality \fref{coercivityproeprty}, provided  the following claim is proved:
\be
\label{intzero} \int_{\scriptsize\begin{array}{@{}l}\\ |x\cdot v|>\eps^3 \\ |x|>\eps\end{array}} \calT \left( (x\cdot v) q^2 \calT(x\cdot v)\right) dx dv \to 0\quad \mbox{as }\eps \to 0.
\ee
{\em Proof of \fref{intzero}:}  We shall in fact deal with the singularity at $x\cdot v =0$ in the integral \fref{intzero}, recalling that $q(x,v)= \xi(x,v)/(x\cdot v).$ We observe that the function $|x|q(x,v)$
is bounded. To see this, let $x\ne 0$ and $R_x$ be the orthogonal transformation of $\RR^3$ such that $R_x \frac{x}{|x|}=e_1$, where $e_1=(1,0,0)^T$. Then, due to the spherical symmetry,
$$|x|q(x,v)=\frac{\xi(|x|e_1,R_xv)}{e_1\cdot R_xv}=\omega(|x|,R_xv),\quad \mbox{with }\omega(r,v)=\frac{\xi(r e_1,v)}{e_1\cdot v}=\frac{\tilde{\xi}(r,|v|,e_1\cdot v)}{e_1\cdot v}.$$
Now, we recall that $\xi$ is odd in $v$ and hence $\tilde{\xi}$ is odd with respect to the last coordinate, thus $\omega$ is bounded and so is $rq$. Note also that $q$ is smooth on $|x\cdot v|>\delta$, for all $\delta >0$.

Let $\eps>0$. We have
\bee
&&\int_{\scriptsize\begin{array}{@{}l}\\ |x\cdot v|>\eps^3 \\ |x|>\eps\end{array}} \calT \left( (x\cdot v) q^2 \calT(x\cdot v)\right) dx dv=\\
&&\quad=-2\eps^3\int_{\scriptsize\begin{array}{@{}l}\\ x\cdot v=\eps^3 \\ |x|>\eps\end{array}}  q^2 \left(\calT(x\cdot v)\right)^2 \frac{d\sigma_{1}(x,v)}{\sqrt{|x|^2+|v|^2}}-\frac{2}{\eps}\int_{\scriptsize\begin{array}{@{}l}\\ x\cdot v>\eps^3 \\ |x|=\eps\end{array}} (x\cdot v)^2 q^2 \calT(x\cdot v)d\sigma_{2}(x) dv,
\eee
where $d\sigma_{1}(x,v)$ is the measure on the set $\{(x,v) \ \mbox{s.t.}\   x\cdot v=\eps^3 \  \mbox{and}\  |x|>\eps\}$ induced by the Lebesgue measure of $\RR^6$, and $d\sigma_{2}(x)$ is the usual measure on the sphere $\{x\in \RR^3;  |x|=\eps\}$. 
Now let $R>0$ such that $Supp(\xi) \subset \{(x,v), \ |x|^2 +|v|^2 \leq R^2\}$, then
\bee
&&\left|\int_{\scriptsize\begin{array}{@{}l}\\ |x\cdot v|>\eps^3 \\ |x|>\eps\end{array}} \calT \left( (x\cdot v) q^2 \calT(x\cdot v)\right) dx dv\right|\leq\\
&&\quad\leq 2\eps\int_{\scriptsize\begin{array}{@{}l}\\ x\cdot v=\eps^3 \\ |x|>\eps\end{array}}  r^2q^2 \left(\calT(x\cdot v)\right)^2 \frac{d\sigma_{1}(x,v)}{\sqrt{|x|^2+|v|^2}}+\frac{2}{\eps}\int_{|x|=\eps} \left(\int_{\RR^3}\xi^2 |\calT(x\cdot v)|^2 dv\right)d\sigma_{2}(x) ,\\
&&\quad \leq 2\eps |rq|_{L^\infty}^2I(\eps,R) + \frac{C}{\eps} \eps^2,
\eee
where we have set
$$
I(\eps,R)=\int_{\scriptsize\begin{array}{@{}l}\\ x\cdot v=\eps^3 \\ |x|^2+|v|^2 < R^2\end{array}}   \left(\calT(x\cdot v)\right)^2 \frac{d\sigma_{1}(x,v)}{\sqrt{|x|^2+|v|^2}}
$$
and where we have used in the last estimate that $rq$  is bounded and that $\xi$ is compactly supported. We claim that 
\be
\label{clfm}
I(\eps,R)\leq C_R,
\ee
where $C_R$ is independent of $\eps$, which concludes the proof of \fref{intzero}. Indeed, we integrate by parts to get:
\bee
I(\eps,R)&=&-\int_{\scriptsize\begin{array}{@{}l}\\ x\cdot v>\eps^3 \\ |x|^2+|v|^2<R^2\end{array}} \calT \left( \calT(x\cdot v)\right) dx dv\\
&&+\int_{\scriptsize\begin{array}{@{}l}\\ x\cdot v>\eps^3 \\ |x|^2+|v|^2 = R^2\end{array}}  \calT(x\cdot v) \left(x\cdot v-v\cdot \na_x\phi_Q\right)\frac{d\sigma_{3}(x,v)}{\sqrt{|x|^2+|v|^2}}\\
&\leq& \int_{|x|^2+|v|^2<R^2} \left|\calT \left( \calT(x\cdot v)\right)\right| dx dv\\
&&+\int_{|x|^2+|v|^2 = R^2}  \left|\calT(x\cdot v) \left(x\cdot v-v\cdot \na_x\phi_Q\right)\right|\frac{d\sigma_{3}(x,v)}{\sqrt{|x|^2+|v|^2}}\leq C_R,
\eee
where $d\sigma_{3}(x,v)$ denotes the usual measure on the sphere  $ |x|^2+|v|^2 = R^2$. This concludes the proof of \fref{clfm} and \fref{coercivityproeprty}.


\section{Proof of Lemma  \ref{lemmaprojection}}
\label{appendixlemma}

Let us first prove that for $h \in \Phi_{rad}$, the projection $\mathcal Ph$ given by \fref{pikernel} is well-defined.  From Lemma \ref{lemmapotentials}, the denominator in the definition \fref{pikernel} is finite and non zero for $(x,v)$ such that $(e(x,v),\ell(x,v))\in \mathcal D$ and we have
$$|({\mathcal P} h)(x,v)|\leq \sup_{r_1\leq r\leq r_2}|h(r)|.$$ Let us now prove that $hF'_e$ belongs to $L^{2,r}_{|F'_e|}$. By a change of variable, we have
$$\begin{array}{l}
\ds \int h^2\,|F'_e|\,dxdv=\\[3mm]
\ds \quad =-4\pi^2\sqrt{2}\int_{L}\int_0^{+\infty} \frac{\pa Q^*}{\pa s}(s,\ell)\int_{r_1}^{r_2}\left(a_{\phi_Q}^{-1}(s,\ell)- \phi_Q(r)-\frac{\ell}{2r^2}\right)^{-1/2}(h(r))^2drdsd\ell\\[5mm]
\ds \quad \leq
-C \int_{L}\int_0^{s_0(\ell)}\frac{\pa Q^*}{\pa s}(s,\ell)\int_{r_1}^{r_2}
\frac{r|h(r)|^2\sqrt{r_1r_2}}{\sqrt{\ell(r-r_1)(r_2-r)}}drdsd\ell,
\end{array}
$$
where we used \fref{ss1}, \fref{ss2} for the support of $Q^*$ (recall the definition \fref{L} of $L$) and \fref{clvm} for $\phi_Q$. From \fref{boundL}, \fref{ss3} and the radial Sobolev bound:
$$r^{1/2}|h(r)|\leq \left(\int_r^{+\infty}s^2(h'(s))^2\,ds\right)^{1/2}\leq |\na h|_{L^2},$$
we deduce that
\bea
\label{houmtsouk}
 \nonumber \int h^2\,|F'_e|\,dxdv & \leq & 
-C |\na h|_{L^2}^2\int_{L}\int_0^{s_0(\ell)}\frac{\pa Q^*}{\pa s}(s,\ell)\int_{r_1}^{r_2}
\frac{1}{\sqrt{\ell(r-r_1)(r_2-r)}}drdsd\ell\\
& \leq & 
C |\na h|_{L^2}^2\int_0^{\ell_0}Q^*(0,\ell)\frac{d\ell}{\sqrt{\ell}}\leq
C |\na h|_{L^2}^2.
\eea
Now we prove that $(\mathcal Ph)F'_e$ belongs to $L^{2,r}_{|F'_e|}$. The same change of variable as above gives
$$\begin{array}{l}
\ds \int (\mathcal Ph)^2\,|F'_e|\,dxdv=\\
\ds =-4\pi^2\sqrt{2}\int_{L}\int_0^{s_0(\ell)}\frac{\pa Q^*}{\pa s}(s,\ell)\frac{\ds \left(\int_{r_1}^{r_2} \left(a_{\phi_Q}^{-1}(s,\ell)- \phi_Q(r)-\frac{\ell}{2r^2}\right)^{-1/2}h(r)dr\right)^2}{\ds \int_{r_1}^{r_2} \left(a_{\phi_Q}^{-1}(s,\ell)- \phi_Q(r)-\frac{\ell}{2r^2}\right)^{-1/2}dr}dsd\ell.
\end{array}
$$
Hence, using the Cauchy-Schwarz inequality, we get
$$
\int (\mathcal Ph)^2\,|F'_e|\,dxdv\leq \int h^2\,|F'_e|\,dxdv
$$
which concludes the proof of \fref{estun}.

Observe now from \fref{pikernel}, \fref{estun} that $(\mathcal Ph)F'_e$ is a $L^{2,r}_{|F'_e|}$ function of $e(x,v)$ and $\ell(x,v)$ and hence belongs to $N(\mathcal T)$ from \fref{kernelt}. It remains to prove that $(h-\mathcal Ph)F'_e$ is orthogonal to $N(\mathcal T)$. Indeed, let
$$\theta=\theta(e(x,v),\ell(x,v))\in N(\mathcal T),$$
then passing from the variables $x,v$ to the variable $r,e,\ell$, we get:
\bee
&&\int (\mathcal Ph)(x,v)\theta(e(x,v),\ell(x,v))dxdv=\\
&&\qquad \qquad =2\pi^2\sqrt{2}\int_0^{+\infty}\int_{e_{\phi_Q,\ell}}^0\int_{r_1}^{r_2}\left(e-\phi_Q(r)-\frac{\ell}{2r^2}\right)^{-1/2}h(r)\theta(e,\ell)drded\ell\\
&&\qquad \qquad =\int h(x,v)\theta(e(x,v),\ell(x,v))dxdv.
\eee
Hence $(h-\mathcal Ph)F'_e$ and $\theta$ are orthogonal for the $L^{2,r}_{|F_e'|}$ scalar product and \fref{projectionkernel} follows.\\ 
This ends the proof of Lemma \ref{lemmaprojection}.


\end{appendix}


\begin{thebibliography}{9}

\bibitem{ATL} Alvino, A.; Trombetti, G.; Lions, P.-L., On optimization problems with prescribed rearrangements.  Nonlinear Anal.  13  (1989),  no. 2, 185--220.

\bibitem{aly} Aly J.-J., On the lowest energy state of a collisionless self-gravitating system under phase volume constraints. MNRAS {\bf 241} (1989), 15.


\bibitem{A1} Antonov, A. V., Remarks on the problem of stability in stellar dynamics. {\em Soviet Astr., AJ.,} {\bf 4},
859-867 (1961).

\bibitem{A2} Antonov, A. V., Solution of the problem of stability of a stellar system with the Emden density law and spherical velocity distribution. {\em J. Leningrad Univ. Se. Mekh. Astro. } {\bf 7}, 135-146 (1962).

\bibitem{arnold1} Arnold, V.I., Conditions for nonlinear stability of stationary [;ane curvilinear fows of an ideal fluid, Sov. Math. Dokl. 6, 773-776 (1965)

\bibitem{arnold2} Arnold, V.I., Sur un principe variationel pour les ecoulements stationaires des liquides parfaits et ses applications aux probl\`emes de stabilit\'e nonlin\'enaire, J. M\'ecanique 5, 29-43 (1966).

\bibitem{arnoldbook} Arnold, V.I., Mathematical models of classical mechanics, New York, Springer Verlag, 1980.

%
\bibitem{Batt} Batt, J.; Faltenbacher, W.; Horst, E., Stationary spherically symmetric models in stellar dynamics, Arch. Rat. Mech. Anal. 93, 159-183 (1986).

\bibitem{Brezis} Br\'ezis, H., Analyse fonctionnelle. Th\'eorie et applications. Collection Math\'ematiques Appliqu\'ees pour la Ma\^{i}trise. Masson, Paris, 1983.
%
\bibitem{binney} Binney, J.; Tremaine, S., Galactic Dynamics, Princeton University Press, 1987.

\bibitem{CL} Cazenave, T.; Lions, P.-L., Orbital stability of standing waves for some nonlinear Schr\"odinger equations, Comm. Math. Phys. 85 (1982), no. 4, 549-561.

\bibitem{chavanis} Chavanis, P.-H., Dynamical stability of collisionless stellar systems and barotropic stars: the nonlinear Antonov first law, Astronomy and Astrophysics {\bf 451} (2006), no. 1, 109-123.

\bibitem{Dol} Dolbeault, J.; S\'anchez, \'O.; Soler, J., Asymptotic behaviour for the Vlasov-Poisson system in the stellar-dynamics case, Arch. Ration. Mech. Anal. 171 (2004), no. 3, 301--327.

\bibitem{doremus} Doremus, J. P.; Baumann, G.; Feix, M. R., Stability of a Self Gravitating System with Phase Space Density Function of Energy and Angular Momentum, Astronomy and Astrophysics {\bf 29} (1973), 401.



%
\bibitem{fridman} Fridmann, A. M.; Polyachenko, V. L., Physics of gravitating systems, Springer-Verlag, 1984.

\bibitem{gardner} Gardner, C.S., Bound on the energy available from a plasma, Phys. Fluids 6, 1963, 839-840.

\bibitem{gillon} Gillon, D.; Cantus, M.; Doremus, J. P.; Baumann, G., Stability of self-gravitating spherical systems in which phase space density is a function of energy and angular momentum, for spherical perturbations, Astronomy and Astrophysics {\bf 50} (1976), no. 3, 467--470.

\bibitem{Guo} Guo, Y., Variational method for stable polytropic galaxies, Arch. Rat. Mech. Anal. {\bf 130} (1999), 163-182. 

\bibitem{GuoLin} Guo, Y.; Lin, Z., Unstable and stable galaxy models,  Comm. Math. Phys.  {\bf 279}  (2008),  no. 3, 789--813.

\bibitem{GuoRein2} Guo, Y.; Rein, G., Stable steady states in stellar dynamics, Arch. Rat. Mech. Anal. {\bf 147} (1999), 225--243.

\bibitem{GuoRein} Guo, Y.; Rein, G., Isotropic steady states in galactic dynamics, Comm. Math. Phys. {\bf 219} (2001), 607--629.

\bibitem{Guo2} Guo, Y., On the generalized Antonov's stability criterion. {\em Contemp. Math.}  {\bf 263}, 85-107 (2000)


\bibitem{GR4} Guo, Y.; Rein, G., A non-variational approach to nonlinear Stability in stellar dynamics applied to the King model, Comm. Math. Phys., 271, 489-509 (2007).

\bibitem{HR} Had\'zi\'c, Mahir; Rein, Gerhard Global existence and nonlinear stability for the relativistic Vlasov-Poisson system in the gravitational case, Indiana Univ. Math. J. 56 (2007), no. 5, 2453--2488.

%

\bibitem{kandrup1} Kandrup, H. E.; Sygnet, J. F., A simple proof of dynamical stability for a class of spherical clusters.  Astrophys. J.  298  (1985),  no. 1, part 1, 27--33.

\bibitem{kandrup2} Kandrup, H. E., A stability criterion for any collisionless stellar equilibrium and some concrete applications thereof, Astrophys. J. 370 (1991) no. 1, 312--317.

\bibitem{kavian} Kavian, O., Introduction \`a la th\'eorie des points critiques et applications aux probl\`emes elliptiques.
Math\'ematiques \& Applications (Berlin), 13. Springer-Verlag, Paris, 1993. 
%

\bibitem{LMR-note} Lemou, M.; M\'ehats, F.; Rapha\"el, P.,  Orbital stability and singularity formation for Vlasov-Poisson systems.  C. R. Math. Acad. Sci. Paris  341  (2005),  no. 4, 269--274. 

%
\bibitem{LMR1} Lemou, M.; M\'ehats, F.; Rapha\"el, P., On the orbital stability of the ground states and the singularity formation for the gravitational Vlasov-Poisson system, Arch. Rat. Mech. Anal. 189 (2008), no. 3, 425--468.

\bibitem{LMR5}  Lemou, M.; M\'ehats, F,; Rapha\"el, P., Stable ground states for the relativistic gravitational Vlasov-Poisson system,  Comm. Partial Diff. Eq. 34 (2009), no. 7, 703--721.
%
\bibitem{LMR6} Lemou, M.; M\'ehats, F,; Rapha\"el, P., Ensemble inequivalence for the gravitational Vlasov-Poisson system, in preparation.

\bibitem{LMR7} Lemou, M.; M\'ehats, F.; Rapha\"el, P., Stable self-similar blow up dynamics for the three dimensional relativistic gravitational Vlasov-Poisson system,  J. Amer. Math. Soc.  {\bf 21}  (2008),  no. 4, 1019--1063. 

\bibitem{lieb-loss} Lieb, E. H.; Loss, Analysis. Second edition. Graduate Studies in Mathematics, 14. American Mathematical Society, Providence, RI, 2001.
%
\bibitem{LiebYau} Lieb, E.H.; Yau, H.T., The Chandrasekhar theory of stellar collapse as the limit of quantum mechanics, Comm. Math. Phys. 112 (1987), no. 1, 147--174.

\bibitem{lin} Lin, Z., Nonlinear instability of periodic BGK waves for Vlasov-Poisson system.  Comm. Pure Appl. Math.  58  (2005),  no. 4, 505--528.

\bibitem{lin-strauss1} Lin, Z.; Strauss, W. A., Linear stability and instability of relativistic Vlasov-Maxwell systems.  Comm. Pure Appl. Math.  60  (2007),  no. 5, 724--787.

\bibitem{lin-strauss2} Lin, Z.; Strauss, W. A., A sharp stability criterion for the Vlasov-Maxwell system.  Invent. Math.  173  (2008),  no. 3, 497--546.

%
\bibitem{PLL1} Lions, P.-L., The concentration-compactness principle in the calculus of variations. The locally compact case. I, Ann. Inst. H. Poincar\'e Anal. Non Lin\'eaire 1 (1984), no. 2, 109--145. 

%
\bibitem{PLL2} Lions, P.-L., The concentration-compactness principle in the calculus of variations. The locally compact case. II, Ann. Inst. H. Poincar\'e Anal. Non Lin\'eaire 1 (1984), no. 4, 223--283.

\bibitem{LP} Lions, P.-L.; Perthame, B., Propagation of moments and regularity for the $3$-dimensional Vlasov-Poisson system,  Invent. Math. 105 (1991), no. 2.

\bibitem{lynden} Lynden-Bell, D., The Hartree-Fock exchange operator and the stability of galaxies, Mon. Not. R. Astr. Soc. {\bf 144}, 1969, 189--217.

\bibitem{lynden2} Lynden-Bell, D., Lectures on stellar dynamics. Galactic dynamics and N-body simulations, 3--31, Lecture Notes in Phys. 433, Springer, Berlin, 1994.

\bibitem{MPbook}  Marchioro, C.; Pulvirenti, M., Mathematical theory of incompressible nonviscous fluids. Applied Mathematical Sciences, 96. Springer-Verlag, New York, 1994.

\bibitem{Mossino} Mossino, J., In\'egalit\'es isop\'erim\'etriques et applications en physique. (French) [Isoperimetric inequalities and applications to physics] Travaux en Cours. [Works in Progress] Hermann, Paris, 1984.

\bibitem{MP} Marchioro, C.; Pulvirenti, M., Some considerations on the nonlinear stability of stationary planar Euler flows, Comm. Math. Phys. 100 (1985), no. 3, 343--354.

%
\bibitem{PA} Perez, J.; Aly, J.-J., Stability of spherical stellar systems -I. Analytical results. Monthly. Not. Royal. Astronomical Soc., {\bf 280}, 689-699 (1996).

%
\bibitem{Pf} Pfaffelmoser, K., Global classical solutions of the Vlasov-Poisson system in three dimensions for general initial data, J. Diff. Eq. {\bf 95} (1992), 281-303.

%
\bibitem{SS} S\'anchez, \'O.; Soler, J., Orbital stability for polytropic galaxies,  Ann. Inst. H. Poincar\'e Anal. Non Lin\'eaire  23  (2006),  no. 6, 781--802.

%
\bibitem{S1} Schaeffer, J., Global existence of smooth solutions to the Vlasov-Poisson system in three dimensions, Comm. Part. Diff. Eq. {\bf 16} (1991), 1313-1335.

%
\bibitem{S2} Schaeffer, J., Steady States in Galactic Dynamics, Arch. Rational, Mech. Anal.  172 (2004), 1--19.

\bibitem{serre} Serre, D., Sur le principe variationnel des \'equations de la m\'ecanique des fluides parfaits. [On the variational principle for the equations of perfect fluid dynamics] RAIRO Mod\'el. Math. Anal. Num\'er. {\bf 27} (1993), no. 6, 739--758. 

\bibitem{SDLP} Sygnet, J.-F.; Des Forets, G.; Lachieze-Rey, M.; Pellat, R., Stability of gravitational systems and gravothermal catastrophe in astrophysics, Astrophys. J. 276 (1984), no. 2, 737--745.

\bibitem{W} Wan, Y-H, On nonlinear stability of isotropic models in stellar dynamics, Arch. Ration. Mech. Anal 147 (1999), no. 3, 245-268.

\bibitem{WZS} Wiechen, H., Ziegler, H.J., Schindler, K. Relaxation of collisionless self gravitating matter: the lowest energy state, Mon. Mot. R. ast. Soc (1988) 223, 623-646.

\bibitem{wol}  Wolansky, G., On nonlinear stability of polytropic galaxies. {\em Ann. Inst. Henri Poincar\'e,} 16, 15-48 (1999).

\bibitem{WG} Wolansky, G., Ghil, M., Nonlinear Stability  for Saddle Solutions of Ideal Flows and Symmetry Breaking.  {\em Commun. Math. Phys} 193, 713-736 (1998).

\end{thebibliography}
\end{document}